\documentstyle[12pt,amssymb,graphicx,amsthm,tabularx]{article}
\def\IC{\hbox{\rm C \kern-.80em \vrule depth 0ex height 1.5ex width .05em \kern.41em}}
\def\IR{\hbox{\rm R \kern-.80em \vrule depth 0ex height 1.5ex width .05em \kern.41em}}
\def\IZ{\hbox{$\,$\parbox[bn]{2.6mm}{\resizebox*{2.0mm}{1.18ex}{/}} \kern-1.14em \rm Z \kern-.35em}}

\def\fn#1{\mathop{{\rm #1}\vphantom{\dim}}}
\def\abb#1#2{\parbox[c][2mm][c]{2mm}{$\put(#2){#1}$}} 

\makeatletter
\renewcommand{\section}{\@startsection{section}{1}{0mm}{15mm}{5mm}{\bf\raggedright}}
\makeatother

\newtheorem{theorem}{\indent Theorem}[section]

\newtheorem{coro}[theorem]{\indent Corollary}
\newtheorem{remark}[theorem]{\indent Remark}
\newtheorem{lemma}[theorem]{\indent Lemma}
\newtheorem{fact}[theorem]{\indent Fact}
\newtheorem{example}{\indent Example}[section]

\newtheorem{thm}{\indent Theorem}

\newtheorem{claim}{\indent Claim}

\begin{document}

\begin{center}
{\large\bf MAXIMAL EXPONENTS OF K-PRIMITIVE MATRICES: THE POLYHEDRAL CONE CASE}\vspace{8mm}\\
{\bf Raphael Loewy}\vspace{1mm}\\
Department of Mathematics\\
Technion\\
Haifa 32000, Israel\vspace{4mm}\\
and\vspace{4mm}\\
{\bf Bit-Shun Tam$^{*\bf 1}$}\vspace{1mm}\\
Department of Mathematics\\
Tamkang University\\
Tamsui, Taiwan 251, R.O.C.
\vspace{2mm}\\

February 18, 2009
\end{center}

{\bf Abstract}. Let $K$ be a proper (i.e., closed, pointed, full
convex) cone in ${\Bbb R}^n$.  An $n\times n$ matrix
 $A$ is said to be $K$-primitive if there exists a positive integer $k$ such that $A^k(K \setminus \{ 0 \})
\subseteq$  int\,$K$; the least such $k$ is referred to as the
exponent of $A$ and is denoted by $\gamma(A)$.  For a polyhedral
cone $K$, the maximum value of $\gamma(A)$, taken over all
$K$-primitive matrices $A$, is denoted by $\gamma(K)$.  It is
proved that for any positive integers $m,n, 3 \le n \le m$, the
maximum value of $\gamma(K)$, as $K$ runs through all
$n$-dimensional polyhedral cones with $m$ extreme rays, equals
$(n-1)(m-1)+1$ when $m$ is even or $m$ and $n$ are both odd, and
is at least $(n-1)(m-1)$ and at most $(n-1)(m-1)+1$ when $m$ is
odd and $n$ is even. For the cases when $m = n, m = n+1$ or $n =
3$, the cones $K$ and the corresponding $K$-primitive matrices $A$
such that $\gamma(K)$ and $\gamma(A)$ attain the maximum value are
identified up to respectively linear isomorphism and
cone-equivalence modulo positive scalar multiplication.

\renewcommand{\thefootnote}{}
\footnote{\vspace*{-6mm}\leftmargini=2mm
\begin{enumerate} \itemsep=-3pt
\item[] {\it AMS classification}: 15A48; 05C50; 47A06. \item[]
\item[{$^*\!\!$}] Corresponding author. \item[] { \it E-mail
addresses}\,: loewy@techunix.technion.ac.il (R. Loewy),
bsm01@mail.tku.edu.tw (B.-S.\ Tam). \item[{$^1\!\!$}] Supported by
National Science Council of the Republic of China.

\end{enumerate}}

\newpage
\section{Introduction}

If $K$ is a polyhedral (proper) cone in ${\Bbb R}^n$ with $m$
extreme rays, what is the maximum value of the exponents of
$K$-primitive matrices ?  This question was posed by Steve
Kirkland in an open problem session at the $8$th ILAS conference
held in Barcelona in July, 1999. Here by a {\it $K$-primitive
matrix} we mean a real square matrix $A$ for which there exists a
positive integer $k$ such that $A^k$ maps every nonzero vector of
$K$ into the interior of $K$; the least such $k$ is referred to as
the {\it exponent of $A$} and is denoted by $\gamma(A)$. In view
of Wielandt's classical sharp bound for exponents of (nonnegative)
primitive matrices of a given order, Kirkland conjectured that
$m^2-2m+2$ is an upper bound for the maximum value considered in
his question. This work is an outcome of our attempt to answer
Kirkland's question.

In the classical nonnegative matrix case, the determination of
upper bounds for the exponents of primitive matrices under various
assumptions has been treated mainly by a graph-theoretic approach.
Here for a $K$-primitive matrix $A$, we work with the digraph
$({\cal E},{\cal P}(A,K))$, which is one of the four digraphs
associated with $A$, introduced by Barker and Tam \cite{B--T},
\cite{T--B}. (Formal definitions will be given later.) Based on
the same digraph, Niu \cite{Niu} has started an initial study of
the exponents of $K$-primitive matrices over a polyhedral cone
$K$. His work has motivated partly the work of Tam \cite{Tam 4}
and our present work.

The study of $K$-primitive matrices in the general polyhedral cone
case differs from the nonnegative matrix case (or, equivalently,
the simplicial cone case) in at least two (not unrelated)
respects. First, in the nonnegative matrix case the (distinct)
extreme vectors of the underlying cone are linearly independent,
whereas in the general polyhedral cone case the extreme vectors of
the underlying cone satisfy certain nonzero (linear) relations.
Second, in the nonnegative matrix case it is always possible to
find a nonnegative matrix with a prescribed digraph as its
associated digraph, whereas in the general polyhedral cone case we
often need to treat first the realization problem, that is, to
determine whether there is a polyhedral cone $K$ for which there
is a $K$-nonnegative matrix $A$ such that the digraph $({\cal
E},{\cal P}(A,K))$ is given by a prescribed digraph.  As expected,
and also illustrated by this work, the study of the polyhedral
cone case is more difficult than the classical nonnegative matrix
case.

We now describe the contents of this paper in some detail.

Section $2$ contains most of the definitions, together with the
relevant known results, which we need in the paper.

In Section $3$ we obtain a Sedl\'{a}\u{c}ek-Dulmage-Mendelsohn
type upper bound for the local exponents (see Section 2 for the
definition), and hence also an upper bound for the exponent, of a
$K$-primitive matrix $A$ in terms of the lengths of circuits in
the digraph $({\cal E},{\cal P}(A,K))$ and the degree of the
minimal polynomial of $A$.

In Section $4$ we single out digraphs on $m\,(\ge 4)$ vertices,
with the length of the shortest circuit equal to $m-1$, that may
be realized as $({\cal E},{\cal P}(A,K))$ for some $K$-primitive
matrix $A$, where $K$ is a polyhedral cone with $m$ extreme rays
(see Lemma \ref{lemma3}). It is found that, up to graph
isomorphism, there are two of them, represented by Figure $1$ and
Figure $2$ respectively. (These figures will be given later in the
paper.) It turns out that they are precisely the two known
so-called primitive digraphs on $m$ vertices with the length of
the shortest circuit equal to $m-1$. When $A$ is a $K$-nonnegative
matrix such that the digraph $({\cal E},{\cal P}(A,K))$ is given
by Figure $1$ or Figure $2$, we make some interesting observations
on $A$, and by delicate manipulation with the relations on the
extreme vectors, we also obtain certain geometric properties of
$K$ (see Lemma \ref{lemma4.2}).  As a consequence, it is proved
that if $K$ is an $n$-dimensional polyhedral cone with $m$ extreme
rays then its exponent $\gamma(K)$, which is defined to be $\max
\{ \gamma(A): A \mbox{ is $K$-primitive} \}$, does not exceed
$(n-1)(m-1)+1$. Thus we answer in the affirmative the
above-mentioned conjecture posed by Kirkland.

In Section 5 we prove that the maximum value of $\gamma(K)$ as $K$
runs through all $n$-dimensional minimal cones (i.e., cones having
$n+1$ extreme rays) is $n^2-n+1$ if $n$ is odd, and is $n^2-n$ if
$n$ is even. We also determine (up to linear isomorphism) the
minimal cones $K$ (and also the corresponding $K$-primitive
matrices $A$) such that $\gamma(K)$ (and $\gamma(A)$) attains the
maximum value. In particular, it is found that every minimal cone
$K$ whose exponent attains the maximum value has a balanced
relation for its extreme vectors and also if $A$ is a
$K$-primitive matrix such that $\gamma(A) = \gamma(K)$ then
necessarily the digraph $({\cal E},{\cal P}(A,K))$ is, up to graph
isomorphism, given by
 Figure $1$ or Figure $2$.

Section $6$ is devoted to the $3$-dimensional cone case.  It is
proved that the maximum value of $\gamma(K)$ as $K$ runs through
all $3$-dimensional polyhedral cones with $m$ extreme rays is
$2m-1$, and also that for any $3$-dimensional polyhedral cone $K$
with $m$ extreme rays and any $K$-primitive matrix $A, \gamma(A) =
2m-1$ if and only if the digraph $({\cal E},{\cal P}(A,K))$ is, up
to graph isomorphism, given by Figure $1$. Indeed, for every
positive integer $m \ge 3$, we can construct, for every real
number $\theta\in (\frac{2\pi}{m},\frac{2\pi}{m-1})$, a
$3$-dimensional polyhedral cone $K_{\theta}$ with $m$ extreme rays
and a $K_{\theta}$-primitive matrix $A_{\theta}$ such that the
digraph $({\cal E},{\cal P}(A_{\theta},K_{\theta}))$ is given by
Figure $1$, and for every positive integer $m \ge 5$, we can also
find a $3$-dimensional polyhedral cone $K$ with $m$ extreme rays
for which there does not exist a $K$-primitive matrix $A$ such
that the digraph $({\cal E}, {\cal P}(A,K))$ is given by Figure
$1$ (up to graph isomorphism). Since any two $3$-dimensional
polyhedral cones with the same number of extreme rays are
combinatorially equivalent, this means that the exponents of
combinatorially equivalent cones may not be the same.  The
construction of the pair $(K_{\theta},A_{\theta})$ makes use of
roots of a polynomial of the form $t^m-ct-(1-c)$, where $0 < c <
1$.  In his study of Leslie matrices Kirkland
(\cite{Kir1},\cite{Kir2}) has considered polynomials of a more
general form, namely, those of the form $t^m-\sum_{k=1}^m
a_kt^{m-k}$, where $a_1, \ldots, a_m$ are nonnegative real numbers
with sum equal to $1$. Evidently, there are some connections
between the work of Kirkland in these two papers and our work in
Section $6$.

In Section $7$ we show that for any positive integers $m,n, 3 \le
n \le m$, the maximum value of $\gamma(K)$, as $K$ runs through
all $n$-dimensional polyhedral cones with $m$ extreme rays, equals
$(n-1)(m-1)+1$ when $m$ is even or $m$ and $n$ are both odd, and
is at least $(n-1)(m-1)$ and at most $(n-1)(m-1)+1$ when $m$ is
odd and $n$ is even.  Our proof involves certain generalized
Vandermonde matrices, the complete symmetric polynomials, the
Jacobi-Trudi determinant, and a nontrivial result about
polynomials with nonnegative coefficients.

In Section 8, for the special cases when $m = n, m = n+1$ and $n =
3$, we settle the question of uniqueness of the cones $K$, in the
class of $n$-dimensional polyhedral cones with $m$ extreme rays,
and the corresponding $K$-primitive matrices $A$ whose exponents
attain the maximum value.
 It is proved that for every positive integer $m
\ge 5$, up to linear isomorphism, the $3$-dimensional cones with
$m$ extreme rays that attain the maximum exponent are precisely
the cones $K_{\theta}$'s introduced in Section $6$, uncountably
infinitely many of them; and, for $m \ge 6$, for each $\theta$
there is, up to multiples, only one  $K_{\theta}$-primitive matrix
whose exponent attains the maximum value.  In contrast,
$n$-dimensional minimal cones whose exponents attain the maximum
value are scanty
--- for every integer $n\ge 3$, there are (up to linear
isomorphism) one or two such cones, depending on whether $n$ is
odd or even.  However, for each of such minimal cones, there are
uncountably infinitely many pairwise non-cone-equivalent linearly
independent primitive matrices whose exponents attain the maximum
value.

In Section $9$, the final section, we give an example, some
further remarks and a few open questions.

\setcounter{equation}{0}
\section{Preliminaries}

We take for granted standard properties of nonnegative matrices,
complex matrices and graphs that can be found in textbooks (see,
for instance, \cite{B--P}, \cite{B--R}, \cite{H}, \cite{H--J},
\cite{L--T}). A familiarity with elementary properties of
finite-dimensional convex sets, convex cones and cone-preserving
maps is also assumed (see, for instance,\cite{Bar 2}, \cite{Roc},
\cite{Tam 2}, \cite{Zie}). To fix notation and terminology, we
give some definitions.

Let $K$ be a nonemtpy subset of a finite-dimensional real vector
space $V$.  The set $K$ is called a {\it convex cone} if $\alpha
x+\beta y \in K$ for all $x,y \in K$ and $\alpha, \beta \ge 0$;
$K$ is {\it pointed} if $K \cap (-K) = \{ 0 \}$; $K$ is {\it full}
if its interior $\fn{int}K$ (in the usual topology of $V$) is
nonempty, equivalently, $K-K = V$.  If $K$ is closed and satisfies
all of the above properties, $K$ is called a {\it proper cone}.

{\it In this paper, unless specified otherwise, we always use $K$
to denote a proper cone in the $n$-dimensional Euclidean space
${\Bbb R}^n$.}

We denote by $\ge^K$ the partial ordering of ${\Bbb R}^n$ induced
by $K$, i.e., $x \ge^K y$ if and only if $x-y \in K$.

A subcone $F$ of $K$ is called a {\it face} of $K$ if $x\ge^K
y\ge^K 0$ and $x \in F$ imply $y \in F$. If $S \subseteq K$, we
denote by $\Phi(S)$ the {\it face of $K$ generated by $S$}, that
is, the intersection of all faces of $K$ including $S$.  If $x \in
K$, we write $\Phi(\{ x \})$ simply as $\Phi(x)$.  It is known
that for any vector $x \in K$ and any face $F$ of $K$, $x \in
\fn{ri}F$ if and only if $\Phi(x) = F$; also, $\Phi(x) = \{ y\in
K\!\!: x\ge^K \alpha y \mbox{ for some }\alpha > 0 \}.$ (Here we
denote by $\fn{ri}F$ the {\it relative interior of $F$}.) A vector
$x \in K$ is called an {\it extreme vector} if either $x$ is the
zero vector or $x$ is nonzero and $\Phi(x) = \{ \lambda x\!\!:
\lambda \ge 0 \}$; in the latter case, the face $\Phi(x)$ is
called an {\it extreme ray}. We use $\fn{Ext}K$ to denote the set
of all nonzero extreme vectors of $K$.  Two nonzero extreme
vectors are said to be {\it distinct} if they are not multiples of
each other. The cone $K$ itself and the set $\{ 0 \}$ are always
faces of $K$, known as {\it trivial faces}.  Other faces of $K$
are said to be ${\it nontrivial}$.

If $S$ is a nonempty subset of a vector space, we denote by
$\fn{pos}S$ the {\it positive hull} of $S$, i.e., the set of all
possible nonnegative linear combinations of vectors taken from
$S$.

A closed pointed cone $K$ is said to be the {\it direct sum} of
its subcones $K_1, \ldots, K_p$, and we write $K = K_1 \oplus
\cdots \oplus K_p$ if every vector of $K$ can be expressed
uniquely as $x_1+x_2+ \cdots + x_p$, where $x_i \in K_i, 1 \le i
\le p$.  $K$ is called {\it decomposable} if it is the direct sum
of two nonzero subcones; otherwise, it is said to be {\it
indecomposable}. It is well-known that every closed pointed cone
$K$ can be written as
$$K = K_1 \oplus \cdots \oplus K_p,$$
where each $K_j$ is an indecomposable cone ($1 \le j \le p$).
Except for the ordering of the summands, the above decomposition
is unique.  We will refer to the $K_j$'s as {\it indecomposable
summands} of $K$.

By a {\it polyhedral cone} we mean a proper cone which has
finitely many extreme rays. By the {\it dimension of a proper
cone} we mean the dimension of its linear span.  A polyhedral cone
is said to be {\it simplicial} if the number of extreme rays is
equal to its dimension.  The nonnegative orthant ${\Bbb R}^n_+ :=
\{ (\xi_1, \ldots, \xi_n)^T \in {\Bbb R}^n\!: \xi_i \ge 0\ \forall
i\}$ is a typical example of a simplicial cone.

 We denote by $\pi(K)$ the set of all $n \times n$ real
matrices $A$ (identified with linear mappings on ${\Bbb R}^n$)
such that $AK \subseteq K$. Members of $\pi(K)$ are said to be
{\it K-nonnegative} and are often referred to as {\it
cone-preserving maps}.  It is clear that $\pi({\Bbb R}^n_+)$
consists of all $n \times n$ (entrywise) nonnegative matrices.

A matrix $A\in \pi(K)$ is said to be {\it $K$-irreducible} if $A$
leaves invariant no nontrivial face of $K$;  $A$ is {\it
$K$-positive} if $A(K\setminus \{0\}) \subseteq {\rm int}~K$ and
is {\it $K$-primitive} if there is a positive integer $p$ such
that $A^p$ is $K$-positive.  If $A$ is $K$-primitive, then the
smallest positive integer $p$ for which $A^p$ is $K$-positive is
called the {\it exponent} of $A$ and is denoted by $\gamma (A)$
(or by $\gamma_K(A)$ if the dependence on $K$ needs to be
emphasized).

It is known that the set $\pi(K)$ forms a proper cone in the space
of $n \times n$ real matrices, the interior of $\pi(K)$ being the
subset consisting of $K$-positive matrices. Also, $\pi(K)$ is
polyhedral if and only if $K$ is polyhedral.  (See \cite{Tam 2},
\cite{S--V} or \cite{Bar 1}.)

A matrix $A$ is said to be {\it an automorphism of $K$} if $A$ is
invertible and $A, A^{-1}$ both belong to $\pi(K)$ or,
equivalently, $AK = K$.

 It is clear that if $K$ is a
simplicial cone with $n$ extreme rays then $K$ is linearly
isomorphic to ${\Bbb R}^n_+$.  The simplicial cones may be
considered as the simplest kind of cones. The next simplest kind
of cones, and also the one with which we will deal considerably in
this work, are the minimal cones.  Minimal cones were first
introduced and studied by Fiedler and Pt\'{a}k \cite{F--P}. We
call an $n$-dimensional polyhedral cone {\it minimal} if it has
precisely $n+1$ extreme rays.  Clearly, if $K$ is a minimal cone
with (pairwise distinct) extreme vectors $x_1, \ldots, x_{n+1}$,
then (up to multiples) these vectors satisfy a unique (linear)
relation. Also, a minimal cone is indecomposable if and only if
the relation for its extreme vectors is full, i.e., in the
relation the coefficient of each extreme vector is nonzero (see
\cite[Theorem 2.25]{F--P}). It is readily shown that every
decomposable minimal cone is the direct sum of an indecomposable
minimal cone and a simplicial cone.

In dealing with (nonzero) relations on (nonzero) extreme vectors
of a polyhedral cone, we find it convenient to write such
relations in the form
$$\alpha_1x_1+ \cdots + \alpha_px_p = \beta_1y_1+ \cdots +
\beta_qy_q,$$ where the extreme vectors $x_1, \ldots x_p, y_1,
\ldots, y_q$ are pairwise distinct and the coefficients $\alpha_1,
\ldots, \alpha_p, \beta_1, \ldots, \beta_q$ are all positive.
Clearly we have $p,q \ge 2$.

We call a relation on extreme vectors of a polyhedral cone {\it
balanced} if the number of nonzero terms on its two sides differ
by at most $1$.

An indecomposable minimal cone is said to be {\it of type
$(p,q)$}, where $p,q$ are positive integers such that $2\le p \le
q$, if the number of (nonzero) terms on the two sides of the
relation for its extreme vectors are respectively $p$ and $q$.
(We do not distinguish a relation with the one obtained from it by
interchanging the left side with the right side.)

Given positive integers with $2\le p\le q$ and $p+q = n+1$, one
can construct as follows an $n$-dimensional indecomposable minimal
cone of type $(p,q)$.  Choose any basis for ${\Bbb R}^n$, say $\{
x_1, \ldots, x_n \}$, and let $K$ be the polyhedral cone
$\fn{pos}\{ x_1, \ldots, x_n,x_{n+1} \}$, where $x_{n+1} = (x_1+
\cdots + x_p)-(x_{p+1}+ \cdots + x_n)$.  Then
$$ x_1 + \cdots + x_p = x_{p+1} + \cdots + x_{n+1}$$
is the (essentially) unique relation for the vectors $x_1, \ldots,
x_{n+1}$. As none of the vectors $x_1, \ldots, x_{n+1}$ can be
written as a nonnegative linear combination of the remaining
vectors, $x_1, \ldots, x_{n+1}$ are precisely (up to nonnegative
scalar multiples) all the extreme vectors of $K$. Therefore, $K$
is the desired indecomposable minimal cone.

Two proper cones $K_1, K_2$ are said to be {\it linearly
isomorphic} if there exists a linear isomorphism $P:~{\rm span}\,
K_2 \longrightarrow {\rm span}\, K_1$ such that $PK_2 = K_1$.

Using the following easy result in linear algebra, one can show
that indecomposable minimal cones of the same type are linearly
isomorphic.

\begin{lemma}\label{lemma 2.3} Let $\{ x_1, \ldots, x_k \}$ and $\{ y_1, \ldots, y_k \}$ be two families of
vectors in finite-dimensional vector spaces $V_1$ and $V_2$
respectively.  In order that there exists a linear mapping $T: V_1
\rightarrow V_2$ that satisfies $T(x_i) = y_i$ for $i = 1, \ldots,
k$, it is necessary and sufficient that $\alpha_1y_1+ \cdots +
\alpha_ky_k = 0$ is a relation for $y_1, \ldots, y_k$ whenever the
corresponding relation holds for $x_1, \ldots x_k$.
\end{lemma}

We also need the following known characterization of maximal faces
of an indecomposable minimal cone (\cite[Theorem 4.1]{Tam 1}):

\begin{thm} Let $K$ be an indecomposable minimal cone generated by
extreme vectors $x_1, \ldots, x_{n+1}$ that satisfy
$$x_1+ \cdots + x_p = x_{p+1} + \cdots + x_{n+1}.$$
Then for each pair ${\rm(}i,j{\rm)}, 1 \le i \le p \mbox{ and }
p+1 \le j \le n+1, \fn{pos}M_{ij}$ is a maximal face of $K$, where
$M_{ij} = \{ x_1, \ldots, x_{n+1} \}\setminus \{ x_i,x_j \}$.
Moreover, each maximal face of $K$ is of this form.
\end{thm}

Note that by the preceding theorem every maximal face, and hence
every nontrivial face, of an indecomposable minimal cone is a
simplicial cone in its own right.

Let $A \in \pi(K)$.  In this work we need the digraph $({\cal
E},{\cal P}(A,K))$, which is one of the four digraphs associated
with $A$ introduced by Barker and Tam (\cite{B--T}, \cite{T--B}).
It is defined in the following way: its vertex set is ${\cal E}$,
the set of all extreme rays of $K$; $(\Phi(x),\Phi(y))$ is an arc
whenever $\Phi(y) \subseteq \Phi(Ax)$.  If there is no danger of
confusion, (in particular, within proofs) we write $({\cal
E},{\cal P}(A,K))$ simply as $({\cal E},{\cal P})$.
 It is readily checked that if $K$ is the nonnegative orthant
${\Bbb R}^n_+$ then
 $({\cal E},{\cal P}(A,K))$ equals the usual digraph associated with $A^T$,
 the transpose of $A$.  (If $B = (b_{ij})$ is an $n\times n$
 matrix then by the usual digraph of $B$ we mean the digraph with
 vertex set $\{1, \ldots, n\}$ such that $(i,j)$ is an arc
 whenever $b_{ij}\ne 0$.)

It is not difficult to show that for any $A, B \in \pi(K)$, if
$\Phi(A) = \Phi(B)$ then $A,B$ are either both $K$-primitive or
both not $K$-primitive, and if they are, then $\gamma(A) =
\gamma(B)$. In Niu \cite{Niu} it is proved that if $K$ is a
polyhedral cone then for any $A, B \in \pi(K)$, we have $({\cal
E},{\cal P}(A,K)) = ({\cal E},{\cal P}(B,K))$ (as labelled
digraphs) if and only if $\Phi(A) = \Phi(B)$. So it is not
surprising to find that the digraph $({\cal E},{\cal
 P}(A,K))$ plays a role in determining a bound for $\gamma(A)$.
(When $K$ is nonpolyhedral, the situation is more subtle.  We
refer the interested readers to Tam \cite{Tam 4} for the details.)

Let $K$ be a polyhedral cone.  Then $\pi(K)$ is also polyhedral
and hence has finitely many faces.  Since $K$-primitive matrices
that belong to the relative interior of the same face of $\pi(K)$
share a common exponent, it follows that there are only finitely
many (integral) values that can be attained by the exponents of
$K$-primitive matrices.

For a proper cone $K$, we say $K$ {\it has finite exponent} if the
set of exponents of $K$-primitive matrices is bounded; then we
denote the maximum exponent by $\gamma(K)$ and refer to it as the
{\it exponent of $K$}.  If $K$ has finite exponent, then a
$K$-primitive matrix $A$ is said to be {\it exp-maximal} if
$\gamma(A) = \gamma(K)$.  By the above discussion, every
polyhedral cone has finite exponent.  An example of a proper cone
in ${\Bbb R}^3$ which does not have finite exponent will be given
in the final section of this paper.

We will make use of the concept of a {\it primitive digraph},
which can be defined as a digraph for which there is a positive
integer $k$ such that for every pair of vertices $i,j$ there is a
directed walk of length $k$ from $i$ to $j$; the least such $k$ is
referred to as the {\it exponent of the digraph}.   It is clear
that a nonnegative matrix is primitive if and only if its usual
digraph is primitive.  It is also well-known that primitive
digraphs are precisely strongly connected digraphs with the
greatest common divisor of the lengths of their circuits equal to
$1$.

It is known that for a $K$-nonnegative matrix $A$, if $A^p$ is
$K$-positive then so is $A^q$ for every $q \ge p$.  This follows
from the fact that if $B$ is a $K$-nonnegative matrix such that
$Bu\in \partial K$ for some $u\in \fn{int}K$ then we have $BK
\subseteq \Phi(Bu) \subseteq \partial K$.  The same fact also
implies that the action of a $K$-nonnegative matrix $A$ on a
vector $x$ in $K$ enjoys a similar property --- if $A^ix$ belongs
to $\fn{int}K$, then so does $A^jx$ for all positive integers $j
> i$.

If $A$ is a $K$-nonnegative matrix and if $p$ is a positive
integer such that $A^pF\subseteq F$ for some nontrivial face $F$,
then $A^{kp}F\subseteq F$ for all positive integers $k$ and hence
$A$ cannot be $K$-primitive.  This shows that the positive powers
of a $K$-primitive matrix are all $K$-irreducible.

To study the exponents of $K$-primitive matrices, we make use of
the concept of local exponent defined in the following way.  (For
definition of local exponent of a primitive matrix, see \cite
[Section 3.5]{B--R}.)   For any $K$-nonnegative matrix $A$, not
necessarily $K$-primitive or $K$-irreducible, and any $0 \ne x \in
K$, by the {\it local exponent of $A$ at $x$}, denoted by
$\gamma(A,x)$,
 we mean the smallest nonnegative integer $k$ such that $A^k x \in {\rm int}~K$. If no such
 $k$ exists, we set $\gamma(A,x)$ equal $\infty$.  (If $A$ is a primitive matrix and $e_j$ is the
$j$th standard unit vector, then $\gamma(A,e_j)$ equals the
smallest integer $k$ such that all elements in column $j$ of $A^k$
are nonzero.) Clearly, $A$ is $K$-primitive if and only if the set
of local exponents of $A$ is bounded; in this case, $\gamma(A)$ is
equal to $\max \{ \gamma(A,x): 0 \ne x \in K \}$, which is also
the same as the maximum taken over all nonzero extreme vectors of
$K$. By a compactness argument Barker \cite{Bar 1} has shown that
the $K$-primitivity of $A$ is equivalent to the apparently weaker
condition --- which is also the definition adopted by him for
$K$-primitivity --- that all local exponents of $A$ are finite.

By the definition of $({\cal E}, {\cal P}(A,K))$, we have

\begin{fact} \rm \label{fact 2.1}
If there is a path in $({\cal E}, {\cal P}(A,K))$ of
length $k$ from $\Phi(x)$ to $\Phi(y)$, then $\Phi(A^kx) \supseteq
\Phi(y)$.
\end{fact}

By Fact \ref{fact 2.1} we obtain the following\,{\rm:}

\begin{fact}\rm \label{fact 2.3}
Let $K$ be a polyhedral cone.  If $A$ is a $K$-nonnegative matrix
such that the digraph $({\cal E},{\cal P}(A,K))$ is primitive,
then $A$ is a $K$-primitive matrix with exponent less than or
equal to that of the primitive digraph $({\cal E},{\cal P}(A,K))$.
\end{fact}

To show the preceding fact, let $k$ denote the exponent of the
primitive digraph $({\cal E},{\cal P})$ and let $x\in \fn{Ext} K$.
Then there is a path of length $k$ in the said digraph from
$\Phi(x)$ to $\Phi(y)$ for any $y\in \fn{Ext}K$.  By Fact
\ref{fact 2.1} we have $\Phi(A^kx)\supseteq \Phi(y)$.  Since this
is true for every $y\in \fn{Ext} K$, it follows that $A^kx\in
\fn{int}K$.  But $x$ is an arbitrary nonzero extreme vector of
$K$, so $A^k$ is $K$-positive.  Hence $A$ is $K$-primitive and
$\gamma(A)\le k$.

Fact \ref{fact 2.1} implies also the following\,{\rm:}

\begin{fact} \rm \label{fact 2.2}
Let $A \in \pi(K)$ and let $x, y \in \fn{Ext}K$. Suppose that
$\gamma(A,y)$ is finite.  If there is a path in $({\cal E}, {\cal
P}(A,K))$ of length $k$ from $\Phi(x)$ to $\Phi(y)$, then
$\gamma(A,x)$ is also finite and we have $\gamma(A,x) \le
k+\gamma(A,y)$.
\end{fact}

 Two
cone-preserving maps $A_1\in \pi (K_1)$ and $A_2 \in \pi(K_2)$ are
said to be {\it cone-equivalent} if there exists a linear
isomorphism $P$ such that $PK_2=K_1$ and $P^{-1}A_1P=A_2$.

\begin{fact}\rm \label{fact2.5} Let $K_1, K_2$ be proper cones in ${\Bbb R}^n$.  Suppose that
$A_1\in \pi(K_1)$ and $A_2\in
\pi(K_2)$ are cone-equivalent.  Then\,{\rm:}
\begin{enumerate}
\item[(i)] $A_1$ and $A_2$ are similar. \item[(ii)] The cones
$K_1,K_2$ are linearly isomorphic. \item[(iii)] The digraphs
$({\cal E},{\cal P}(A_1,K_1)), ({\cal E},{\cal P}(A_2,K_2))$ are
isomorphic. \item[(iv)] Either $A_1$ is $K_1$-primitive and $A_2$
is $K_2$-primitive or they are not, and if they are, then
$\gamma_{K_1}(A_1) = \gamma_{K_2}(A_2)$. \item[(v)] For any $x\in
K_2$, $\gamma(A_2,x) = \gamma(A_1,Px)$.
\end{enumerate}
\end{fact}

  Also, it is clear that if
$K_1$ and $K_2$ are linearly isomorphic cones, then either $K_1,
K_2$ both have finite exponent or they both do not have, and if
they both have, then $\gamma(K_1) = \gamma(K_2)$.

Under inclusion as the partial order, the set of all faces of $K$,
denoted by ${\cal F}(K)$, forms a lattice with meet and join given
respectively by: $F\wedge G = F \cap G$ and $F\vee G = \Phi(F \cup
G)$. Two proper cones $K_1$, $K_2$ are said to be {\it
combinatorially equivalent} if their face lattices ${\cal F}(K_1)$
and ${\cal F}(K_2)$ are isomorphic (as lattices).

For minimal cones, the concepts of ``linearly isomorphic" and
``combinatorially equivalent" are equivalent.

\begin{theorem}\label{theorem2.5} Let $K_1, K_2$ be minimal cones of dimension
$n_1,n_2$ respectively.  Suppose that for $j = 1,2, K_j =
K_j^\prime \oplus K_j^{\prime\prime}$, where $K_j^\prime$ is a
simplicial cone and $K_j^{\prime\prime}$ is an indecomposable
minimal cone of type $(p_j,q_j)$.  The following conditions are
equivalent\,{\rm:}
\begin{enumerate}
\item[{\rm(i)}] $n_1 = n_2$ and $(p_1,q_1) = (p_2,q_2)$.
\item[{\rm(ii)}] $K_1, K_2$ are linearly isomorphic.
\item[{\rm(iii)}] $K_1, K_2$ are combinatorially equivalent.
\end{enumerate}
\end{theorem}

{\it Proof}.  For $j = 1,2$, let $d_j$ be the dimension of
$K_j^\prime$.

 (i)$\Rightarrow$(ii): We have
  $$d_1 = n_1-\dim K_1^{\prime\prime} = n_1-(p_1+q_1-1) = n_2-(p_2+q_2-1) = n_2-\dim K_2^{\prime\prime} = d_2;$$
so $K_1^\prime$ and $K_2^\prime$ are linearly isomorphic, being
simplicial cones of the same dimension.  On the other hand,
$K_1^{\prime\prime}$ and $K_2^{\prime\prime}$ are also linearly
isomorphic, as they are indecomposable minimal cones of the same
type.  Therefore, $K_1$ and $K_2$ are linearly isomorphic.

(ii)$\Rightarrow$(iii):  Obvious.

(iii)$\Rightarrow$(i):  As is well-known, if $K$ is a polyhedral
cone then $K$ has faces of all possible dimensions from $0$ to
$\fn{dim}K$.  So $\fn{dim}K$ is equal to the length of a maximal
chain in the face lattice ${\cal F}(K)$ of $K$.  Since the face
lattices ${\cal F}(K_1)$ and ${\cal F}(K_2)$ are isomorphic, they
have maximal chains of the same length.  Hence, we have $n_1 =
n_2$.

To proceed further, we need to establish two assertions first:\\

{\bf Assertion 1}.  Let $K$ be an $n$-dimensional minimal cone, $n
\ge 3$.  Suppose that $K = K^\prime \oplus K^{\prime\prime}$,
where $K^\prime$ is a $d$-dimensional simplicial cone ($0 \le d
\le n-3$) and $K^{\prime\prime}$ is an ($n-d$)-dimensional
indecomposable minimal cone of type $(p,q)$.  Then $K$ has $d+pq$
maximal faces, $d$ of which are minimal cones and the remaining
are simplicial cones.

{\it Proof}.  Since $K$ is the direct sum of $K^\prime$ and
$K^{\prime\prime}$, the maximal faces of $K$ are precisely those
of the form $M^\prime \oplus K^{\prime\prime}$ or $K^\prime \oplus
M^{\prime\prime}$, where $M^\prime, M^{\prime\prime}$ denote
respectively a maximal face of $K^\prime$ and a maximal face of
$K^{\prime\prime}$.  There are $d$ maximal faces of the first
kind, each of which, being the direct sum of a simplicial cone and
a minimal cone, is a minimal cone in its own right.  In view of
Theorem A, maximal faces of the indecomposable minimal cone
$K^{\prime\prime}$ are themselves simplicial cones and there are
altogether $pq$ of them; hence, $K$ has $pq$ maximal faces of the
second
kind, each of which is a simplicial cone. \hfill$\blacksquare$\\

{\bf Assertion 2}.  A minimal cone cannot be combinatorially
equivalent to a simplicial cone.

{\it Proof}. As we have already noted, any two combinatorially
equivalent polyhedral cones have the same dimension.  Now an
$n$-dimensional simplicial cone has $n$ extreme rays, whereas an
$n$-dimensional minimal cone has $n+1$ extreme rays.  So a minimal
cone and a simplicial cone cannot be combinatorially equivalent.
\hfill$\blacksquare$\\

Now back to the proof of the theorem.  Let $\Psi$ be a lattice
isomorphism between ${\cal F}(K_1)$ and ${\cal F}(K_2)$.  Clearly
$\Psi$ provides a one-to-one correspondence between the maximal
faces of $K_1$ and those of $K_2$.  Note that if $M_1$ is a
maximal face of $K_1$ that corresponds to the maximal face $M_2$
of $K_2$, then $M_1$ and $M_2$ are combinatorially equivalent, as
$\Psi$ induces a lattice isomorphism between ${\cal F}(M_1)$ and
${\cal F}(M_2)$.  In view of Assertion $2$, under $\Psi$, maximal
faces of $K_1$ which are themselves minimal (respectively,
simplicial) cones correspond to maximal faces of $K_2$ which are
themselves minimal (respectively, simplicial) cones.  By Assertion
$1$, we have $d_1 = d_2$ and $p_1q_1 = p_2q_2$.  But we have
already shown that $n_1 = n_2$ and also we have $p_j+q_j =
n_j-d_j+1$ for $j = 1,2$, it follows that we have $(p_1,q_1) =
(p_2,q_2)$. \hfill$\blacksquare$\\

\setcounter{equation}{0}
\section{Upper bounds for exponents}

Hereafter, for every pair of positive integers $m,n$ with $3\le n
\le m$, we denote by ${\cal P}(m,n)$ the set of all
$n$-dimensional polyhedral cones with $m$ extreme rays. Note that
we start with $n = 3$ as the cases $n = 1$ or $2$ are trivial.
Also, when $m = 3$, in order that ${\cal P}(m,n)$ is nonvacuous,
$n$ must be $3$.

The following theorem of Sedl\'{a}\u{c}ek \cite{Sed} and Dulmage
and Mendelsohn \cite{D--M} (see, for instance, \cite[Theorem
3.5.4]{B--R}) gives an upper bound for the exponent of a primitive
matrix $A$ in terms of lengths of circuits in the digraph of $A$.

\begin{thm}  Let $A$ be an $n\times n$ primitive matrix.  If $s$ is the length of
the shortest circuit in the digraph of $A$, then
$\gamma(A) \le n+s(n-2).$
\end{thm}

By setting $s=n-1$ in Theorem B, one recovers the sharp general
upper bound $(n-1)^2+1$, due to Wielandt \cite{Wie}, for exponents
of $n\times n$ primitive matrices.

The next lemma gives an analogous result on the local exponents of
a cone-preserving map, which is essential to our work.

If $D$ is a digraph, $v$ is a vertex of $D$ and $W$ is a nonempty
subset of the vertex set of $D$, then we say $v$ {\it has access
to} $W$ if there is a path from $v$ to a vertex of $W$. In this
case, the length (i.e., the number of edges) of the shortest path
from $v$ to a vertex of $W$ is referred to as the {\it distance}
from $v$ to $W$.  If $v$ belongs to $W$, the distance is taken to
be zero.

 For a square matrix $C$, we denote by $m_C$ the
{\it degree of the minimal polynomial of $C$}.

\begin{lemma}\label{lemma1}
Let $K$ be a proper cone and let $A \in \pi(K)$. Let $\Phi(x)$ be
a vertex of $({\cal E}, {\cal P}(A,K))$ which is at a distance
$w\, (\ge 0)$ to a circuit ${\cal C}$ of length $l$.  Suppose that
$A^l$ is $K$-irreducible, or that the circuit ${\cal C}$ contains
a vertex $\Phi(u)$ for which $\gamma(A,u)$ is finite.  Then
$\gamma(A,x)$ is finite and
$$\gamma(A,x) \le w+(m_{A^l}-1)l \le
w+(m_A-1)l \le w+(n-1)l.$$
\end{lemma}

{\it Proof}.  Let ${\cal C}\!\!: \Phi(x_1) \rightarrow \Phi(x_2)
\rightarrow \cdots \rightarrow \Phi(x_l) \rightarrow \Phi(x_1)$ be
the circuit under consideration.  (Here, for convenience, we
represent the arc $(\Phi(x),\Phi(y))$ by $\Phi(x) \rightarrow
\Phi(y)$.) Without loss of generality,  we may assume that the
distance from $\Phi(x)$ to $\Phi(x_1)$ is $w$. Since there is a
path of length $l$ from $\Phi(x_1)$ to itself, by Fact \ref{fact
2.1} we have $\Phi(A^lx_1) \supseteq \Phi(x_1)$, which implies the
following chain of inclusions:
$$
\Phi(x_1) \subseteq \Phi(A^lx_1) \subseteq \Phi(A^{2l}x_1)
\subseteq \cdots \subseteq \Phi(A^{jl}x_1) \subseteq
\Phi(A^{(j+1)l}x_1) \subseteq \cdots.
$$
Let $p$ denote the dimension of the subspace $\fn{span}\{
(A^l)^jx_1\!\!: j = 0, 1, \ldots \}$. By the above chain of
inclusions, clearly the face $\Phi((A^l)^{ p-1}x_1)$ contains the
vectors $x_1, A^lx_1, \ldots, (A^l)^{p-1}x_1$, which are linearly
independent and hence form a basis for the said subspace; so
$\Phi((A^l)^{p-1}x_1)$ includes, and hence is equal to,
$\Phi(\fn{span} \{ (A^l)^jx_1: j = 0, 1, \ldots \} \cap K)$.  Note
that the latter face is the smallest $A^l$-invariant face of $K$
that contains $x_1$. If $A^l$ is $K$-irreducible, the latter face
is clearly $K$.  On the other hand, if the circuit ${\cal C}$
contains a vertex $\Phi(u)$ for which $\gamma(A,u)$ is finite,
then by Fact \ref{fact 2.2} $\gamma(A,x_1)$ is also finite.  Hence
$A^jx_1\in \fn{int}K$ for all positive integers $j$ sufficiently
large and, as a consequence, the smallest $A^l$-invariant face of
$K$ that contains $x_1$ is $K$. In either case, we have,
$\Phi((A^l)^{p-1}x_1) = K$ and so $(A^l)^{p-1}x_1 \in \fn{int}K$;
hence $\gamma(A,x_1) \le (p-1)l$. Then by Fact \ref{fact 2.2} we
have
$$\gamma(A,x) \le w+ \gamma(A,x_1) \le w+(p-1)l.$$
It is clear that $p \le m_{A^l}$.  But we also have $m_{A^l} \le
m_A\le n$, so the desired inequalities follow.
\hfill$\blacksquare$\\

\begin{lemma}\label{lemma3.2} Let $K\in {\cal P}(m,n)$ and let $A$
be a $K$-primitive matrix.  If $\Phi(x)$ is a vertex of $({\cal
E},{\cal P}(A,K))$ which has access to a circuit of length $l$,
then
$$\gamma(A,x) \le m+(m_A-2)l \le m+(n-2)l.$$
\end{lemma}

{\it Proof}.  This follows from Lemma \ref{lemma1}, as the
distance from $\Phi(x)$ to ${\cal C}$ is at most $m-l$ and also
$A^l$ is $K$-irreducible. \hfill$\blacksquare$

Using Lemma \ref{lemma1}, one can also readily deduce the
following result.

\begin{coro}\label{coro 3.3}
Let $K\in {\cal P}(m,n)$, and let $A \in \pi(K)$.  Suppose that
the digraph $({\cal E}, {\cal P}(A,K))$ is strongly connected.
Let $s$ be the shortest circuit length of the digraph.  If $A^s$
is $K$-irreducible, then $A$ is $K$-primitive and $\gamma(A) \le
m+s(m_A-2)$.
\end{coro}

It is known that if $K$ is a polyhedral cone with $m$ extreme
rays, then a $K$-nonnegative matrix $A$ is $K$-primitive if $A^j$
are $K$-irreducible for $j = 1, \ldots, 2^m-1$ (see \cite[Theorem
2]{Bar 1}).  The preceding corollary tells us that when the
digraph $({\cal E},{\cal P})$ is strongly connected, to show the
$K$-primitivity of $A$, it suffices to check the
$K$-irreducibility of only one positive power of $A$.

Clearly, the following result of Niu [Niu] is a consequence of
Corollary \ref{coro 3.3}:

\begin{thm}
Let $K\in {\cal P}(m,n)$ and let $A$ be $K$-primitive. If the
digraph $({\cal E}, {\cal P}(A,K))$ is strongly connected and $s$
is the length of the shortest circuit in $({\cal E}, {\cal
P}(A,K))$, then $\gamma(A) \le m+s(m-2)$.
\end{thm}

In Theorem C, by choosing $K=\IR^n_+$ we recover Theorem B.

It is known (see \cite{S--V}) that if $A$ is $K$-irreducible, then
$(I+A)^{n-1}$ is $K$-positive (where $n$ is the dimension of $K$).
 Hartwig and Neumann \cite{H--N} have shown that in the nonnegative matrix case
 the result can be strengthened by
replacing $n$ by $m_A$, the degree of the minimal polynomial of
$A$. Now we can show that the latter improvement is also valid for
a cone-preserving map on a proper cone.

\begin{coro} If $A \in \pi(K)$ is $K$-irreducible, then
$(I+A)^{m_A-1}$ is $K$-positive.
\end{coro}

{\it Proof}. If $A$ is $K$-irreducible, then clearly $I+A$ is also
$K$-irreducible and in the digraph $({\cal E},{\cal P}(I+A,K))$
there is a loop at each vertex.  By Lemma \ref{lemma1},
$\gamma(I+A,x) \le m_{I+A}-1$ for every $x\in\fn{Ext} K$.
But $m_{I+A} = m_A$, so $(I+A)^{m_A-1}$ is $K$-positive. \hfill$\blacksquare$\\

It is also possible to provide a direct proof for Corollary 3.4,
one that does not involve the digraph $({\cal E},{\cal P})$.

We denote by ${\cal N}(A)$ the nullspace of $A$.  It is easy to
show that for any $A \in \pi(K), {\cal N}(A) \cap K = \{ 0 \}$ if
and only if the outdegree of each vertex of $({\cal E},{\cal P})$
is positive.  As a consequence, for any $K$-primitive matrix $A$,
the digraph $({\cal E},{\cal P})$ has at least one circuit.

In contrast with the nonnegative matrix case, the digraph $({\cal
E},{\cal P})$ associated with a $K$-primitive matrix $A$ (where
$K$ is polyhedral) need not be strongly connected.  (Many such
examples can be found in \cite{Tam 4}.)  Nevertheless, every
vertex of $({\cal E},{\cal P})$ has access to some circuit of
$({\cal E}, {\cal P})$.  This makes it possible to apply Lemma
\ref{lemma3.2} to obtain bounds for the exponents of $K$-primitive
matrices.

As yet another application of Lemma \ref{lemma1}, we obtain the
following result, which is an extension of the corresponding
result for a symmetric primitive matrix (cf. \cite [Theorem
3.5.3]{B--R}).  Recall that a digraph $D$ is said to be {\it
symmetric} if for every pair of vertices $u,v$ of $D$, $(u,v)$ is
an arc if and only if $(v,u)$ is an arc.

\begin{coro}
Let $A \in \pi(K)$.  If $A$ is $K$-primitive and the digraph
$({\cal E},{\cal P}(A,K))$ is symmetric, then
$$\gamma(A) \le 2(m_{A^2}-1) \le 2(m_A-1).$$
\end{coro}
{\it Proof}.  Since $A$ is $K$-primitive, ${\cal N}(A) \cap K = \{
0 \}$; hence the digraph $({\cal E},{\cal P})$ has an outgoing
edge (possibly a loop) at each vertex. As the digraph $({\cal
E},{\cal P})$ is symmetric, it follows that $({\cal E},{\cal
P}(A^2,K))$ has a loop at each vertex.  By Lemma \ref{lemma1}  we
have $\gamma(A^2) \le m_{A^2}-1$ and hence $\gamma(A) \le
2(m_{A^2}-1) \le 2(m_A-1)$. \hfill$\blacksquare$\\

It is clear that for any $K$-primitive matrix $A, m_A \ge 2$. When
$m_A = 2$, more can be said.

\begin{lemma}\label{lemma2}
Let $A$ be $K$-primitive.  If $m_A = 2$ then $\gamma(A) = 1$ or $2$.
\end{lemma}
{\it Proof}. Since $m_A = 2$ (and $A$ is a real matrix), there
exist real numbers $a, b$ such that $A^2+aA+bI = 0$.  Clearly, $a,
b$ cannot be both zero, as $A$ is not nilpotent.  By the
pointedness of the cone $\pi(K)$, at least one of $a,b$ is
negative.  If $b < 0$ and $a \ge 0$, then $A^2$ belongs to the
face $\Phi(I)$ (of $\pi(K)$) and so it must be $K$-reducible,
which is a contradiction.  If $a < 0$ and $b \ge 0$, then $A^2 \in
\Phi(A)$ or $A^2 \le \alpha A$ for some $\alpha > 0$, which
implies that all positive powers of $A$ lie in $\Phi(A)$.  But
$A^p$ is $K$-positive or, equivalently, belongs to
$\fn{int}\pi(K)$ for $p$ sufficiently large, it follows that in
this case we must have $\Phi(A) = \pi(K)$, or in other words,
$\gamma(A) = 1$.  In the remaining case when $a, b$ are both
negative, $A^2$ is a positive linear combination of $A$ and $I$
and hence lies in $\fn{ri}\Phi(A+I)$. Then one readily shows that
all positive powers of $A$ also lie in $\fn{ri}\Phi(A+I)$.  By the
$K$-primitivity of $A$, $A^p$ belongs to $\fn{int}\pi(K)$ for $p$
sufficiently large. This implies that $\Phi(A+I) = K$.  As $A^2$
is a positive linear combination of $A$ and $I$, it follows that
 $A^2$ also belongs to $\fn{int}\pi(K)$; so we have
$\gamma(A) \le 2$. This completes the proof. \hfill$\blacksquare$

\setcounter{equation}{0}
\section{Special digraphs for $K$-primitive matrices}

The results of Section 3 may suggest that for a $K$-primitive
matrix $A$ the longer the shortest circuit in $({\cal E},{\cal
P})$ is, the larger is the value of $\gamma(A)$.  Given a pair of
positive integers $m,n$ with $3\le n \le m$, what can we say about
digraphs on $m$ vertices, with the length of the shortest circuit
equal to $m-1$, which can be identified (up to graph isomorphism)
with $({\cal E},{\cal P}(A,K))$ for some pair $(A,K)$ where $K\in
{\cal P}(m,n)$ and $A$ is a $K$-primitive matrix ?  It turns out
that such digraphs must be primitive.  When $m \ge 4$, apart from
the labelling of its vertices (or, in other words, up to graph
isomorphism), there are two such digraphs, which are given by
Figure $1$ and Figure $2$. When $m = 3$, there is one more
digraph.  (See Lemma \ref{lemma3} below.)

Note that if $K$ is a polyhedral cone with $m$ extreme rays then
for any $K$-primitive matrix $A$, the length of the shortest
circuit in $({\cal E},{\cal P})$ is at most $m-1$.  This is
because, if the length of the shortest circuit is $m$, then the
digraph must be a circuit of length $m$ and, as a consequence,
$A^m$ is $K$-reducible, which is impossible.

In what follows when we say the digraph $({\cal E},{\cal P})$ is
given by Figure $1$ (or by other figures), we mean the digraph is
given either by the figure up to graph isomorphism or by the
figure as a labelled digraph.  In most instances, we mean it in
the former sense but in a few instances we mean it in the latter
sense.  It should be clear from the context in what sense we mean.
(For instance, in parts (i) and (iii) of Lemma \ref{lemma4.2} we
mean the former sense, but in part (ii) we mean the latter sense.)

We will obtain certain geometric properties of $K$ when
 $K$ is a non-simplicial
polyhedral cone with $m\,(\ge 4)$ extreme rays for which there
exists a $K$-primitive matrix $A$ such that $({\cal E},{\cal P})$
is given by Figure $1$ or Figure $2$.  More precisely, we show
that if $({\cal E},{\cal P})$ is given by Figure $1$ then $K$ is
indecomposable; if $({\cal E},{\cal P})$ is given by Figure $2$
and if $K$ is decomposable, then $m$ is odd and $K$ is the direct
sum of a ray and an indecomposable minimal cone with a balanced
relation on its extreme vectors.  We also obtain some properties
on the corresponding $K$-primitive matrix $A$.  It will be shown
that
 if $K$ is a polyhedral cone
with $m$ extreme rays, then for any $K$-primitive matrix $A$,
$\gamma(A) \le (m_A-1)(m-1)+1$, where the equality holds only if
the digraph $({\cal E}, {\cal P})$ is given by Figure $1$. As a
consequence, $(n-1)(m-1)+1$ is an upper bound for $\gamma(K)$ when
$K\in {\cal P}(m,n)$.

\begin{lemma}\label{lemma3}
Let $K\in {\cal P}(m,n)$ {\rm(}$3\le n \le m${\rm)} and let $A$ be
a $K$-primitive matrix. Then the length of the shortest circuit in
the digraph $({\cal E},{\cal P}(A,K))$ equals $m-1$ if and only if
the digraph $({\cal E},{\cal P}(A,K))$ is, apart from the
labelling of its vertices, given by Figure $1$ or Figure $2$, or
{\rm(}in case $m = n = 3${\rm)} by the digraph of order $3$ whose
arcs are precisely all possible arcs between every pair of
distinct vertices\,{\rm:}
\newsavebox{\figone}
\savebox{\figone}{ \put(-145,62){$x_1$}
\put(-167,36){\vector(1,1){20}} \put(-134.5,56){\vector(1,-1){19}}
\put(-180,27){$x_m$} \put(-162,30){\vector(1,0){40}}
\put(-114,27){$x_2$} \put(-99,30){\vector(1,0){40}}
\put(-53,27){$x_3$} \put(-163,-5){\vector(-1,3){8}}
\put(-52,20){\vector(-1,-3){8.5}} \put(-169,-16){$x_{m-1}$}
\put(-68,-16){$x_4$} \put(-134,-40){\vector(-1,1){17}}
\put(-117,-43){$\cdots$} \put(-70,-22){\vector(-1,-1){17}} }
\[
\put(0,0){\usebox{\figone}} \put(-132,-65){\rm Figure 1.}
\put(212,0){\usebox{\figone}} \put(84,60){\vector(3,-1){71}}
\put(80,-65){\rm Figure 2.}
\]
$($For simplicity, we label the vertex $\Phi(x_i)$ simply by
$x_i$.$)$.
\end{lemma}

{\it Proof}.  We treat only the ``only if" part,  the ``if" part
being obvious.

 It is not difficult to show that there are
precisely three non-isomorphic primitive digraphs of order three
with shortest circuit length two, namely, the digraphs given by
Figure $1$, Figure $2$ and the one with all possible arcs between
every pair of distinct vertices.  So there is no problem when $m =
3\,(= n)$. Hereafter, we assume that $m \ge 4$.

Let $x_1, \ldots, x_m$ denote the pairwise distinct extreme
vectors of $K$.
 Let $A$ be a $K$-primitive matrix such that the length of the
 shortest circuit in the digraph $({\cal E},{\cal P})$ is $m-1$. Without loss of generality, we may
assume that the digraph $({\cal E},{\cal P})$ contains the circuit
${\cal C}$ (of length $m-1$) that is made up of the arcs
$(\Phi(x_m),\Phi(x_2))$ and $(\Phi(x_j),\Phi(x_{j+1}))$ for $j =
2, 3, \ldots, m-1$. Being a circuit of shortest length,  ${\cal
C}$ cannot contain any chord, nor can it have loops at its
vertices. If there is no arc from a vertex of ${\cal C}$ to the
remaining vertex $\Phi(x_1)$, then we have $A\Phi(x_m) =
\Phi(x_2)$ and $A\Phi(x_j) = \Phi(x_{j+1})$ for $j = 2, 3, \ldots,
m-1$ and it will follow that $A^{m-1}x_m$ is a positive multiple
of $x_m$, hence $A^{m-1}$ is $K$-reducible, which contradicts the
assumption that $A$ is $K$-primitive. So there is at least one arc
from a vertex of ${\cal C}$ to $\Phi(x_1)$, say,
$(\Phi(x_m),\Phi(x_1))$ is one such arc.  Similarly, there is also
an arc from $\Phi(x_1)$ to a vertex of ${\cal C}$.  Since the
length of the shortest circuit in the digraph $({\cal E},{\cal
P})$ is $m-1$, there cannot be an arc of the form
$(\Phi(x_1),\Phi(x_j))$ with $4 \le j \le m$.  So
$(\Phi(x_1),\Phi(x_2))$ and $(\Phi(x_1),\Phi(x_3))$ are the only
possible arcs with initial vertex $\Phi(x_1)$, and at least one of
them must be present.

We treat the case when the arc $(\Phi(x_1),\Phi(x_2))$ is present
first.  Clearly, none of the arcs $(\Phi(x_j),\Phi(x_1)), \mbox{
for }j = 2, \ldots, m-2$, can be present, as the length of the
shortest circuit in $({\cal E},{\cal P})$ is $m-1$.  However, it
is possible that $(\Phi(x_{m-1}),\Phi(x_1))$ is an arc, provided
that $(\Phi(x_1),\Phi(x_3))$ is not an arc.   Note that the
digraph that consists of the circuit ${\cal C}$ and the arcs
$(\Phi(x_m),\Phi(x_1)),(\Phi(x_1),\Phi(x_2)),(\Phi(x_{m-1}),\Phi(x_1))$
is isomorphic with the one given by Figure 2.  So, in this case,
up to graph isomorphism, the digraph $({\cal E},{\cal P})$ is
given by Figure 1 or Figure 2.

Now we consider the case when the arc $(\Phi(x_1),\Phi(x_2))$ is
absent.  Then the arc $(\Phi(x_1),\Phi(x_3))$ is present and the
digraph $({\cal E},{\cal P})$ contains two circuits of length
$m-1$.  As each of these two circuits cannot contain a chord,
$(\Phi(x_2),\Phi(x_1))$ is the only possible remaining arc.  If
the arc $(\Phi(x_2),\Phi(x_1))$ is absent, then a little
calculation shows that $A^{m-1}x_m$ is a positive multiple of
$x_m$, which violates the assumption that $A$ is $K$-primitive. On
the other hand, if the arc $(\Phi(x_2),\Phi(x_1))$ is present,
then the digraph $({\cal E},{\cal P})$ is isomorphic with the one
given by Figure $2$.  This completes the proof.
\hfill$\blacksquare$

Note that Figure $1$ is the same as the (unique) digraph
associated with an $m\times m$ primitive matrix whose exponent
attains Wielandt's bound $m^2-2m+2$ (see \cite{B--R}).

To proceed further, we need to manipulate with the relations on
the extreme vectors of a polyhedral cone.  We now explain the
relevant terminology.

Let $R$ be a relation on the extreme vectors of a polyhedral cone
$K$. Suppose that the vectors that appear in $R$ come from $p\,
(\ge 2)$ different indecomposable summands of $K$, say, $K_1,
\ldots, K_p$. To be specific, let $R$ be given by: $\sum_{i \in M}
\alpha_ix_i = \sum_{j \in N} \beta_jy_j$, where $M,N$ are finite
index sets, each with at least two elements and the $\alpha_i$s,
$\beta_j$s are all positive real numbers.  For each $r = 1,
\ldots, p$, let $M_r = \{i \in M: x_i \in K_r \}$ and $N_r = \{j
\in N: y_j \in K_r \}$.  Then for each fixed $r$, rewriting
relation $R$, we obtain
$$\sum_{i\in M_r}\alpha_ix_i - \sum_{j\in N_r}\beta_jy_j = \sum_{j\in N \setminus N_r}\beta_jy_j - \sum_{i\in M \setminus M_r}\alpha_ix_i.$$
Now the vector on the left side of the above relation belongs to
$\fn{span}K_r$, whereas the one on the right side belongs to
$\sum_{s\ne r}\fn{span}K_s$.  But $\fn{span}K_r \cap \sum_{s\ne
r}\fn{span}K_s = \{ 0 \}$ (as $K_1, \ldots, K_p$ are pairwise
distinct indecomposable summands of $K$), so it follows that we
have the relation
$$\sum_{i\in M_r}\alpha_ix_i = \sum_{j\in N_r}\beta_jy_j,$$
which we denote by $R_r$.  This is true for each $r$.  It is clear
that relation $R$ can be obtained by adding up relations $R_1,
\ldots, R_p$.  In this case, we say relation $R$ {\it splits into
the subrelations $R_1, \ldots, R_p$.}  Note that each $R_r$ has at
least four (nonzero) terms. So when we pass from the relation $R$
to one of its subrelations $R_r$, the number of terms involved in
the relation decreases by at least four.

Recall that an $n\times n$ complex matrix $A$ is said to be {\it
non-derogatory} if every eigenvalue of $A$ has geometric
multiplicity $1$ or, equivalently, if the minimal and
characteristic polynomials of $A$ are identical.  (See, for
instance, \cite[Theorem 3.3.15]{H--J}.)

\begin{lemma}\label{lemma4.2}
Let $K\in {\cal P}(m,n)$ {\rm(}$3\le n \le m${\rm)}.  Let $A$ be a
$K$-nonnegative matrix. Suppose that the digraph $({\cal E},{\cal
P}(A,K))$ is given by Figure $1$ or Figure $2$.  Then\,{\rm:}
\begin{enumerate}
\item [{\rm (i)}] $A$ is $K$-primitive, nonsingular,
non-derogatory, and has a unique annihilating polynomial of the
form $t^m-ct-d$, where $c,d > 0$.

\item [{\rm (ii)}] $\gamma(A)$ equals $\gamma(A,x_1)$ or
$\gamma(A,x_2)$ depending on whether the digraph $({\cal E},{\cal
P}(A,K))$ is given by Figure $1$ or Figure $2$.  In either case,
$\fn{max}_{1\le i \le m}\gamma(A,x_i)$ is attained at precisely
one $i$.

\item[{\rm (iii)}] Assume, in addition, that $K$ is
non-simplicial. If $({\cal E},{\cal P}(A,K))$ is given by Figure
$1$ then $K$ must be indecomposable.  If the digraph is given by
Figure $2$ then either $K$ is indecomposable or $m$ is odd and $K$
is the direct sum of a ray and an indecomposable minimal cone with
a balanced relation for its extreme vectors.
\end{enumerate}
\end{lemma}

{\it Proof}.  (i) Since the digraphs given by Figure $1$ and
Figure $2$ are primitive, by Fact \ref{fact 2.3} $A$ is
$K$-primitive.  To show that $A$ is nonsingular, we treat the case
when the digraph $({\cal E},{\cal P})$ is given by Figure $1$, the
argument for the other case being similar.  Then for $j = 1,
\ldots, m-1, Ax_j$ is a positive multiple of $x_{j+1}$. So $x_2,
x_3, \ldots, x_m$ all belong to ${\cal R}(A)$, the range space of
$A$.  On the other hand, since $x_2 \in {\cal R}(A)$ and $Ax_m$ is
a positive linear combination of $x_1$ and $x_2$, we also have
$x_1 \in {\cal R}(A)$. Therefore, regarded as a linear map $A$ is
onto and hence is nonsingular.

To establish the second half of this part, we may assume that
$\rho(A) = 1$ as $\rho(A) > 0$, $A$ being $K$-primitive.  We first
deal with the case when $({\cal E},{\cal P}(A,K))$ is given by
Figure $1$. Since $A$ is $K$-primitive, $A^T$ is $K^*$-primitive.
Let $v$ denote the Perron vector of $A^T$.  As $v\in
\fn{int}K^*,\, C\, := \{ x\in K: \langle x,v \rangle\ = 1 \}$ is a
complete (and hence compact) cross-section of $K$ and indeed it is
a polytope with $m$ extreme points. Replacing the extreme vectors
$x_1, \ldots, x_m$ of $K$ by suitable positive multiples, we may
assume that $x_1, \ldots, x_m$ are precisely all the extreme
points of $C$.  It is clear that $AC \subseteq C$, as $A\in
\pi(K)$ and $\rho(A) = 1$. Since the digraph $({\cal E},{\cal P})$
is given by Figure $1$, we have $Ax_j = x_{j+1}$ for $j = 1,
\ldots, m-1$ and also $Ax_m = (1-c)x_1+cx_2$ for some $c\in
(0,1)$.  The latter condition can be rewritten as
$(A^m-cA-(1-c)I_n)x_1 = 0$. It is clear that the $A$-invariant
subspace of ${\Bbb R}^n$ generated by $x_1$ is ${\Bbb R}^n$
itself; so $A$ is non-derogatory and also it follows that
$t^m-ct-(1-c)$ is an annihilating polynomial for $A$. Therefore,
$A$ has an annihilating polynomial of the desired form.

When $({\cal E},{\cal P}(A,K))$ is given by Figure $2$, we proceed
in a similar way. In this case we have
$$Ax_i = x_{i+1} \mbox{ for } i = 2, \ldots, m-1,$$
and
$$Ax_m = (1-a)x_1+ax_2 \mbox{ and } Ax_1 = (1-b)x_2+bx_3$$
for some $a,b\in (0,1)$.  Then after a little calculation we
obtain
$$[A^m-\left((1-a)b+a\right)A-(1-a)(1-b)I]x_2 = 0.$$
Since the $A$-invariant subspace of ${\Bbb R}^n$ generated by
$x_2$ is ${\Bbb R}^n$ itself, it follows that $A$ is
non-derogatory and $t^m-((1-a)b+a)t-(1-a)(1-b)$ is an annihilating
polynomial for $A$.  The latter polynomial can be rewritten as
$t^m-ct-(1-c)$, where $c = (1-a)b+a \in (0,1)$.

The uniqueness of the annihilating polynomial for $A$ of the
desired form is obvious, because $\{ A, I_n\}$ is a linearly
independent set.

(ii)  Note that for any $0 \ne x \in K$ and $j = 0,1, \ldots,
\gamma(A,x)-1, \gamma(A,x) = \gamma(A,A^jx)+j$.

First, consider the case when the digraph $({\cal E},{\cal P})$ is
given by Figure $1$.  Since $A^{j-1}x_1$ is a positive multiple of
$x_j$ for $j = 2, \ldots, m$, we have $\gamma(A,x_1) > m$ and
$$\gamma(A,x_j) = \gamma(A,A^{j-1}x_1) = \gamma(A,x_1)-j+1$$
for $j = 2, \ldots, m$; hence
$$\gamma(A) = \fn{max}_{1\le j \le
m}\gamma(A,x_j) = \gamma(A,x_1).$$

When the digraph $({\cal E},{\cal P})$ is given by Figure $2$, we
apply a similar but slightly more elaborate argument.  For $j = 3,
\ldots, m$, we have
 $$\gamma(A,x_2) = \gamma(A,A^{j-2}x_2)+j-2 =
\gamma(A,x_j)+j-2;$$ hence $\gamma(A,x_2) > \gamma(A,x_j)$ for
each such $j$.  A little calculation shows that $A^mx_2$ is a
positive linear combination of $x_2$ and $x_3$.  But $Ax_1$ is
also a positive linear combination of $x_2$ and $x_3$, hence
$\Phi(Ax_1) = \Phi(A^mx_2)$.  So we have
$$\gamma(A,x_2) = \gamma(A,A^mx_2)+m = \gamma(A,Ax_1)+m =
\gamma(A,x_1)-1+m,$$ which implies $\gamma(A,x_2) >
\gamma(A,x_1)$.  Therefore, we have $\gamma(A) = \gamma(A,x_2)$.

(iii) In the following argument, unless specified otherwise, we
assume that the digraph $({\cal E},{\cal P})$ is given by Figure
$1$ or Figure $2$. First, we show that each of the extreme vectors
$x_1, \ldots, x_m$, except possibly $x_2$, is involved in at least
one relation on $\fn{Ext}K$. For the purpose, it suffices to show
that $x_3$ is involved in one such relation; because, by applying
$A$ or its positive powers to a relation on $\fn{Ext}K$ involving
$x_3$, we can obtain for each of the vectors $x_4, \ldots,
x_{m-1}, x_m, x_1$ a relation that involves the vector. Suppose
that $x_3$ is not involved in any (nonzero) relation on
$\fn{Ext}K$. Take any relation $S$ on $\fn{Ext}K$; as $K$ is
non-simplicial, such relation certainly exists. Note that, since
$x_3$ does not appear in $S$, $x_4$ (and also $x_3$) cannot appear
in the (necessarily nonzero) relation obtained from $S$ by
applying $A$.  Similarly, the vectors $x_3,x_4,x_5$ all do not
appear in the relation obtained from $S$ by applying $A^2$.
Continuing the argument, we can show that the only vectors that
can appear in the nonzero relation obtained from $S$ by applying
$A^{m-3}$ are $x_1,x_2$. This contradicts the hypothesis that
$x_1,x_2$ are distinct nonzero extreme vectors of $K$.

Next, we note that if $x_2$ is not involved in any relation on
$\fn{Ext}K$ and if $({\cal E},{\cal P})$ is given by Figure $1$
then, by applying $A$ repeatedly to a nonzero relation on
$\fn{Ext}K$ sufficiently many times, we would obtain a nonzero
relation on $\fn{Ext}K$ that involves less than four vectors,
which is a contradiction.  So if $x_2$ is not involved in any
relation on $\fn{Ext}K$ then $({\cal E},{\cal P})$ must be given
by Figure $2$.

Now we contend that $K$ is either indecomposable or is the direct
sum of a ray and an indecomposable cone.  By what we have done
above, each of the extreme vectors $x_1, \ldots, x_m$, except
possibly $x_2$, belongs to an indecomposable summand of $K$ that
is not a ray.  Let $K_1$ be the indecomposable summand of $K$ that
contains $x_m$.  To establish our contention, it remains to show
that $K$ has no indecomposable summand which is not a ray and is
different from $K_1$.  It is known that every indecomposable
polyhedral cone of dimension greater than $1$ has a full relation
for its extreme vectors (see \cite{F--P}, p.37, (2.14)).  So it
suffices to show that there is no relation on $\fn{Ext}K$ that
involves vectors not belonging to $K_1$. Assume to the contrary
that there are such relations. Let $T_0$ be one such shortest
relation (i.e., one having the minimum number of terms).   Note
that, since $x_m$ is not involved in $T_0$, the relation obtained
from $T_0$ by applying $A$
 has the same number of terms as $T_0$, unless $x_1$
appears in $T_0$ but $x_2$ does not and $({\cal E},{\cal P})$ is
given by Figure $2$, in which case the said relation may have one
term more than $T_0$.  Suppose that the extreme vectors that
appear in $T_0$ are $x_{k_1}, x_{k_2}, \ldots, x_{k_s}$, where $1
\le k_1 < k_2 < \cdots < k_s \le m-1$. It is readily seen that the
relations obtained from $T_0$ by applying $A^i$ for $i = 1,
\ldots, m-k_s$ all have the same number of terms, as $x_1$ and
$x_m$ are both not involved in the first $m-k_s-1$ of these
relations.

 Let $q$ denote the
least positive integer such that $x_{m-q} \notin K_1$.  Certainly,
we have $k_s \le m-q$.  Denote by $\tilde{T_0}$ the relation
obtained from $T_0$ by applying $A^{m-q-k_s}$.  (If $k_s = m-q$,
we take $\tilde{T_0}$ to be $T_0$.)  Then $\tilde{T_0}$ either has
the same number of terms as $T_0$ or has one term more.  Note that
now $x_{m-q}$ is involved in $\tilde{T_0}$ and, by our choice of
$q$, $x_{m-q}\notin K_1$.  If $\tilde{T_0}$ involves also extreme
vectors of $K_1$, then $\tilde{T_0}$ splits and we would obtain a
relation for extreme vectors not belonging to $K_1$, which has at
least four terms fewer than that of $\tilde{T_0}$ and hence is a
relation shorter than the shortest relation $T_0$, which is a
contradiction. So we assume that all the vectors appearing in
$\tilde{T_0}$ do not belong to $K_1$ (and in fact they all lie in
the same indecomposable summand of $K$).

For $j = 1, 2, \ldots$, let $T_j$ denote the relation obtained
from $\tilde{T_0}$ by applying $A^j$.  By considering the cases
when $\tilde{T_0} = T_0$ and $\tilde{T_0} \ne T_0$ separately, one
readily sees that relation $T_1$ has at most one term more than
$T_0$.  Also, $T_1$ involves $x_{m-q+1}$ which, by the definition
of $q$, belongs to $K_1$. If $T_1$ involves also extreme vectors
not belonging to $K_1$, then $T_1$ splits and we would obtain a
contradiction.  So $T_1$ is a relation on $\fn{Ext}K_1$ and
$x_{\tilde{k}_1+1}, \ldots, x_{\tilde{k}_s+1}$ all belong to
$K_1$, where $\tilde{k}_j = (m-q-k_s)+k_j$.  By the same argument
we may assume that for $j = 1, \ldots, q, T_j$ is a relation on
$\fn{Ext}K_1$ with at most one term more than $T_0$.  So, we have
$x_{\tilde{k}_j+r} \in K_1$ for $r = 1, \ldots, q$ and $j = 1,
\ldots, s$.  Note that $x_m$ is involved in $T_q$ but $x_1$ is not
(as $x_m$ is not involved in $T_{q-1}$). So $x_1, x_2$ are both
involved in $T_{q+1}$ and lie on the same side of it.  As a
consequence, $x_2,x_3$ are both involved in $T_{q+2}$, $x_3,x_4$
are both involved in $T_{q+3}$, and so forth. Clearly, $T_{q+1}$
has one term more than $T_q$ and hence at most two terms more than
$T_0$.

If $T_{q+1}$ involves vectors not belonging to the same
indecomposable summand of $K$, then the relation splits and the
minimality of $T_0$ would be violated.  So we assume that
$T_{q+1}$ is a relation on $\fn{Ext}K_2$ --- here $K_2$ may or may
not be the same as $K_1$.  Note that now we have $x_1,x_2,
x_{\tilde{k}_j+q+1}\in K_2$ for $j = 1, \ldots, s-1$.  Let $l$
denote the smallest positive integer such that at least one of the
vectors $A^lx_2, A^lx_{\tilde{k}_j+q+1}, j = 1, \ldots, s-1$ does
not belong to $K_2$.   It is readily seen that $A^lx_2$ is always
a positive multiple of $x_{2+l}$ and $A^lx_1$ is a positive
multiple of $x_{1+l}$ or a positive linear combination of
$x_{1+l}$ and $x_{2+l}$, depending on whether $({\cal E},{\cal
P})$ is given by Figure $1$ or Figure $2$.  If $A^lx_2\notin K_2$,
then, from the definition of $l$, necessarily we have
$x_{2+l}\notin K_2$ and $x_{1+l}\in K_2$.  On the other hand, if
$A^lx_2\in K_2$, then there must exist $j = 1, \ldots, s-1$ such
that $A^lx_{\tilde{k}_j+q+1} \notin K_2$.  In any case, the
relation obtained from $T_{q+1}$ by applying $A^l$ involves at
least one extreme vector in $K_2$ and at least one extreme vector
not in $K_2$. Then the relation splits and we would obtain a
relation for extreme vectors not belonging to $K_1$, which is
shorter than the shortest such relation, which is a contradiction.
(A cautious reader may wonder whether there is a positive integer
$r < l$ such that $A^rx_{\tilde{k}_{s-1}+q+1} = x_m$.  If this is
so, then at one of the steps in applying $A$ $l$ times to relation
$T_{q+1}$ the number of terms in the relation is increased by one.
However, it is possible to show that such positive integer $r$
does not exist by analyzing carefully the cases
$\tilde{k}_{s-1}+q+1 < \tilde{k}_s$ and $\tilde{k}_{s-1}+q+1 =
\tilde{k}_s$ separately.)

It remains to show that if $({\cal E},{\cal P})$ is given by
Figure $2$ and $K$ is the direct sum of the ray $\fn{pos}\{ x_2
\}$ and the indecomposable polyhedral cone $\fn{pos}\{
x_1,x_3,x_4, \ldots, x_m \}$, then $m$ is odd and the latter cone
is an indecomposable minimal cone with a balanced relation for its
extreme vectors.  Note that the assumption that $x_2$ does not
appear in any relation on $\fn{Ext}K$ guarantees that every
relation obtained from a shortest relation on $\fn{Ext}K$ by
applying $A$ or its positive powers is still a shortest relation.
We contend every shortest relation involves each of the vectors
$x_1,x_3, \ldots, x_m$. Assume that the contrary holds. Take a
shortest relation $R$. Since $R$ has at least four terms, one of
the vectors $x_3, x_4, \ldots, x_m$ must appear in $R$.  On the
other hand, $R$ cannot involve all of these vectors; otherwise,
$x_1$ does not appear in $R$, and so the relation obtained from
$R$ by applying $A$ involves the vector $x_2$, which is a
contradiction. Thus we can find an $i, 4 \le i \le m$, such that
$x_i$ appears in $R$ and $x_{i-1}$ does not or the other way
round.  Then the relation obtained from $R$ by applying
$A^{m-i+1}$ involves one of the vectors $x_m,x_1$ but not both,
and so the relation obtained from $R$ by applying $A^{m-i+2}$ must
involve the vector $x_2$, which is a contradiction.  This proves
our contention.  Since every shortest relation on $\{ x_1,x_3,
\ldots, x_m \}$ is a full relation, it is clear that any two
relations on the latter set are multiples of each other; else, by
subtracting an appropriate multiple of one relation from another
we would obtain a shorter nonzero relation. This proves that the
cone $\fn{pos}\{ x_1,x_3, \ldots, x_m \}$ is minimal. Let $R$
denote the relation on $\fn{Ext}K$. Since $x_2$ does not appear in
any relation on $\fn{Ext}K, x_m, x_1$ must appear on opposite
sides of $R$.  So $x_1,x_3$ also appear on opposite sides of the
relation obtained from $R$ by applying $A$ and hence on opposite
sides of relation $R$. Continuing the argument, we infer that for
$j = 3, \ldots, m-1, x_j$ and $x_{j+1}$ lie on opposite sides of
$R$.  It follows that $m$ is odd and $R$ has the same number of
terms on its two sides, i.e., $R$ is a balanced relation.
\hfill$\blacksquare$

It is easy to show the following\,:

\begin{remark}\rm\label{remark3}
For any real numbers $p,l$ with $p \ge 3$, we have
$$(p-1)(l-2)+2 \le (p-1)(l-1),$$ where
the inequality becomes equality if and only if $p = 3$.
\end{remark}

Part (iii) and (iv) of Theorem \ref{theorem1}, our next result,
are not needed in the rest of the paper.  They are included for
the sake of completeness.

\begin{theorem}\label{theorem1}
Let $K\in {\cal P}(m,n)$, where $m \ge 4$, and let $A$ be a
$K$-primitive matrix. Then\,{\rm:}
\begin{enumerate}
\item[{\rm(i)}]  $\gamma(A) \le (m_A-1)(m-1)+1$, where the
equality holds only if the digraph $({\cal E},{\cal P}(A,K))$ is
 given by Figure $1$, in which
case $\gamma(A) = (n-1)(m-1)+1$. \item[{\rm(ii)}] $\gamma(A) =
(m_A-1)(m-1)$ only if either $({\cal E},{\cal P}(A,K))$ is given
by Figure $1$ or Figure $2$, in which case $\gamma(A) =
(n-1)(m-1)$, or $m_A = 3$. \item[{\rm(iii)}] $\gamma(A) =
(m_A-1)(m-2)+2$ only if $m_A \ge 3$ and either $({\cal E},{\cal
P}(A,K))$ is given by Figure $1$, Figure $2$, Figure $3$, Figure
$4$ or Figure $5$, or $({\cal E},{\cal P}(A,K))$ is obtained from
Figure $5$ by deleting any one or two of
the three arcs $(\Phi(x_{m-1}),\Phi(x_1))$,\\
\noindent$(\Phi(x_m),\Phi(x_1))$ and $(\Phi(x_m),\Phi(x_2))$, or
from Figure $3$ with $m=4$ by adding the arc
$(\Phi(x_3),\Phi(x_1))$ or substituting it for the arc
$(\Phi(x_4),\Phi(x_1))$. \item[{\rm(iv)}] If $m_A = 2$ or $({\cal
E},{\cal P}(A,K))$ is not given by Figures $1$-- $5$, nor is
derived from Figure $5$ or from Figure $3$ {\rm(}with $m =
4${\rm)} in the way as described in part {\rm(iii)}, then
 $$\gamma(A) \le
(m_A-1)(m-2)+1.$$
\end{enumerate}
\end{theorem}

\newsavebox{\figthree}
\savebox{\figthree}{ \put(-145,100){$x_1$} \put(-142,94){\vector(0,-1){21}}
\put(-145,62){$x_2$}
\put(-134.5,56){\vector(1,-1){19}}
\put(-180,27){$x_m$} \put(-162,30){\vector(1,0){40}}
\put(-114,27){$x_3$} \put(-99,30){\vector(1,0){25}}
\put(-68,27){$x_4$} \put(-163,-5){\vector(-1,3){8}}
\put(-67,20){\vector(-1,-3){8.5}} \put(-169,-16){$x_{m-1}$}
\put(-83,-16){$x_5$} \put(-134,-40){\vector(-1,1){17}}
\put(-125,-43){$\cdots$} \put(-85,-22){\vector(-1,-1){17}} }
\[
\put(-30,0){\usebox{\figthree}} \put(-203,37){\vector(1,2){28}} \put(-160,62){\vector(2,-1){58}}
\put(-170,-65){\rm Figure 3.}
\put(121,0){\usebox{\figthree}} \put(-26,56){\vector(-1,-1){20}}
\put(-90,45){\qbezier(38,-53)(-8,10)(50,45)} \put(-39,90){\vector(2,1){10}}
\put(-24,-65){\rm Figure 4.}
\put(270,0){\usebox{\figthree}} \put(97,37){\vector(1,2){28}} \put(103,36){\vector(1,1){20}}
\put(59,45){\qbezier(38,-53)(-8,10)(50,45)} \put(110,90){\vector(2,1){10}}
\put(125,-65){\rm Figure 5.}
\]

\noindent {\it Proof}.  When $m_A = 2$, by Lemma \ref{lemma2} we
have $\gamma(A) \le 2$. As $m \ge 4$, in this case, the inequality
$\gamma(A) \le (m_A-1)(m-2)+1$ is clearly satisfied and none of
the equalities $\gamma(A) = (m_A-1)(m-1)$ or $\gamma(A) =
(m_A-1)(m-2)+2$ can be attained. Hereafter, we assume that $m_A
\ge 3$.

As explained before, the length of the shortest circuit in $({\cal
E},{\cal P})$ is at most $m-1$.

(i) Since $A$ is $K$-primitive, $A$ is non-nilpotent.  So the
outdegree of each vertex of $({\cal E},{\cal P})$ is positive.
Consider any vertex $\Phi(x)$ of the digraph $({\cal E},{\cal
P})$.  It is clear that $\Phi(x)$ lies on or has access to a
circuit of length $l \le m-1$.  By Lemma \ref{lemma3.2} we have
\[
    \gamma(A,x) \le (m_A-2)l+m \le (m_A-2)(m-1)+m = (m_A-1)(m-1)+1.
\]
Since this is true for every nonzero extreme vector $x$ of $K$, the
inequality $\gamma(A) \le (m_A-1)(m-1)+1$ follows.

To establish the desired necessary condition for the inequality to
become equality, first we dispense with the case when the length
of the shortest circuit in $({\cal E},{\cal P})$ is less than or
equal to $m-2$. We contend that in this case every vertex of the
digraph lies on or has access to a circuit of length less than or
equal to $m-2$. Consider any vertex $\Phi(x)$ of the digraph. As
we have explained before, $\Phi(x)$ lies on or has access to some
circuit, say ${\cal C}$. Choose such a ${\cal C}$ of shortest
length, say length $l$.  By the definition of $l$, $C$ contains no
chords or loops (unless ${\cal C}$ is itself a loop).  If $l = m$
then $A$ is not primitive, contradiction.  If $l = m-1$ let $z$ be
the unique vertex of the digraph not on ${\cal C}$.  Then, there
is an access from ${\cal C}$ to $z$ and vice versa, or else $A$ is
not primitive.  Hence the graph is strongly connected, so
$\Phi(x)$ has access to a circuit of length $m-2$ or less,
contradicting $l = m-1$.  Hence, necessarily, $l\le m-2$.  This
proves our contention.  So, in this case, we have
\begin{equation}\label{eq?.2}
\gamma(A) \le (m_A-2)(m-2)+m = (m_A-1)(m-2)+2 \le (m_A-1)(m-1),
\end{equation}

\noindent where the first inequality holds by Lemma 3.2 (with $l
\le m-2$) and the second inequality follows from Remark 4.3 (with
$p = m_A$ and $l = m$).

When the length of the shortest circuit in $({\cal E},{\cal P})$
is $m-1$, by Lemma \ref{lemma3} the digraph is given by Figure $1$
or Figure $2$. If the digraph is given by Figure $2$, then each
vertex lies on a circuit of length $m-1$ and by Lemma \ref{lemma1}
we obtain $\gamma(A) \le (m_A-1)(m-1)$.  So when the equality
$\gamma(A) = (m_A-1)(m-1)+1$ holds, the digraph $({\cal E},{\cal
P})$ must be given by Figure $1$.  In that case, by Lemma
\ref{lemma4.2}(i), we have $m_A = n$ and hence $\gamma(A) =
(n-1)(m-1)+1$.

(ii) Suppose that the equality $\gamma(A) = (m_A-1)(m-1)$ holds.
The length of the shortest circuit in $({\cal E},{\cal P})$ is
either $m-1$ or less. In the former case, by Lemma \ref{lemma3}
the digraph is given by Figure 1 or Figure 2; then by Lemma
\ref{lemma4.2}(i), as $m_A = n$ the said equality becomes
$\gamma(A) = (n-1)(m-1)$.  In the latter case, by the proof of
part (i) the inequalities in (\ref{eq?.2}) both hold as equality,
and by Remark \ref{remark3} we have $m_A = 3$ and hence also
$\gamma(A) = 2(m-1)$.

(iii)  Suppose that $\gamma(A) = (m_A-1)(m-2)+2$.  If $({\cal
E},{\cal P})$ is given by Figure $1$ or Figure $2$ we are done, so
assume that this is not the case.  Note that it is not possible
that every vertex of $({\cal E},{\cal P})$ lies on or has access
to a circuit of length $\le m-3$ or is at a distance at most $1$
to a circuit of length $m-2$, because then by Lemma \ref{lemma3.2}
(with $l = m-3$) or by Lemma \ref{lemma1} (with $w = 1$ and $l =
m-2$) it will follow that $\gamma(A) \le (m_A-1)(m-2)+1$. It
remains to show that if $({\cal E},{\cal P})$ has a vertex which
is at a distance $2$ to a circuit of length $m-2$ and which does
not lie on or has access to a circuit of length $m-3$ or less, nor
is it at a distance at most $1$ to another circuit of length
$m-2$, then the digraph is given by Figure $3$, Figure $4$ or
Figure $5$, or is derived from them in the manner as described in
the theorem.

To treat the remaining case we assume that the digraph $({\cal
E},{\cal P})$ contains the circuit ${\cal C}: \Phi(x_3)
\rightarrow \Phi(x_4) \rightarrow \cdots \rightarrow \Phi(x_{m-1})
\rightarrow \Phi(x_m) \rightarrow \Phi(x_3)$ and also the path
$\Phi(x_1) \rightarrow \Phi(x_2) \rightarrow \Phi(x_3)$.  For $i =
2, 3, \ldots, m$, since $\Phi(x_i)$ is at a distance at most $1$
to the circuit ${\cal C}$, which is of length $m-2$, by Lemma 3.1
we have $\gamma(A,x_i) \le (m_A-1)(m-2)+1$.  This forces
$\gamma(A,x_1) = \gamma(A) = (m_A-1)(m-2)+2$, which, in turn,
implies that $\Phi(x_1)$ does not lie on or has access to a
circuit of length $m-3$ or less, nor is $\Phi(x_1)$ at a distance
at most $1$ to a circuit of length $m-2$.  Therefore, ${\cal C}$
does not contain any chords or loops.  Besides the arcs on the
circuit ${\cal C}$ and the above-mentioned path, $({\cal E},{\cal
P})$ certainly has other arcs.  We want to find out what possible
additional arcs there can be.

Note that there is at least one arc from a vertex of ${\cal C}$ to
either one of the vertices $\Phi(x_1)$ or $\Phi(x_2)$; else,
$A^{m-2}$ maps the extreme ray $\Phi(x_3)$ of $K$ onto itself,
which contradicts the $K$-primitivity of $A$.  Since $x_1$ is not
allowed to lie on a circuit of length $m-2$ or less, none of the
arcs $(\Phi(x_j),\Phi(x_1))$, for $j = 2, \ldots, m-2$, can be
present. Similarly, since $\Phi(x_1)$ is not allowed to be at a
distance $1$ to a circuit of length $m-2$ or less, the arcs
$(\Phi(x_j),\Phi(x_2))$, for $j = 3, \ldots, m-1$, also cannot be
present.  So $(\Phi(x_{m-1}),\Phi(x_1)), (\Phi(x_m),\Phi(x_1))$
and $(\Phi(x_m),\Phi(x_2))$ are the only possible arcs from a
vertex of ${\cal C}$ to either $\Phi(x_1)$ or $\Phi(x_2)$; also,
at least one of these three arcs is present.

There cannot exist an arc from $\Phi(x_1)$ to a vertex of ${\cal C}$, because
in the presence of any such arc the distance from $\Phi(x_1)$ to the circuit
 ${\cal C}$ becomes $1$.  Similarly, for $m \ge 6$, each of the arcs
 $(\Phi(x_2),\Phi(x_j)), j = 5, \ldots, m-1$, also cannot exist, because the arc $(\Phi(x_2),\Phi(x_j))$
and one of the arcs $(\Phi(x_{m-1}),\Phi(x_1)),
(\Phi(x_m),\Phi(x_1))$ and $(\Phi(x_m),\Phi(x_2))$ (which must be
present), together with some of the arcs in ${\cal C}$, form
either a circuit of length $m-3$ or less, which is at a distance
$1$ from $\Phi(x_1)$, or a circuit of length $m-2$ or less that
contains $\Phi(x_1)$, but this is not allowed.
 So $(\Phi(x_2),\Phi(x_4))$ and $(\Phi(x_2),\Phi(x_m))$ are the only
possible arcs from $\Phi(x_1)$ or $\Phi(x_2)$ to a vertex of
${\cal C}$ when $m \ge 6$.  The preceding argument does not cover
the cases when $m = 4 \mbox{ or }5$.  However, we have not ruled
out the possibility of the existence of the arcs
$(\Phi(x_2),\Phi(x_4))$ and $(\Phi(x_2),\Phi(x_m))$ in these
cases.

\newsavebox{\figp}
\savebox{\figp}{ \put(40,63.5){$x_1$} \put(40,17){$x_2$} \put(40,-29){$x_3$}
\put(43,58){\vector(0,-1){30}} \put(43,11.5){\vector(0,-1){30}}
\put(35,18){\qbezier(1,-40)(-30,0)(0,40)} \put(32.5,55){\vector(1,1){5}} }
\[
\put(-170,-75){\usebox{\figp} \put(95,-29){$x_4$} \put(55,-27){\vector(1,0){34}}
\put(53.5,13.5){\vector(1,-1){37}} \put(92.5,-19.5){\vector(-1,2){39}}
\put(73,-32){\qbezier(-19,0)(0,-15)(19,0)} \put(75,-39.5){\vector(-1,0){4}} }
\put(-132,-135){\rm Figure $3^{'}$}
\put(40,-75){\usebox{\figp} \put(95,-29){$x_4$} \put(55,-27){\vector(1,0){34}}
\put(91,-23){\vector(-1,1){37}} \put(92.5,-19.5){\vector(-1,2){39}}
\put(73,-32){\qbezier(-19,0)(0,-15)(19,0)} \put(75,-39.5){\vector(-1,0){4}} }
\put(80,-135){\rm Figure $5^{'}$}
\]

Consider the case when the arc $(\Phi(x_2),\Phi(x_4))$ is present.
If $m \ge 5$, then the last two of the three arcs
$(\Phi(x_m),\Phi(x_1)), (\Phi(x_m),\Phi(x_2))$ and
$(\Phi(x_{m-1}),\Phi(x_1))$ cannot be present, else $\Phi(x_1)$ is
at a distance at most 1 to a circuit of length $m-2$, which is not
allowed.  So in this case the arc $(\Phi(x_m),\Phi(x_1))$ must be
present and, furthermore, the arc $(\Phi(x_2),\Phi(x_m))$ also
cannot be present (otherwise, we have a circuit of length three
containing $\Phi(x_1)$). Therefore, the digraph $({\cal E},{\cal
P})$ is given by Figure 3. If $m = 4$, we find that the arc
$(\Phi(x_4),\Phi(x_2))$ cannot be present, but the arcs
$(\Phi(x_3),\Phi(x_1))$ and $(\Phi(x_4),\Phi(x_1))$ may be present
and, indeed, at least one of them must be present. So the digraph
$({\cal E},{\cal P})$ is given by Figure $3^{'}$ or is obtained
from it by deleting one of the arcs $(\Phi(x_3),\Phi(x_1)),
(\Phi(x_4),\Phi(x_1))$.  In other words, the digraph is given by
Figure 3 (with $m = 4$) or is derived from it in the manner as
described in the theorem.

Now suppose that the arc $(\Phi(x_2),\Phi(x_m))$ is present. Using
the same kind of argument as before, for $m \ge 5$, one readily
rules out the presence of the arcs $(\Phi(x_m),\Phi(x_1))$ and
$(\Phi(x_m),\Phi(x_2))$.  So in this case the arc
$(\Phi(x_{m-1}),\Phi(x_1))$ must be present.  Then we can show
that the arc $(\Phi(x_2),\Phi(x_4))$ cannot be present. Therefore,
the digraph $({\cal E},{\cal P})$ is given by Figure $4$.  For $m
= 4$, we are dealing with the situation when
$(\Phi(x_2),\Phi(x_4))$ is an arc, but this has already been
treated above (for arbitrary $m \ge 4$).

It remains to consider the case when the arcs
$(\Phi(x_2),\Phi(x_4))$ and $(\Phi(x_2),\Phi(x_m))$ are both
absent. Then the presence of any one, two or three of the arcs
$$(\Phi(x_{m-1}),\Phi(x_1)), (\Phi(x_m),\Phi(x_1)),
(\Phi(x_m),\Phi(x_2))$$ will produce only circuits of length at
least $m-1$, but that causes no problem.  Then the digraph $({\cal
E},{\cal P})$ is given by Figure $5$ (which becomes Figure $5^{'}$
when $m = 4$) or is obtained from it by deleting any one or two of
the above-mentioned three arcs.

(iv)  Now this is obvious.
\hfill$\blacksquare$\\

\begin{remark}\label{remarkZ}\rm Let $K\in {\cal P}(3,3)$, and let $A$ be a $K$-primitive
matrix.  Then $\gamma(A) \le 2m_A-1$, where the equality holds
only if $({\cal E},{\cal P}(A,K))$ is given by Figure $1$, in
which case $\gamma(A) = 5$.
\end{remark}

The preceding remark, in fact, says that part(i) of Theorem
\ref{theorem1} still holds when $m = n = 3$.  It holds by what is
known in the $3\times 3$ nonnegative matrix case. However, parts
(ii)--(iv) of Theorem \ref{theorem1} cannot be extended to the
case $m = n = 3$. This is mainly because in that case the equality
in part (ii) or (iii) can hold when $m_A = 2$.

By Lemma \ref{lemma3}, for any polyhedral cone $K$ with $m \ge 3$
extreme rays and any $K$-primitive matrix $A$, the length of the
shortest circuit in $({\cal E},{\cal P})$ is $m-2$ or less if and
only if the digraph $({\cal E},{\cal P})$ is not given by Figure
$1$ or Figure $2$ or by a digraph of order $3$ whose arc set
consists of all possible arcs between every pair of distinct
vertices.  The proof of Theorem \ref{theorem1} (i) shows that in
this case $(m_A-1)(m-2)+2$ is an upper bound for $\gamma(A)$. (The
case $m = n = 3$ can be treated separately.)  So we have the
following

\begin{remark}\label{remark extra}\rm  For any $K\in {\cal P}(m,n)$ and any $K$-primitive matrix $A$, if the digraph of
$({\cal E},{\cal P}(A,K))$ is not given by Figure $1$ or Figure
$2$ or by a digraph of order $3$ whose arc set consists of all
possible arcs between every pair of distinct vertices, then
$\gamma(A) \le (n-1)(m-2)+2$.
\end{remark}

In below we give another bound for $\gamma(A)$ in terms of $m_A,
m$ and $s$, where $s$ is the length of the shortest circuit in
$({\cal E},{\cal P})$.  Before we do that, we need to obtain a
general result on a digraph first.

\begin{remark}\label{remark4.5} \rm Let $D$ be a digraph on $m \ge 3$ vertices, each of which
has positive out-degree. If the length of the shortest circuit in
$D$ is greater than $\lfloor\frac{m-1}{2}\rfloor$, then every
vertex of $D$ lies on or has access to a circuit of $D$ of
shortest length.
\end{remark}

{\it Proof}. Since each vertex of $D$ has positive out-degree,
each vertex lies on or has access to a circuit. Denote by $s(D)$
the length of the shortest circuit in $D$. If there is a vertex
that does not lie on or has access to a circuit of length $s(D)$,
then such vertex must lie on or has access to a circuit, say
${\cal C}$, of length $s(D)+1$ or more.  It is clear that the
circuit ${\cal C}$ is vertex disjoint from every circuit of
shortest length. Consequently, we have $m \ge s(D)+(s(D)+1)$ or
$\lfloor\frac{m-1}{2}\rfloor \ge s(D)$, which is a contradiction.
\hfill$\blacksquare$\\

By Lemma \ref{lemma3.2} and Remark \ref{remark4.5} we have

\begin{remark}\label{remarkA}\rm  Let $K\in {\cal P}(m,n)$ and let $A$ be a $K$-primitive matrix.  Let $s$ be the length
of the shortest circuit in $({\cal E},{\cal P}(A,K))$.  If $s >
\lfloor\frac{m-1}{2}\rfloor$, then $\gamma(A) \le s(m_A-2)+m.$
\end{remark}

It is interesting to note that the digraphs given by Figure $3$,
Figure $4$ and Figure $5$ are all primitive, like the digraphs
given by Figure $1$ and Figure $2$.  Moreover, if $A$ is a
$K$-primitive matrix such that $({\cal E},{\cal P})$ is given by
Figure $3$, Figure $4$ or Figure $5$ then $A$ is necessarily
nonsingular --- this can be proved using the argument given in the
proof of Lemma \ref{lemma4.2} (i).  However, the digraph obtained
from Figure $5'$ (i.e., Figure $5$ with $m = 4$) by removing the
arcs $(x_4,x_2)$ and $(x_3,x_1)$ is strongly connected but not
primitive, whereas the one obtained from Figure $3'$ (i.e., Figure
$3$ with $m = 4$) by removing the arcs $(x_4,x_1)$ and $(x_3,x_1)$
is not even strongly connected.  Also, $A$ is singular if it is a
$K$-primitive matrix such that its digraph $({\cal E},{\cal P})$
is derived from Figure $5$ by deleting the arcs
$(\Phi(x_m),\Phi(x_1))$ and $(\Phi(x_m),\Phi(x_2))$.

\begin{coro}\label{coro4}
For any $K\in {\cal P}(m,n)$ with $m = n+k$, we have $\gamma(K)
\le (n-1)(m-1)+1 = m^2-(k+2)m+k+2$. The equality holds only if
there exists a $K$-primitive matrix $A$ such that the digraph
$({\cal E},{\cal P}(A,K))$ is given by Figure $1$.
\end{coro}

{\it Proof}.  Follows from part (i) of Theorem \ref{theorem1}\,(i)
and Remark \ref{remarkZ}. \hfill$\blacksquare$\\

By Corollary \ref{coro4} the answer to Kirkland's conjecture
mentioned at the beginning of Section $1$ is in the affirmative.

\begin{coro}\label{coro5}  For any positive integer $m \ge 3$,
$$ \max \{ \gamma(K): K \mbox{ is a polyhedral cone with } m \mbox{ extreme rays } \} = m^2-2m+2.$$
\end{coro}
{\it Proof}.  Let $K$ be an $n$-dimensional polyhedral cone with $m$ extreme rays.  Since $n
\le m$, by Corollary \ref{coro4}, $\gamma(K) \le (m-1)^2+1$.  So we have
$$ \max \{ \gamma(K): K \mbox{ is a polyhedral cone with } m \mbox{ extreme rays} \} \le m^2-2m+2.$$
On the other hand, by Wielandt's bound we also have $\gamma({\Bbb
R}^m_+) = m^2-2m+2$.  Hence, the desired equality follows.
\hfill$\blacksquare$\\

We would like to emphasize that in Corollary \ref{coro5} the
number of extreme rays (i.e., $m$) for the polyhedral cones $K$
under consideration is fixed but there is no restriction on their
dimensions (i.e., $n$).

Hereafter, we call a polyhedral cone $K_0\in {\cal P}(m,n)$ an
{\it exp-maximal cone} if $\gamma(K_0) = \fn{max}\{ \gamma(K):
K\in {\cal P}(m,n) \}$.

\setcounter{equation}{0}
\section{The minimal cones case}

 For a non-simplicial polyhedral cone $K$ and a $K$-primitive
matrix $A$, if the digraph $({\cal E},{\cal P})$ is given by
Figure $1$ or Figure $2$, then Lemma \ref{lemma4.2} gives us some
information about $K$ and $A$.  The first lemma of this section
shows that if, in addition, $K$ is a minimal cone then the
relation for the extreme vectors of $K$ and the action of $A$ on
the extreme vectors can be described completely.

\begin{lemma}\label{lemma 5.1}
Let $K\in {\cal P}(n+1,n), n \ge 3$. Let $A$ be a $K$-nonnegative
matrix.
\begin{enumerate}
\item[{\rm (i)}]  Suppose that $({\cal E},{\cal P}(A,K))$ is given
by Figure $1$ {\rm(}with $m = n+1$, and $K$ is necessarily
indecomposable{\rm)}. If $n$ is odd then, after normalization
{\rm(}on the extreme vectors of $K$ and on $A${\rm)}, the relation
on $\fn{Ext}K$ and the matrix $A$ are given respectively by
{\rm(\ref{rel2})} and {\rm(\ref{eq2})}\,{\rm :}
\begin{eqnarray}
 x_1+x_3+\cdots+x_{m-3}+x_{m-1} &=& x_2+x_4+ \cdots + x_{m-2}+x_m.  \label{rel2}
\end{eqnarray}

\begin{eqnarray}
 Ax_1 &=&  (1+\alpha)x_2,  \nonumber\\
 Ax_i &=& x_{i+1}\mbox{ for } i = 2, 3, \ldots, m-1, \label{eq2}\\
 Ax_m &=& x_1+\alpha x_2, \nonumber
\end{eqnarray}
where $\alpha$ is some positive scalar.  If $n$ is even then,
after normalization, the relation on $\fn{Ext}K$ and the matrix
$A$ are given respectively by {\rm(\ref{rel3})} and
{\rm(\ref{eq3})}\,{\rm :}
\begin{eqnarray}
 x_1+x_2+x_4+\cdots+x_{m-3}+x_{m-1} &=& x_3+x_5+ \cdots + x_{m-2}+x_m. \label{rel3}
\end{eqnarray}

\begin{eqnarray}
Ax_1 &=& \alpha x_2,  \nonumber \\
Ax_i &=& x_{i+1}\mbox{ for } i = 2, 3, \ldots, m-1, \label{eq3}\\
Ax_m &=& x_1+(1+\alpha)x_2, \nonumber
\end{eqnarray}
where $\alpha > 0$.

 \item[{\rm(ii)}] Suppose $K$ is indecomposable
and $({\cal E},{\cal P}(A,K))$ is given by Figure $2$.  If $n$ is
even then, after normalization, the relation on $\fn{Ext}K$ and
the matrix $A$ are given respectively by {\rm(\ref{rel3})} and
{\rm(\ref{eq4})},  or by {\rm(\ref{rel5})} and
{\rm(\ref{eq5})}\,{\rm:}

\begin{eqnarray}
 Ax_1 &=& \alpha x_2+(1-\beta)x_3,  \nonumber \\
 Ax_2 &=& \beta x_3, \nonumber \\
 Ax_i &=& x_{i+1}, \mbox{ for } i = 3, \ldots, m-1 \label{eq4},\\
 Ax_m &=& x_1+(1+\alpha)x_2, \nonumber
\end{eqnarray}
where $\alpha > 0, 0 < \beta < 1$.

\begin{eqnarray}
 &\!\!x_2+x_3+x_5+ \cdots +x_{m-2}+x_m = x_1+x_4+x_6+ \cdots +x_{m-3}+x_{m-1}.& ~\label{rel5}
\end{eqnarray}

\begin{eqnarray}
   Ax_1 &=& (1+\alpha)x_2+(1+\beta)x_3, \nonumber \\
   Ax_2 &=& \beta x_3, \nonumber \\
   Ax_i &=& x_{i+1} \mbox{ for } i = 3, \ldots, m-1, \label{eq5} \\
   Ax_m &=& x_1+\alpha x_2, \nonumber
\end{eqnarray}
where $\alpha,\beta > 0$.  If $n$ is odd then, after
normalization, the relation on $\fn{Ext}K$ and the matrix $A$ are
given respectively by {\rm(\ref{rel2})} and {\rm(\ref{eq6})}\,{\rm
:}

\begin{eqnarray}
 Ax_1 &=& (1+\alpha) x_2+\beta x_3, \nonumber \\
 Ax_2 &=& (1+\beta) x_3, \nonumber \\
 Ax_i &=& x_{i+1}, i = 3, \ldots, m-1, \label{eq6} \\
 Ax_m &=& x_1+\alpha x_2, \nonumber
\end{eqnarray}
where $\alpha, \beta > 0$.

\item[{\rm (iii)}]  Suppose $K$ is decomposable and $({\cal
E},{\cal P}(A,K))$ is given by Figure $2$. Then, after
normalization, the relation on $\fn{Ext}K$ and the matrix $A$ are
given respectively by {\rm (\ref{rel1})} and {\rm
(\ref{eq1})}\,{\rm :}

\begin{eqnarray}
 x_1+x_4+x_6+ \cdots + x_{m-3}+ x_{m-1} = x_3+x_5+ \cdots +x_{m-2}+ x_m. \label{rel1}
\end{eqnarray}

\begin{eqnarray}
   Ax_1 &=& \alpha x_2+x_3,  \nonumber \\
   Ax_2 &=& \beta x_3, \nonumber \\
   Ax_i &=& x_{i+1} \mbox{ for } i = 3, \ldots, m-1, \label{eq1}\\
   Ax_m &=& x_1+\alpha x_2, \nonumber
\end{eqnarray}
where $\alpha, \beta > 0$.
\end{enumerate}
\end{lemma}

{\it Proof}.  (i) Since $K$ is non-simplicial and $({\cal E},{\cal
P})$ is given by Figure $1$, by Lemma \ref{lemma4.2}(iii) $K$ is
indecomposable.  So the (essentially unique) relation on
$\fn{Ext}K$, which we denote by $R$, is full.  We contend that
$x_m,x_1$ lie on different sides of $R$. Suppose not. Then
$x_1,x_2$ lie on the same side of the relation obtained from $R$
by applying $A$. But the latter relation, which is nonzero (as the
coefficients of $x_1,x_2$ are both nonzero), is just a multiple of
$R$, so $x_1,x_2$ also lie on the same side of relation $R$. By
applying $A$ to the latter relation, we deduce that $x_2, x_3$
also lie on the same side of relation $R$.  Continuing the
argument, we can then show that $x_1, x_2, \ldots, x_m$ all lie on
the same side of $R$, which is impossible, as $K$ is a pointed
cone. This proves our contention. The same argument, in fact, also
shows that for $j = 2, 3, \ldots, m-1$, $x_j, x_{j+1}$ lie on
different sides of $R$. Now a simple parity count shows that $x_2,
x_m$ lie on the same side or opposite sides of $R$, depending on
whether $n$ is odd or even. So when $n$ is odd (i.e., $m$ is
even), after normalizing the extreme vectors of $K$ we may assume
that relation $R$ is given by (\ref{rel2}).

As the digraph $({\cal E},{\cal P})$ is given by Figure $1$, we
have
$$Ax_1 = \beta x_2, Ax_i = \lambda_{i+1}x_{i+1} \mbox{ for }
i = 2, \ldots, m-1 \mbox{ and } Ax_m = \lambda_1x_1+\alpha x_2,$$
where $\alpha, \beta$ and $\lambda_1, \lambda_3, \lambda_4,
\ldots, \lambda_m$ are some positive numbers.  Substituting the
values of the $Ax_i$'s into the relation obtained from
(\ref{rel2}) by applying $A$, we obtain the relation:
$$\beta x_2+\lambda_4x_4+\lambda_6x_6+ \cdots + \lambda_mx_m
= \lambda_3x_3+\lambda_5x_5+ \cdots +
\lambda_{m-1}x_{m-1}+(\lambda_1x_1+\alpha x_2).$$ But relation
(\ref{rel2}) and the above relation are positive multiples of each
other, so it follows that we have  $\lambda_1 = \lambda_3 =
\lambda_4 = \cdots = \lambda_m$ and $\beta = \lambda + \alpha$,
where we use $\lambda$ to denote the common value of the
$\lambda_j$'s.
 Replacing $A$ by a positive multiple, we may assume that
 $\lambda = 1$.  Then $A$ is
 given by (\ref{eq2}).

When $n$ is even, we can show in a similar way that the relation
on $\fn{Ext}K$ and the matrix $A$ are given by (\ref{rel3}) and
(\ref{eq3}) respectively.

(ii)  We consider the case when $n$ is even first.  By the same
kind of argument that we have used for part (i) we can show that
for $j = 3, \ldots, m$, the vectors $x_j, x_{j+1}$ lie on
different sides of the relation on $\fn{Ext}K$ (where $x_{m+1}$ is
taken to
 be $x_1$).  Hence, the vectors $x_3, x_5, \ldots, x_m$ lie on one side of the relation and the
vectors $x_1, x_4, x_6, \ldots, x_{m-1}$ lie on the other side. As
for the vector $x_2$ it can be on either side. If $x_2$ is on the
same side as $x_1$ then, after normalizing the extreme vectors of
$K$, we may assume that the relation on $\fn{Ext}K$ is given by
(\ref{rel3}); if $x_2$ lies on the side opposite to $x_1$, we may
assume that the relation is given by (\ref{rel5}).

We treat the subcase when the relation is given by (\ref{rel3}),
the argument for the remaining subcase being similar. Since the
digraph $({\cal E},{\cal P})$ is given by Figure $2$, we have
$$Ax_1 = \alpha x_2+\gamma x_3, Ax_2 = \beta x_3, Ax_i = \lambda_{i+1}x_{i+1} \mbox{ for }
i = 3, \ldots, m-1 \mbox{ and } Ax_m = \lambda_1x_1+\delta x_2,$$
where $\alpha, \beta, \delta, \gamma$ and $\lambda_1, \lambda_4,
\ldots, \lambda_m$ are some positive numbers.  Applying $A$ to
relation (\ref{rel3}), we obtain the relation:
$$\lambda_4x_4+\lambda_6x_6+ \cdots +
\lambda_{m-1}x_{m-1}+(\lambda_1x_1+\delta x_2) = (\alpha
x_2+\gamma x_3)+\beta x_3+\lambda_5x_5+ \cdots +
\lambda_{m-2}x_{m-2}+\lambda_mx_m.$$ But relation (\ref{rel3}) and
the above relation are positive multiples of each other, it
follows that we have $\lambda_4 = \lambda_5 = \cdots = \lambda_m =
\lambda$, say, and $\lambda_1 = \lambda, \delta = \lambda+\alpha$
and $\gamma + \beta = \lambda$.   Replacing $A$ by a positive
multiple, we may assume that $\lambda = 1$.  Then $A$ is
 given by equation (\ref{eq4}).

Now we consider the case when $n$ is odd.  Again, we can show that
for $j = 3, \ldots, m$, the vectors $x_j, x_{j+1}$ lie on
different sides of the relation on $\fn{Ext}K$ (where $x_{m+1}$ is
taken to be $x_1$). Hence, $x_1,x_3,x_5, \ldots, x_{m-3},x_{m-1}$
lie on one side of the unique
 relation and $x_4,x_6, \ldots, x_{m-2},x_m$ lie on the other side.  If $x_1,x_2$ lie on the
 same side of the relation then, since $x_1, x_3$ also lie on the same side, the same is
 true for the pair $x_2,x_3$.  Then by applying $A$ we find that $x_3,x_4$ also lie
 on the same side of the relation, which contradicts what we have observed above.
 So $x_2$ lies on the same side as $x_4, x_6, \ldots, x_m$, and after normalization we may assume that
 the unique relation is given by (\ref{rel2}).  In a similar way as before we can
 also show that after normalization $A$ is given by
 (\ref{eq6}).

(iii)  Suppose $K$ is decomposable and $({\cal E},{\cal P}(A,K))$
is given by Figure $2$.  By Lemma \ref{lemma4.2}(iii), $m$ is odd
and $K$ is the direct sum of a ray and an indecomposable minimal
cone with a balanced relation for its extreme vectors.  The last
part of the proof for Lemma \ref{lemma4.2}(iii) shows that in the
relation on $\fn{Ext}(K)$ the vectors $x_1, x_4, x_6, \ldots,
x_{m-1}$ lie on one side and the vectors $x_3, x_5, \ldots, x_m$
lie on the other side. After normalizing the extreme vectors of
$K$, we may assume that the relation on $\fn{Ext}K$ is given by
relation (\ref{rel1}).  By the same kind of argument as before, we
can also show that $A$, after normalization, is given by equation
(\ref{eq1}). \hfill$\blacksquare$

\begin{lemma}\label{lemma?}  Let $K\in {\cal P}(m,n)$ be a minimal cone with extreme vectors
 $x_1, \ldots, x_m$ {\rm(}where
$m = n+1${\rm)}, and let $A$ be a $K$-nonnegative matrix.  Let the
relations {\rm (\ref{rel2}), (\ref{rel3}), (\ref{rel5}), (\ref{rel1})}
on the extreme vectors of $K$ and the equations {\rm
(\ref{eq2}), (\ref{eq3}), (\ref{eq4}), (\ref{eq5}), (\ref{eq6}), (\ref{eq1})}
on the action of $A$ be as given in Lemma \ref{lemma 5.1}.
 Then{\rm:}
\begin{enumerate}
\item[{\rm (i)}] $({\cal E},{\cal P}(A,K))$ is given by Figure $1$
{\rm(}and $K$ is indecomposable{\rm)} if and only if after
normalization the relation on $\fn{Ext}K$ and the matrix $A$ are
given respectively by relation {\rm(\ref{rel2})} and equation
{\rm(\ref{eq2})}, in which case $n$ is odd and $\gamma(A) =
n^2-n+1$, or by relation {\rm(\ref{rel3})} and equation
{\rm(\ref{eq3})}, in which case $n$ is even and $\gamma(A) =
n^2-n$. \item[{\rm(ii)}] $({\cal E},{\cal P}(A,K))$ is given by
Figure $2$ and $K$ is indecomposable if and only if after
normalization the relation on $\fn{Ext}K$ and the matrix $A$ are
given respectively by relation {\rm(\ref{rel5})} and equation
{\rm(\ref{eq5})}, in which case $n$ is even and $\gamma(A) =
n^2-n$, or by relation {\rm(\ref{rel3})} and equation
{\rm(\ref{eq4})}, in which case $n$ is even and $\gamma(A) =
n^2-n-1$, or by relation {\rm(\ref{rel2})} and equation
{\rm(\ref{eq6})}, in which case $n$ is odd and $\gamma(A) =
n^2-n$. \item[{\rm (iii)}]  $({\cal E},{\cal P}(A,K))$ is given by
Figure $2$ and $K$ is decomposable if and only if after
normalization the relation on $\fn{Ext}K$ and the matrix $A$ are
given respectively by relation {\rm(\ref{rel1})} and equation
{\rm(\ref{eq1})}, in which case $n$ is even and $\gamma(A) =
n^2-n$.
\end{enumerate}
\end{lemma}

{\it Proof}.  The ``only if" parts of (i), (ii) and (iii) are done
in Lemma \ref{lemma 5.1}.  It remains to treat the ``if" parts and
the parts concerning the value of $\gamma(A)$.

(i)  First, suppose the relation on $\fn{Ext}K$ and the matrix $A$
are given respectively by relation {\rm(\ref{rel2})} and equation
{\rm(\ref{eq2})}. It is clear that $m$ is even, and so $n$ is odd.
Note that $A$ is well-defined, as it preserves the relation on
$\fn{Ext}K$. We contend that the digraph $({\cal E},{\cal P})$ is
given by Figure $1$.  It is clear that for $i = 1, \ldots, m-1,
(\Phi(x_i),\Phi(x_{i+1}))$ is the only outgoing arc from vertex
$\Phi(x_i)$.  By definition $Ax_m = x_1+\alpha x_2$, so
$\Phi(Ax_m)$ equals $\Phi(x_1+x_2)$, which is the smallest face of
$K$ containing $x_1,x_2$.  Since relation (\ref{rel2}), the
(unique) relation on the extreme vectors of $K$, is full, $K$ is
indecomposable.  As $x_1, x_2$ lie on opposite sides of
(\ref{rel2}), by Theorem A, $x_1, x_2$ lie on a common maximal
face of $K$, i.e., $\Phi(x_1+x_2)$ is a nontrivial face. But every
nontrivial face of an indecomposable minimal cone is simplicial,
so $x_1,x_2$ are the only extreme vectors of $\Phi(x_1+x_2)$.  It
follows that $(\Phi(x_m),\Phi(x_1))$ and $(\Phi(x_m),\Phi(x_2))$
are the only outgoing arcs from vertex $\Phi(x_m)$.  This proves
our contention.  Now a straightforward calculation yields the
following: $A^{m-1}x_1 = (1+\alpha)x_m$; $A^m x_1 =
(1+\alpha)(x_1+\alpha x_2)$, i.e., $\Phi(A^m x_1) =
\Phi(x_1+x_2)$; and $\Phi(A^{j(m-1)}x_1) =
\Phi(x_{m-j+1}+x_{m-j+2}+ \cdots +x_{m-1}+x_m)$ for $j = 1,\ldots,
m-2$.  So $A^{(n-1)(m-1)}x_1$ is a positive linear combination of
$x_3, x_4, \ldots, x_m$ and by Theorem A it belongs to the
relative interior of a maximal face of $K$.  On the other hand,
$A^{(n-1)(m-1)+1}x_1$ belongs to $\fn{int}K$ as it can be written
as a positive linear combination of all $x_i$ except $x_3$. Thus
$\gamma(A,x_1) = (n-1)(m-1)+1 = n^2-n+1$.  But $({\cal E},{\cal
P})$ is given by Figure $1$, so by Lemma \ref{lemma4.2}(ii),
$\gamma(A) = \gamma(A,x_1) = n^2-n+1$.

Next, suppose that the relation on $\fn{Ext}K$ and the matrix $A$
are given by relation {\rm(\ref{rel3})} and equation
{\rm(\ref{eq3})} respectively.  Then $K$ is indecomposable. Note
the left side of relation (\ref{rel3}) has one term more than its
right side and it contains both $x_1,x_2$.  Since $m (\ge 4)$ is
odd, the left side has at least three terms; so $Ax_m =
x_1+(1+\alpha)x_2\in
\partial K$.  As before one can verify that $({\cal E},{\cal P})$ is given by
Figure $1$.  Also, a straightforward calculation shows that
$A^{(n-1)(m-1)-1}x_1$, being a positive linear combination of
$x_2, x_3, \ldots, x_{m-1}$, belongs to $\partial K$, whereas
$A^{(n-1)(m-1)}x_1$, being a positive linear combination of $x_3,
x_4, \ldots, x_m$, belongs to $\fn{int}K$.  Since $({\cal E},{\cal
P})$ is given by Figure $1$, we have $\gamma(A) = \gamma(A,x_1) =
(n-1)(m-1) = n^2-n$.

 (ii)  Suppose the relation on $\fn{Ext}K$ and the matrix $A$ are
given by relation (\ref{rel5}) and equation (\ref{eq5})
respectively.  Since relation (\ref{rel5}) is full, $K$ is
indecomposable.  Note that the left side of (\ref{rel5}) has at
least three terms and it contains both $x_2$ and $x_3$, so $Ax_1$,
which is a positive linear combination of $x_2,x_3$, belongs to
$\partial K$.  Similarly, $Ax_m (= x_1+\alpha x_2)$ also belongs
to $\partial K$, as $x_1, x_2$ lie on opposite sides of
(\ref{rel5}).  By the same kind of argument as given in the proof
for part (i), one readily verifies that $({\cal E},{\cal P})$ is
given by Figure $2$.  Also, $A^{(n-1)(m-1)-1}x_2$, being a
positive linear combination of $x_3,x_4, \ldots, x_m$, belongs to
$\partial K$, whereas $A^{(n-1)(m-1)}x_2$ belongs to $\fn{int}K$,
as it can be written as a positive linear combination of all the
$x_i$'s except $x_3$.
 Since $({\cal E},{\cal
P})$ is given by Figure $2$, we have $\gamma(A) = \gamma(A,x_2) =
(n-1)(m-1) = n^2-n$.

When the relation on $\fn{Ext}K$ and the matrix $A$ are given by
relation (\ref{rel3}) and equation (\ref{eq4}) respectively, we
can show that $K$ is indecomposable and also that $({\cal E},{\cal
P})$ is given by Figure $2$.  In this case, $A^{(n-1)(m-1)-2}x_2$,
being a positive linear combination of $x_2,x_3, \ldots, x_{m-1}$,
belongs to $\partial K$, whereas $A^{(n-1)(m-1)-1}x_2$, being a
positive linear combination of $x_3,x_4, \ldots, x_m$, belongs to
$\fn{int}K$.  It follows that we have $\gamma(A) = \gamma(A,x_2) =
(n-1)(m-1)-1 = n^2-n-1$.

Similarly, we can show that when the relation and the matrix $A$
are given by relation (\ref{rel2}) and equation (\ref{eq6})
respectively, $K$ is indecomposable, $({\cal E},{\cal P})$ is
given by Figure $2$ and $\gamma(A) = n^2-n$ .

(iii)
 Suppose the relation on $\fn{Ext}K$ and the matrix
$A$ are given by relation (\ref{rel1}) and equation (\ref{eq1})
respectively.  In this case $K = \fn{pos} \{ x_2 \}\oplus \fn{pos}
\{ x_1,x_j, 3 \le j \le m \}$.  We readily check that
 $({\cal E},{\cal P})$ is given by Figure $2$.
A straightforward calculation shows that $A^{(m-1)(m-2)-1}x_2$ is
a positive linear combination of $x_3, x_4, \ldots, x_m$.  So
$A^{(n-1)(m-1)-1}x_2$ belongs to the indecomposable summand
$\fn{pos}\{ x_1,x_j, 3 \le j \le m \}$ of $K$ and hence lies in
$\partial K$. On the other hand, $A^{(n-1)(m-1)}x_2$ belongs to
$\fn{int}K$, as it can be written as a positive linear combination
of all $x_i$'s except $x_3$.  Hence we have $\gamma(A) =
\gamma(A,x_2) = (n-1)(m-1) = n(n-1)$. \hfill$\blacksquare$

\begin{theorem}\label{theorem2}  Let $n \ge 3$ be a given positive integer.
\begin{enumerate}
\item[{\rm (I)}] The quantity $\fn{max}\{\gamma(K): K\in {\cal
P}(n+1,n)\}$ equals $n^2-n+1$ if $n$ is odd and equals $n^2-n$ if
$n$ is even. \item[{\rm(II)}]  Suppose $n$ is odd.
\begin{enumerate}
\item[{\rm(i)}]  An $n$-dimensional minimal cone is exp-maximal if
and only if the cone is indecomposable and the relation for its
extreme vectors is balanced. \item[{\rm(ii)}] Let $K$ be an
$n$-dimensional exp-maximal minimal cone.  A $K$-primitive matrix
$A$ is exp-maximal if and only if the digraph $({\cal E},{\cal
P}(A,K))$ is given by Figure $1$.
\end{enumerate}
\item[{\rm(III)}]  Suppose $n$ is even.
\begin{enumerate}
\item[{\rm(i)}] An $n$-dimensional minimal cone is exp-maximal if
and only if either the cone is indecomposable and has a balanced
relation for its extreme vectors, or it is the direct sum of a ray
and an indecomposable minimal cone with a balanced relation for
its extreme vectors.  \item[{\rm(ii)}] Let $K$ be an
indecomposable exp-maximal minimal cone.   A $K$-primitive matrix
$A$ is exp-maximal if and only if the digraph $({\cal E},{\cal
P}(A,K))$ is, upon relabelling its vertices suitably, given by
Figure $1$ or Figure $2$, and in the latter case $x_1,x_2$ are
required to appear on opposite sides of the relation for the
extreme vectors of $K$. \item[{\rm(iii)}] Let $K$ be a
decomposable exp-maximal minimal cone.  A $K$-primitive matrix $A$
is exp-maximal if and only if the digraph $({\cal E},{\cal
P}(A,K))$ is given by Figure $2$.
\end{enumerate}
\end{enumerate}
\end{theorem}

{\it Proof}.  We first observe that when $n$ is even, there is no
minimal cone $K$ such that $\gamma(K) = n^2-n+1$. Assume to the
contrary that there is one such $K$.  Choose a $K$-primitive
matrix $A$ that satisfies $\gamma(A) = n^2-n+1$.  By Corollary
\ref{coro4} $({\cal E},{\cal P})$ is given by Figure $1$. Since
$n$ is even, by the second half of Lemma \ref{lemma 5.1}(i), after
normalization, the relation on $\fn{Ext}K$ and the matrix $A$ are
given by relation (\ref{rel3}) and equation (\ref{eq3})
respectively. So by Lemma \ref{lemma?}(i), we have $\gamma(A) =
n^2-n$, which is a contradiction.

 For any positive integer $n$, by Corollary \ref{coro4}, $\gamma(K) \le n^2-n+1$ for every
$n$-dimensional minimal cone $K$.

Let $n$ be odd.  Take any $n$-dimensional indecomposable cone $K$
for which the relation on its extreme vectors has the same number
of terms on its two sides. After re-indexing and normalizing the
extreme vectors $x_1, \ldots, x_m$ of $K$, we may assume that the
relation on $\fn{Ext}K$ is given by (\ref{rel2}). Let $A$ be the
$n \times n$ real matrix given by (\ref{eq2}). By Lemma
\ref{lemma?}(i) $A$ is $K$-primitive and $\gamma(A) = n^2-n+1$. So
we have
 $\gamma(K) = n^2-n+1$.  This establishes (I) for odd $n$ as well as
 the "if" part of (II)(i).

Now let $n$ be even.  In view of the above observations, the
maximum value of $\gamma(K)$ as $K$ runs through all
$n$-dimensional minimal cones is at most $n^2-n$. We are going to
show that the value $n^2-n$ can be attained.

Take any indecomposable minimal cone $K$ such that in the relation
on $\fn{Ext}K$ the number of vectors on its two sides differs by
$1$.  Normalizing the extreme vectors of $K$, we may assume that
the relation is given by (\ref{rel3}). Let $A$ be the matrix given
by (\ref{eq3}). By Lemma \ref{lemma?}(i) we have $\gamma(A) =
n^2-n$.  For this $K$, certainly we have $\gamma(K) = n^2-n$. This
establishes (I) for even $n$ and completes the proof for (I).

If $K$ is the direct sum of a ray and an indecomposable minimal
cone for which the relation on its extreme vectors has same number
of terms on its two sides, then necessarily $n$ is even and after
normalization we may assume that the relation is given by
(\ref{rel1}). By Lemma \ref{lemma?}(iii), the matrix $A$ defined
by (\ref{eq1}) satisfies $\gamma(A) = n^2-n$.

So we have also established the ``if" part of (III)(i).

To prove the ``only if" part of (II)(i), assume that $n$ is odd
and let $K$ be an $n$-dimensional minimal cone that satisfies
$\gamma(K) = n^2-n+1$. By Corollary \ref{coro4} there exists a
$K$-primitive matrix $A$ such that the digraph $({\cal E},{\cal
P})$ is given by Figure 1. By
 Lemma \ref{lemma4.2}(iii), $K$ is indecomposable.  By
 Lemma \ref{lemma 5.1}(i), after normalizing the
 extreme vectors of $K$, we may assume that the relation on $\fn{Ext}K$ is given by relation (\ref{rel2}).
So the relation has the same number of terms on its two
sides.

The ``only if" part of (II)(ii) follows from part(i) of Theorem
\ref{theorem1} (by taking $m = n+1$ and $m_A = n$), whereas its
``if" part is a consequence of Lemma \ref{lemma?}(i).

The proof of part(II) is complete.

To prove the ``only if" part of (III)(i), assume that $n$ is even
and let $K$ be an $n$-dimensional minimal cone such that
$\gamma(K) = n^2-n$. Choose a $K$-primitive matrix $A$ such that
$\gamma(A) = \gamma(K)$.  By part (ii) of Theorem \ref{theorem1},
in this case we have $m_A = n$ and either $({\cal E},{\cal P})$ is
given by Figure $1$ or Figure $2$, or $m_A = 3$.  The case $m_A =
3$ cannot happen; otherwise, $n = 3$, contradicting the assumption
that $n$ is even. Then, by Lemma \ref{lemma4.2}(iii), $K$ is
either indecomposable  or is the direct sum of a ray and an
indecomposable minimal cone for which the relation on its extreme
vectors has the same number of terms on its two sides.  In the
latter case, we are done.  In the former case, $({\cal E},{\cal
P})$ is given by Figure $1$ or Figure $2$.  If the digraph is
given by Figure $1$
 then, since $K$ is indecomposable minimal, by Lemma \ref{lemma 5.1}(i), after
   normalization, the relation on $\fn{Ext}K$ is given by
(\ref{rel3}).  If the digraph is given by Figure $2$, then by part
(ii) of the same lemma, after normalization, the relation on
$\fn{Ext}K$ and the matrix $A$ are given respectively by
(\ref{rel3}) and (\ref{eq4}) or (\ref{rel5}) and (\ref{eq5}). (We
have to rule out the former possibility, because by Lemma
 \ref{lemma?}(ii) we have $\gamma(A) = n^2-n-1$, which is a
contradiction.)  In any case, in the relation on $\fn{Ext} K$ the
number of terms on its two sides differ by $1$.

To prove the ``only if" part of (III)(ii), let $K$ be an
indecomposable minimal cone such that in the relation on its
extreme vectors the number of terms on its two sides differ by
$1$, and suppose that $A$ is a $K$-primitive matrix such that
$\gamma(A) = n^2-n$.  By the above proof for the ``only if" part
of (III)(i), $({\cal E},{\cal P})$ is given by Figure $1$ or
Figure $2$.  If it is given by Figure $2$ then after normalization
the relation on $\fn{Ext}K$ is given by relation (\ref{rel5}),
hence $x_1, x_2$ lie on opposite sides of the relation.

 To prove the ``if" part of (III)(ii), first suppose that $({\cal
E},{\cal P})$ is given by Figure $1$.  Since $n$ is even, by Lemma
\ref{lemma?}(i), after normalization, the relation on $\fn{Ext} K$
and the matrix $A$ are given respectively by (\ref{rel3}) and
(\ref{eq3}), and we have $\gamma(A) = n^2-n$.  If the digraph is
given by Figure $2$ and $x_1,x_2$ appear on opposite sides of the
relation on $\fn{Ext} K$, then since $K$ is indecomposable and $n$
is even, by Lemma \ref{lemma?}(ii) after normalization the
relation on $\fn{Ext} K$ and the matrix $A$ are given respectively
by (\ref{rel5}) and (\ref{eq5}), and we have $\gamma(A) = n^2-n$.

The ``if" part of (III)(iii) follows from  Lemma
\ref{lemma?}(iii). To prove its ``only if" part, let $A$ be a
$K$-primitive matrix such that $\gamma(A) = n^2-n$. As explained
in the above proof for the ``only if" part of (III)(i), the
digraph $({\cal E},{\cal P})$ is given by either Figure $1$ or
Figure $2$. By Lemma \ref{lemma4.2}(iii), if it is given by Figure
$1$ then $K$ is necessarily indecomposable. But now $K$ is
decomposable, so the digraph is given by Figure $2$.

The proof is complete. \hfill$\blacksquare$\\

By Corollary \ref{coro4} and Theorem \ref{theorem2}(I) we readily
deduce the following result, which is an improvement of the
already-proved Kirkland's conjecture.

\begin{coro}\label{coro5.4}  For any positive integer $m \ge 4$,
the maximum value of $\gamma(K)$ as $K$ runs through
 non-simplicial polyhedral cones with $m$
extreme rays and of all possible dimensions is $m^2-3m+3$ if $m$
is even, and is $m^2-3m+2$ if $m$ is odd.
\end{coro}

\setcounter{equation}{0}
\section{The 3-dimensional case}

Two distinct extreme rays $\Phi(x),\Phi(y)$ (or, distinct extreme vectors $x,y$) of $K$ are said to be {\it neighborly} if
$x+y \in \partial K$.

\begin{lemma}\label{lemma6.1} Let $K$ be a
$3$-dimensional polyhedral cone with extreme vectors $x_1, \ldots,
x_m$, where $m \ge 3$.  Let $A \in \pi(K)$ and suppose that
$({\cal E},{\cal P}(A,K))$ is given by Figure $1$ or Figure $2$.
Then\,{\rm:}
\begin{enumerate}
\item[{\rm(i)}]
 For $i = 1, \ldots, m,
\Phi(x_i)$ and $\Phi(x_{i+1})$ {\rm (}where $\Phi(x_{m+1})$ is
taken to be $\Phi(x_1)${\rm )} are neighborly extreme rays of $K$.
\item[{\rm(ii)}] $\gamma(A)$ equals $2m-1$ or $2m-2$ depending on
whether $({\cal E},{\cal P}(A,K))$ is given by Figure $1$ or by
Figure $2$.
\end{enumerate}
\end{lemma}

Note that for $i = 1, \ldots, m, \Phi(x_i)$ and $\Phi(x_{i+1})$
are adjacent vertices of the digraph $({\cal E},{\cal P})$ (when
it is given by Figure $1$ or Figure $2$). However, it is not clear
that $\Phi(x_i)$ and $\Phi(x_{i+1})$ are neighborly extreme
rays of the cone $K$.\\

 {\it Proof of Lemma 6.1}.  First, we consider the case when the digraph
$({\cal E},{\cal P})$ is given by
 Figure $1$.  It is not difficult to establish the following:\\

{\bf Assertion}.  Let $C$ be a convex polygon in ${\Bbb R}^2$ with
extreme points $w_1, \ldots, w_m$ and edges
$\overline{w_iw_{i+1}}, i = 1, \ldots, m$, where $w_{m+1}$ is
taken to be $w_1$ and $\overline{w_iw_{i+1}}$ denotes the line
segment joining the points $w_i$ and $w_{i+1}$.  Let $\tilde{C}$
be the polygon with extreme points $w_1', w_2, w_3, \ldots, w_m$,
where $w_1' = (1-\lambda)w_1+ \lambda w_2$ for some $\lambda, 0 <
\lambda < 1$. Then the edges of $\tilde{C}$ are
$\overline{w_1'w_2},
\overline{w_iw_{i+1}}, i = 2, \ldots, m-1$ and $\overline{w_mw_1'}$.\\
 The fact that $({\cal
E},{\cal P})$ is given by Figure $1$ implies that $Ax_i$ is a
positive multiple of $x_{i+1}$ for $i = 1, \ldots, m-1$ and,
moreover, $\Phi(Ax_m)$ is a $2$-dimensional face of $K$ with
extreme rays $\Phi(x_1)$ and $\Phi(x_2)$.  Hence, $\Phi(x_1),
\Phi(x_2)$ are neighborly extreme rays of $K$ and
 we have $Ax_m = \alpha_1x_1+\alpha_2x_2$ for some positive numbers $\alpha_1, \alpha_2$.
So $AK$ is generated by the (pairwise distinct) extreme vectors
$x_1', x_2,x_3, \ldots, x_m$, where $x_1' := Ax_m$.
 Using an equivalent formulation of the above assertion in terms of $3$-dimensional
polyhedral cones, we see that
 for all $i, j, 2 \le i,j \le m, x_i, x_j$ are
neighborly extreme vectors of $AK$ if and only if they are
neighborly extreme vectors of $K$.  (Note that here we are not
using (i), something that we have not yet established.)

By Lemma \ref{lemma4.2}(i) $A$ is nonsingular; so the cones $K$
and $AK$ are linearly isomorphic under $A$.  Since $x_1,x_2$ are
neighborly extreme vectors of $K, Ax_1,Ax_2$ are neighborly
extreme vectors of $AK$. But $Ax_1, Ax_2$ are respectively
positive multiples of $x_2$ and $x_3$, so $x_2,x_3$ are neighborly
extreme vectors of $AK$ and, in view of what we have done above,
$x_2,x_3$ are also neighborly extreme vectors of $K$. By repeating
the argument, we can show that for $i = 2, 3, \ldots, m-1$, $x_i$
and $x_{i+1}$ (also, $x_i$ and $x_{i-1}$) are neighborly extreme
vectors of $K$. It is clear that the remaining extreme vector
$x_1$ is neighborly to $x_m$ (and $x_2$).

 By direct calculation, $A^{2(m-1)}x_1$ is a positive linear combination of $x_{m-1}$ and
$x_m$, and as $x_{m-1}$ and $x_m$ are neighborly extreme vectors,
$A^{2(m-1)}x_1 \in \partial K$. On the other hand, $A^{2m-1}x_1$
is a positive linear combination of $x_m,x_1$ and $x_2$, so it
belongs to $\fn{int}K$.  This shows that $\gamma(A,x_1) = 2m-1$.
In view of Lemma \ref{lemma4.2}(ii), we have $\gamma(A) =
\gamma(A,x_1) = 2m-1$.

When the digraph $({\cal E},{\cal P})$ is given by Figure $2$, we
employ a similar argument.  For convenience, we denote the extreme
vectors of $AK$ by $y_1, \ldots, y_m$, where $y_1 = Ax_m, y_2 =
Ax_1$ and $y_i = x_i$ for $i = 3, \ldots, m$.  In this case, $y_1$
(respectively, $y_2$) is a positive linear combination of $x_1$
and $x_2$ (respectively, $x_2$ and $x_3$), and the extreme vectors
$x_2, x_1$ of $K$ are neighborly, and so are $x_2$ and $x_3$.
Moreover, $K$ and $AK$ are still linearly isomorphic under $A$.
Using an assertion similar to the one given above, we can show
that for $i, j = 3, \ldots, m, y_i,y_j$ are neighborly extreme
vectors of $AK$ if and only if $x_i,x_j$ are neighborly extreme
vectors of $K$.  Inductively we can show that $x_i,x_{i+1}$ (also,
$x_i$ and $x_{i-1}$) are neighborly extreme vectors of $K$ for $i
= 3, \ldots, m-1$.  Also, we can conclude that the remaining
extreme vector $x_1$ of $K$ is neighborly to $x_m$ (and $x_2$).

By direct calculation, $A^{2m-3}x_2$ is a positive linear
combination of $x_{m-1}$ and $x_m$, whereas $A^{2m-2}x_2$ is a
positive linear combination of $x_m,x_1$ and $x_2$; so
$A^{2m-3}x_2 \in \partial K$ and $A^{2m-2}x_2 \in \fn{int}K$.  By
Lemma \ref{lemma4.2}(ii) we have $\gamma(A) = \gamma(A,x_2) =
2m-2.$ \hfill$\blacksquare$\\

According to Lemma \ref{lemma4.2}(i), for any $K\in {\cal P}(m,n),
A\in \pi(K)$, if $({\cal E},{\cal P}(A,K))$ is given by Figure $1$
or Figure $2$, then $A$ has an annihilating polynomial of the form
$t^m-ct-d$, where $c,d > 0$.  Note that if, in addition, $\rho(A)
= 1$ then necessarily $c+d = 1$.  Conversely, if $A$ is
$K$-nonnegative and has an annihilating polynomial of the form
$t^m-ct-(1-c)$, where $0 < c < 1$, then necessarily $\rho(A) = 1$.
This follows from an application of the Perron-Frobenius theory to
$A$, because then $1$ is the only positive real root of the
polynomial, in view of Descartes' rule of signs, which says that
the number of positive roots of a polynomial either is equal to
the number of its variations of sign or is less than that number
by an even integer, a root of multiplicity $k$ being counted as
$k$ roots.  (The fact that $1$ is the only positive real root of
the polynomial can also be shown directly by proper
factorization.)

We will see that in constructing examples of $K$ and $A$ such that
$({\cal E},{\cal P}(A,K))$ is given by Figure $1$, polynomials of
the form $t^m-ct-(1-c)$, where $0 < c < 1$, play a role.  In our
next result, we study the roots of a polynomial of the said form.

\begin{lemma}\label{lemma6.2}
Consider the polynomial
$$h(t) = t^m-ct-(1-c),$$ where $m\ge 3$ and $c\in (0,1)$.
\begin{enumerate}
\item[{\rm(i)}] The roots of $h(t)$
 are all simple unless $m$ is odd and $c = c_m$,
where $c_m$ denotes the unique real root in $(0,1)$ of the
equation
$$\frac{(m-1)^{m-1}}{m^m}t^m = (t-1)^{m-1}.$$ \item[{\rm(ii)}]
When $m$ is even, $h(t)$ has precisely one real root other than
$1$.  When $m$ is odd, besides the root $1$, $h(t)$ has precisely
two real roots if $c > c_m$, one real root {\rm(}which is a double
root{\rm)} if $c = c_m$ and no real roots if $c < c_m$.  In all
cases, each real root of $h(t)$ other than $1$ lies in $(-1,0)$.
 \item[{\rm(iii)}] When $m \ge 4$ or when $m = 3$ and $c < \frac{3}{4} {\rm (= c_3)}$, the polynomial
$h(t)$ has a unique complex root of the form $re^{i\theta}$, where
$r > 0$ and $\theta\in (\frac{2\pi}{m},\frac{2\pi}{m-1})$.
Moreover, $r$ and $\theta$ are related in that $r$ is the unique
positive real root, which is less than $1$, of the polynomial
$g_{\theta}(t)$ given by $\!${\rm:}
$$g_\theta(t) = \frac{\sin(m-1)\theta}{\sin\theta}t^m-\frac{\sin m\theta}{\sin\theta}t^{m-1}+1.$$
\end{enumerate}
\end{lemma}

Three remarks are in order. First, by Descartes' rule of signs,
one can show that the polynomial $h(t)$, considered in the lemma,
has exactly one positive real root, and also that it has exactly
one negative real root if $m$ is even and either two (counting
multiplicities) or no negative real roots if $m$ is odd. Moreover,
since the coefficients of the polynomial $h(-t-1)$ are all
positive if $m$ is even, and all negative if $m$ is odd, the
polynomial $h(t)$ always has no real root less than $-1$.  Of
course, this agrees with part (ii) of the lemma,
 but the lemma contains more information.  Second, actually every root of $h(t)$ other than $1$ has
modulus strictly less than $1$. Here is a one-line proof: If
$\lambda$ is a root of $h(t)$, then $|\lambda|^m = |c\lambda +
(1-c)| \le c|\lambda|+(1-c)$, which is possible only if $\lambda =
1$ or $|\lambda| < 1$.  Third, part (iii) of the lemma follows
from properly combining Theorem 2 in \cite{Kir2} and Lemma 3 in
\cite{Kir1}, at least for $m \ge 4$. (The assumptions made in
\cite{Kir2}, namely (2.2) there, seems to rule out the case $m =
3$ of our theorem.  In his notation we have $d = n, k = 1$ and $s
= 1$ and his $n$ is our $m$.) We give a proof for the sake of
completeness.\\

{\it Proof of Lemma \ref{lemma6.2}}. (i)  Clearly, the roots of
$h'(t)$ are precisely all the $(m-1)$th roots of $\frac{c}{m}$.  A
little calculation shows that if $t_0$ is a common root of $h(t)$
and $h'(t)$ then $\frac{ct_0}{m} = t_0^m = ct_0+(1-c)$ and so
$ct_0(\frac{1}{m}-1) = 1-c$, which implies that $t_0$ is the
negative $(m-1)$th real root of $\frac{c}{m}$ and $m$ is odd.  So
$h(t)$ and $h'(t)$ have a common root if and only if $m$ is odd
and the negative $(m-1)$th real root of $\frac{c}{m}$ is a root of
$h(t)$.  By calculation one finds that the latter condition, in
turn, is equivalent to the condition that $m$ is odd and $c$ is a
root of the polynomial $\alpha_mt^m - (t-1)^{m-1}$, where
$\alpha_m = \frac{(m-1)^{m-1}}{m^m}$.  When $m$ is odd, by
considering the said polynomial and its derivative we readily show
that the polynomial has a unique real root in $(0,1)$, which we
denote by $c_m$. So we can conclude that
 the roots of $h(t)$ are all simple unless $m$
is odd and $c=c_m$.

(ii)  Rewriting $h(t)$, we have, $h(t) = (t^m-1)-c(t-1) =
(t-1)(p(t)-c)$, where $p(t) = t^{m-1}+t^{m-2}+ \cdots + t+1.$ So
$1$ is always a root of $h(t)$ and for any complex number $w \ne
1$, $w$ is a root of $h(t)$ if and only if $w$ is a root of the
equation $p(t) = c$.

When $m$ is even, a consideration of the derivative of $p(t)$
shows that $p(t)$ is a strictly increasing function on the real
line. But $p(-1) = 0, p(0) = 1$ and $c\in (0,1)$, so the equation
$p(t) = c$ has exactly one real root and that real root belongs to
$(-1,0)$. Hence $h(t)$ has precisely one real root other than $1$
and this root lies in $(-1,0)$.

When $m$ is odd, on the real line $p(t)$ is a strictly convex
function, since its second derivative always takes positive
values, as can be shown by some calculation.  It is
straightforward to show that a complex number $z_0$ is a common
root of $h(t)$ and $h'(t)$ if and only if $p(z_0) = c$ and
$p'(z_0) = 0$.
 On the other hand,
by part (i) and its proof, $h(t)$ and $h'(t)$ have a common root
if and only if $c = c_m$; in that case, the common root is unique
and is equal to the negative $(m-1)$th real root of
$\frac{c_m}{m}$.  So $c_m$ is, in fact, the absolute minimum value
of $p(t)$.  Hence, the equation $p(t) = c$ has two distinct real
roots if $c
> c_m$, one real root (which is a double root) if $c = c_m$, and
no real roots if $c < c_m$. As $p(-1) = p(0) = 1$ and $p(t)$ is
strictly convex, each real root of the equation $p(t) = c$ or, in
other words, each real root of $h(t)$ other than $1$, must belong
to $(-1,0)$.

(iii)  First, we establish the uniqueness of the root of $h(t)$ in
the desired polar form.  Let $z_1 = r_1e^{i\theta_1}, z_2 =
r_2e^{i\theta_2}$, where $r_1, r_2
> 0$ and $\theta_1, \theta_2\in
(\frac{2\pi}{m},\frac{2\pi}{m-1})$, be two different roots of
$h(t)$. Then from $z_1^m-cz_1-(1-c) = 0$ and $z_2^m-cz_2-(1-c) =
0$ we obtain $(z_1-z_2)(\sum_{k=0}^{m-1}z_1^{m-1-k}z_2^k) =
c(z_1-z_2)$, which implies that $c =
\sum_{k=0}^{m-1}z_1^{m-1-k}z_2^k$.  For any nonzero complex number
$z$, denote by $\fn{arg}(z)$ the argument of $z$ that belongs to
$(0, 2\pi]$.  Since $\theta_1, \theta_2\in
(\frac{2\pi}{m},\frac{2\pi}{m-1})$, we have
$\fn{arg}(z_1^{m-1-k}z_2^k)\in (\frac{2\pi(m-1)}{m},2\pi)$ for $k
= 0, \ldots, m-1$.  Thus, for each $k$, the complex number
$z_1^{m-1-k}z_2^k$ belongs to the relative interior of the convex
cone in the complex plane generated by $1$ and
$e^{-\frac{2\pi}{m}i}$, and hence so does the sum
$\sum_{k=0}^{m-1}z_1^{m-1-k}z_2^k$, which is a contradiction as
$c$ is a positive real number.

Next, we contend that for any real number $\theta\in
(\frac{2\pi}{m},\frac{2\pi}{m-1})$, the polynomial $g_\theta(t)$
has a unique positive real root and this root is less than $1$.

We have $g_\theta(0) = 1 > 0$ and
$$g_\theta(1) = \frac{\sin (m-1)\theta-\sin m\theta}{\sin \theta}+1 = -\frac{\cos (m-\frac{1}{2})\theta}{\cos \frac{\theta}{2}}+1 < 0,$$
where the second equality follows from the trigonometric identity
$\sin\alpha-\sin\beta =
2\cos\frac{\alpha+\beta}{2}\sin\frac{\alpha-\beta}{2}$ and the
inequality holds as $2\pi-\frac{1}{2}\theta <
(m-\frac{1}{2})\theta < 2\pi+\frac{\theta}{2}$ and $0 <
\frac{\theta}{2} < \frac{\pi}{2}$.  In addition, we also have
$$g'_\theta(t) = \frac{t^{m-2}}{\sin\theta}\left[ mt\sin(m-1)\theta-(m-1)\sin m\theta \right] < 0$$
for all $t \in (0,\infty)$, as $\sin m\theta > 0$ and
$\sin(m-1)\theta < 0$.  So it is clear that the polynomial
$g_\theta(t)$ has a unique positive real root and this root is
less than $1$.  (The latter assertion can also be established by
using Descartes' rule of signs.)

\indent Now define a real-valued function $\zeta$ on
$(\frac{2\pi}{m},\frac{2\pi}{m-1})$ by $\zeta(\theta) =
r_\theta^{m-1}\frac{\sin m\theta}{\sin \theta},$ where $r_\theta$
denotes the unique positive real root of the polynomial
$g_\theta(t)$.  Note that $\zeta$ is a continuous function, as
$r_\theta$ depends on $\theta$ continuously.  Also, we have $0 <
\zeta(\theta) < 1$, as $2\pi < m\theta < 2\pi+\theta$ and $0 <
r_\theta < 1$.

It is readily checked that a complex number $re^{i\theta}$ (in
polar form) is a root of the polynomial $h(t)$ if and only if we
have
\begin{eqnarray}
       (1-c)+cr\cos \theta &=& r^m\cos m\theta, \label{okay1}\\
                    cr\sin \theta &=& r^m\sin m\theta.
                    \label{okay2}
\end{eqnarray}
Rewriting (\ref{okay2}) and adding $\cos \theta$ times
(\ref{okay2}) to $-\sin \theta$ times (\ref{okay1}) (and noting
that $\sin \theta, \cos \theta \ne 0 \mbox{ for }\theta\in
(\frac{2\pi}{m},\frac{2\pi}{m-1})$), we find that (\ref{okay1})
and (\ref{okay2}) holds if and only if
$$ c = r^{m-1}\frac{\sin m\theta}{\sin \theta} = \frac{r^m\sin (m-1)\theta}{\sin \theta}+1.$$
\noindent So, for any $\theta\in
(\frac{2\pi}{m},\frac{2\pi}{m-1})$, $r_{\theta}e^{i\theta}$ is a
root of $h(t)$ if $c = \zeta(\theta)$.  To complete our proof, it
remains to show that the function $\zeta$ maps the open interval
$(\frac{2\pi}{m},\frac{2\pi}{m-1})$ onto $(0,1)$ when $m \ge 4$
and onto $(0,\frac{3}{4})$ when $m = 3$.

 Note that the function
$\zeta$ is one-to-one.  Otherwise, there exist
$\theta_1,\theta_2\in (\frac{2\pi}{m},\frac{2\pi}{m-1})$,
$\theta_1\ne \theta_2$, such that $\zeta(\theta_1) =
\zeta(\theta_2)$.  But then $r_{\theta_1}e^{i\theta_1}$ and
$r_{\theta_2}e^{i\theta_2}$ are distinct roots of the polynomial
$t^m-\zeta(\theta_1)t-(1-\zeta(\theta_1))$, both with argument
belonging to $(\frac{2\pi}{m},\frac{2\pi}{m-1})$, which
contradicts what we have obtained at the beginning of the proof of
part(iii). As $\zeta$ is one-to-one and continuous, it is either
strictly increasing or strictly decreasing on
$(\frac{2\pi}{m},\frac{2\pi}{m-1})$.  In view of the intermediate
value property of a real-valued continuous function, it suffices
to show that $\lim_{\theta \rightarrow \frac{2\pi}{m}^+}
\zeta(\theta) = 0$ and $\lim_{\theta \rightarrow
\frac{2\pi}{m-1}^-} \zeta(\theta)$ equals $1$ when $m\ge 4$ and
equals $\frac{3}{4}$ when $m = 3$.

 By the
definition of $\zeta$ and the fact that $0 < r_{\theta} < 1$ it is
readily seen that we have $\lim_{\theta \rightarrow
\frac{2\pi}{m}^+} \zeta(\theta) = 0$. On the other hand, from the
condition $g_\theta(r_\theta) = 0$ we also have $\zeta(\theta) =
\frac{r_\theta^m\sin \ (m-1)\theta}{\sin \theta}+1$, which implies
that $\lim_{\theta \rightarrow \frac{2\pi}{m-1}^-}\zeta(\theta) =
1$, provided that $m \ge 4$. When $m = 3$, the polynomial equation
given in (i) becomes $\frac{2^2}{3^3}t^3 = (t-1)^2.$ Since
$\frac{3}{4}$ is a root of the latter equation, we have $c_3 =
\frac{3}{4}$. Note that the polynomial $g_\theta(t)$ tends to
$-2t^3-3t^2+1$ as $\theta$ tends to $\pi$ from the left.  Also,
$\frac{1}{2}$ is a root of the latter polynomial.  So
$\lim_{\theta\rightarrow \pi^-}r_{\theta} = \frac{1}{2}$, and we
have
$$\lim_{\theta\rightarrow \pi^-}\zeta(\theta) =
\lim_{\theta\rightarrow \pi}r_{\theta}^2\frac{\sin 3\theta}{\sin
\theta} = \frac{3}{4},$$ as desired.  \hfill$\blacksquare$\\

In the course of the proof of Lemma \ref{lemma6.2} we have also
established the following (see also \cite[Lemma 3]{Kir1}):

\begin{coro}\label{coro6.3}\begin{enumerate}
\item[{\rm(i)}] For any $\theta\in
(\frac{2\pi}{m},\frac{2\pi}{m-1}), r_{\theta}e^{i\theta}$, where
$r_{\theta}$ is the unique positive real root of the polynomial
$g_{\theta}(t)$ given in Lemma \ref{lemma6.2}, is a root of the
polynomial $t^m-ct-(1-c)$ where $c = r_{\theta}^{m-1} \frac{\sin
m\theta}{\sin \theta}$. Moreover, the pair $(r_{\theta},\theta)$
satisfies the relations {\rm(\ref{okay1})} and {\rm(\ref{okay2})},
which appear in the proof of Lemma \ref{lemma6.2}.
\item[{\rm(ii)}]  For every positive integer $m \ge 4$
{\rm(}respectively, $m = 3${\rm)}, every real number $c$ in
$(0,1)$ {\rm(}respectively, in $(0,\frac{3}{4})${\rm)} can be
expressed uniquely as $r_{\theta}^{m-1} \frac{\sin m\theta}{\sin
\theta}$, where $\theta\in (\frac{2\pi}{m},\frac{2\pi}{m-1})$.
\end{enumerate}
\end{coro}

\begin{theorem}\label{coro6}
\begin{enumerate}
\item[{\rm(i)}] For every positive integer $m \ge 3$, $$\fn{max}\{
\gamma(K): K\in {\cal P}(m,3)\} = 2m-1.$$  \item[{\rm(ii)}] For
every $K\in {\cal P}(m,3)$, $K$ is exp-maximal if and only if
there exists a $K$-primitive matrix $A$ such that the digraph
$({\cal E},{\cal P}(A,K))$ is given by Figure $1$.
\item[{\rm(iii)}]  Let $m \ge 3$ be a positive integer. For any
$\theta\in (\frac{2\pi}{m},\frac{2\pi}{m-1})$, let $r_\theta$
denote the unique positive real root of the polynomial
$g_\theta(t)$ as defined in Lemma \ref{lemma6.2}{\rm(iii)}. Let
$K_{\theta}$ be the polyhedral cone in ${\Bbb R}^3$ generated by
the vectors
$$ x_j(\theta) := \left[\begin{array}{c}
r_\theta^{j-1}\cos(j-1)\theta\\
r_\theta^{j-1}\sin(j-1)\theta\\1 \end{array}\right],\ j = 1,
\ldots, m.$$ Also, let $A_\theta = r_\theta\left[\begin{array}{rr}\cos \theta&-\sin \theta\\
\sin \theta & \cos \theta \end{array}\right]\oplus [1]$.  Then
$x_1(\theta), \ldots, x_m(\theta)$ are the extreme vectors of
$K_{\theta}$,
  $K_\theta$ is an exp-maximal polyhedral cone and $A_\theta$ is
  an exp-maximal
$K_\theta$-primitive matrix.
\end{enumerate}
\end{theorem}

\noindent Note that, for simplicity, we suppress the dependence of
$K_{\theta}$ on $m$. \\

{\it Proof}.  For every $K\in {\cal P}(m,3)$, by Corollary
\ref{coro4} (with $n = 3$) we have $\gamma(K) \le 2m-1$, where the
equality holds only if there exists a $K$-primitive matrix $A$
such that the digraph $({\cal E},{\cal P})$ is given by Figure
$1$.  In view of Lemma \ref{lemma6.1}(ii), parts (i) and (ii) will
follow if we can
 construct a polyhedral cone $K\in {\cal P}(m,3)$
for which there exists a $K$-primitive matrix $A$ such that
$({\cal E},{\cal P}(A,K))$ is given by Figure $1$. To complete the
proof, we are going to establish part (iii) and at the same time
show that the digraph $({\cal E}, {\cal
P}(A_{\theta},K_{\theta}))$ is given by Figure $1$.

Consider any fixed $\theta\in (\frac{2\pi}{m},\frac{2\pi}{m-1})$.
Hereafter we denote $r_{\theta}$ simply by $r$. By Corollary
\ref{coro6.3}(i) $re^{i\theta}$ is a root of the polynomial
$t^m-ct-(1-c)$ where $c = r^{m-1} \frac{\sin m\theta}{\sin
\theta}$.  Consider the $m$ points $y_1, \ldots, y_m$ in ${\Bbb
R}^2$ given by:
$$y_1 = \left( \begin{array}{c}1\\0\end{array}\right) \mbox{ and } y_j =
\left( \begin{array}{c}r^{j-1}\cos(j-1)\theta\\
r^{j-1}\sin(j-1)\theta\end{array}\right) \mbox{ for } j = 2,
\ldots, m.$$  Take $B$ to be the $2 \times 2$ matrix $r\left[
\begin{array}{rr}\cos \theta&-\sin \theta \\ \sin\theta&\cos\theta
\end{array}\right].$ Note that $y_j = B^{j-1}y_1$ for $j = 2,
\ldots, m$.  By Corollary \ref{coro6.3}(i) the equations
(\ref{okay1}) and (\ref{okay2}), which appear in the proof of
Lemma \ref{lemma6.2}, are satisfied.  So we have
$$By_m = {r^m \cos m\theta \choose r^m \sin m\theta} = c{r\cos
\theta \choose r \sin \theta}+(1-c){1 \choose 0} =
cy_2+(1-c)y_1.$$ We contend that $y_1, \ldots, y_m$
 are precisely all the extreme points of the convex polygon
$C:=\fn{conv}\{ y_1, \ldots, y_m \}$, and moreover
$\overline{y_my_1}$ and $\overline{y_{i-1}y_i}$, for $i = 2,
\ldots, m$, are precisely all its sides.  Since $x_j(\theta) =
\left(
\begin{array}{c}y_j\\1 \end{array}\right)$ for $j = 1, \ldots, m$,
 $A_{\theta} = B \oplus [1]$ and $BC\subseteq C$, once the
contention is established, it is clear that $x_1(\theta), \ldots,
x_m(\theta)$ are all the extreme vectors of $K_{\theta}$, and
$A_{\theta}$ is $K_{\theta}$-nonnegative.  Then we are done, as it
can be
 readily checked that
the digraph $({\cal E},{\cal P}(A_{\theta},K_{\theta}))$ is given
by Figure $1$.

Clearly the extreme points of $C$ are among $y_1, \ldots, y_m$.
Since the Euclidean norm of $y_1$ is $1$, whereas that of $y_j$,
for $j = 2, \ldots, m$, is less than $1$, $y_1$ is certainly an
extreme point of $C$.

For $j = 2, \ldots, m$, proceeding inductively, we are going to
show that $y_j$ is an extreme point and $\overline{y_{j-1}y_j}$
forms a side of $C$.  Once this is proved, it will follow that
$\overline{y_my_1}$ is also a side of $C$ and we are done.

We begin with $j = 2$.  If $y_2$ is not an extreme point of $C$,
then $y_2$ lies in the relative interior of a line segment that
joins two distinct points of $C$.  But then
 $y_3$, which is $By_2$, also lies in
the relative interior of a line segment joining two distinct
points of $C$ (as $B$ is nonsingular and maps $C$ into itself) and
hence is not an extreme point of $C$. Continuing the argument, we
conclude that $y_2, \ldots, y_m$ are all not extreme points of
$C$; so $y_1$ is the only extreme point of $C$, which is
impossible.  This proves that $y_2$ must be an extreme point of
$C$.

To proceed further, we need the following:\\

{\bf Assertion}.   Let $M$ be a compact convex set in ${\Bbb R}^2$
which contains the origin as an interior point.  Let $u_1,u_2$ be
two distinct extreme points of $M$ which are not negative
multiples of each other.  If the line segment $\overline{u_1u_2}$
is not a face of $M$, then $M$ has an extreme point of the form
$\alpha_1u_1+\alpha_2u_2$, where
$\alpha_1, \alpha_2 > 0$.\\

{\it Proof}.  Let $L$ denote the straight line joining the points
$u_1, u_2$.  Since $u_1,u_2$ are not negative multiples of each
other, the origin is on one side of the line.  Move $L$ parallel
to itself but away from the origin until $L$ becomes a supporting
line $L^\prime$ for $M$.  (Of course, the argument can be given
more precisely in terms of a (continuous) linear functional.) Then
$F: = L^\prime \cap M$ is a (boundary) face of $M$ and hence is
either an extreme point or a line segment.
 Note that the point
$u_\lambda := \lambda(\frac{u_1+u_2}{2})$ belongs to $\fn{int}M$,
provided that $\lambda > 1$, sufficiently close to $1$.  Suppose
that $F$ has a point, say $v$, that does not belong to the
interior of the cone $\fn{pos}\{ u_1,u_2 \}$.  Since $u_1, u_2$
are extreme points of $M$ and $0\in M$, it is obvious that $v$
cannot lie on the boundary of the said cone.  So it lies outside
the cone. But then the line segment $\overline{u_{\lambda}v}$ will
meet one of the boundary rays of the said cone at a point of the
form $\alpha u_j$, where $\alpha > 1$ and $j = 1$ or $2$, which
contradicts the assumption that $u_1,u_2$ are extreme points of
$M$.  This shows that $F\subseteq \fn{int}(\fn{pos}\{ u_1, u_2
\})$. As each extreme point of $F$ is also an extreme point of
$M$, now it is clear that $M$ has an extreme point of the desired
form. \hfill$\blacksquare$

Now back to the proof of our theorem.  In view of the definition
of the $y_i$'s, it is clear that none of the points $y_3, \ldots,
y_m$ can be written as a positive linear combination of  $y_1$ and
$y_2$.  (This can be seen, for instance, by considering the
arguments of the complex numbers corresponding to these points.)
So, by the above Assertion, the line segment $\overline{y_1y_2}$
forms a side of the polygon $C$.

Note that for each $j, j = 3, \ldots, m-1$, by the preceding
argument we have used for $y_2$, we can show that if $y_j$ is not
an extreme point of $C$, then none of the points $y_{j+1}, \ldots,
y_m$ is an extreme point of $C$, and so we must have $C =
\fn{conv}\{ y_1,\ldots, y_{j-1} \}$.  But $C$ is not the convex
hull of all the $y_j$'s that lie on the upper closed half-plane,
so it follows each of points $y_1, \ldots, y_k$ must be an extreme
point of $C$, where we denote by $k$ the first positive integer
such that $k\theta
> \pi$ (i.e., $k = \lfloor\frac{\pi}{\theta}\rfloor+1$).
Moreover, for $j = 3, \ldots, k$, once we have obtained that
$y_{j-1}$ and $y_j$ are extreme points, like the case $j=2$, using
the above Assertion, we can infer that $\overline{y_{j-1}y_j}$
forms a side of the polygon $C$.

If our inductive argument should fail at one step, then there must
exist some $s, k\le s \le m-1$, such that $y_1, \ldots, y_s$ are
precisely all the extreme points of $C$ and moreover each of the
line segments $\overline{y_{j-1}y_j}$, for $j = 2, \ldots, s$,
forms a side of $C$.  Then the line segment $\overline{y_sy_1}$ is
clearly the remaining side of $C$.   Since $y_s$ is in the lower
open half-plane, a consideration of the arguments of the complex
numbers corresponding to the points $y_1, \ldots, y_m$ shows that
$y_m$ must be a positive linear combination of $y_1$ and $y_s$. So
we have either $y_m\in \fn{int}\fn{conv}\{y_1,y_s,0\}$ or $y_m\in
\fn{ri}\overline{y_1y_s}$.  In any case, $y_s$ belongs to
$\Phi(y_m)$, the face of $C$ generated by $y_m$.  Hence $y_{s+1} =
By_s\in \Phi(By_m) = \overline{y_1y_2}$, where the last equality
holds as $By_m = (1-c)y_1+cy_2 (0 < c < 1)$ and it has already
been shown that $\overline{y_1y_2}$ is a side (i.e., a face) of
$C$.  As $k+1 \le s+1 \le m$, we arrive at a contradiction.

So proceeding inductively, we show that $y_1, \ldots, y_m$ are all
the extreme points of $C$ and $\overline{y_{j-1}y_j}$ forms a side
of $C$ for $j = 2, \ldots, m$.  The proof is complete.
\hfill$\blacksquare$\\

In contrast to Theorem \ref{coro6}, we have the following:

\begin{theorem}\label{theorem6.4}
For every positive integer $m \ge 5$ there exists a
$3$-dimensional polyhedral cone $K$ with $m$ extreme rays for
which there is no $K$-primitive matrix $A$ with the digraph
$({\cal E},{\cal P}(A,K))$ given by Figure $1$.
\end{theorem}

{\it Proof}.   Let $K$ be the polyhedral cone in ${\Bbb R}^3$ with
extreme vectors $$y_j = (\cos \frac{2j\pi}{m}, \sin
\frac{2j\pi}{m},1)^T, j = 1, \ldots, m.$$  We contend that there
is no $K$-primitive matrix $A$ for which $({\cal E},{\cal P})$ is
given by Figure $1$, where $x_1, x_2, \ldots, x_m$ is a
rearrangement of $y_1, \ldots, y_m$.

We assume to the contrary that there is one such $A$.   By Lemma
\ref{lemma6.1}, for $i = 1, \ldots, m, x_i$ and $x_{i+1}$ are
neighborly extreme vectors of $K$ (where $x_{m+1}$ is taken to be
$x_1$).  On the other hand, for each $j$, the extreme vectors
neighborly to $y_j$ are $y_{j+1}$ and $y_{j-1}$.  [We adopt the
convention that for each integer $j, y_j$ equals $y_k$ where $k$
is the unique integer that satisfies $1 \le k \le m, k \equiv j~
(\fn{mod} m)$.] Suppose $x_1 = y_{j_1}$, where $1 \le j_1 \le m$.
Since $x_2$ is neighborly to $x_1$, $x_2$ must be either
$y_{j_1+1}$ or $y_{j_1-1}$. Consider first the case when $x_2 =
y_{j_1+1}$. Since $x_3$ is neighborly to $x_2$, it is equal to
either $y_{j_1+2}$ or $y_{j_1}$.  But we already have $x_1 =
y_{j_1}$, so $x_3$ must be $y_{j_1+2}$.  Continuing the argument,
we can show that $x_j = y_{j_1+j-1}$ for $j = 1, \ldots, m$.  Then
we take $\hat A$ to be $\left[
\begin{array}{cc}
\cos \frac{2\pi}{m} & -\sin \frac{2\pi}{m}\\
\sin \frac{2\pi}{m} & \cos \frac{2\pi}{m}
\end{array} \right] \oplus [1]$.
If $x_2 = y_{j_1-1}$, we can show in a similar manner that $x_j
=y_{j_1-j+1}$ for $j = 1, \ldots, m$.  In this case, we take
$\hat{A}$ to be $\left[
\begin{array}{cc}
\cos \frac{2\pi}{m} & \sin \frac{2\pi}{m}\\
-\sin \frac{2\pi}{m} & \cos \frac{2\pi}{m}
\end{array} \right] \oplus [1].$  As we are going to show,
in either case,  $\hat{A}$ and $A$ are positive multiples of each
other.

It is known that if $C_1,C_2$ are proper cones in ${\Bbb R}^{n_1}$
and ${\Bbb R}^{n_2}$ respectively then the set $\pi(C_1,C_2)$
which
 consists of all $n_2 \times n_1$ matrices $B$ such that $BC_1 \subseteq C_2$ is
a proper cone in the space of $n_2 \times n_1$ real matrices.

 Now let $K_1$ denote the $3$-dimensional polyhedral cone
$\fn{pos}\{ x_1,x_2,x_3,x_4 \}$. It is clear that $A \in
\pi(K_1,K)$.  As the digraph $({\cal E},{\cal P})$ is given by
Figure $1$ and $m \ge 5, Ax_i$ is a positive multiple of $x_{i+1}$
for $i = 1,2,3,4$.  So $A$ maps every extreme vector of $K_1$ to
an extreme vector of $K$.  But $K_1$ is indecomposable and $A$ is
nonsingular, so by a variant of a sufficient condition for a
cone-preserving map to be extreme due to Loewy and Schneider \cite
[Theorem 3.3]{L--S} (see \cite [Theorem 1]{Loe} or \cite [Theorem
5.2]{Tam 3}), it follows that $A$ is an extreme matrix of the
proper cone $\pi(K_1,K)$.  Our definition of $\hat{A}$ guarantees
that $\hat{A}x_i = x_{i+1}$ for $i = 1, \ldots, 4$; hence
$\hat{A}x_i$ is a positive multiple of $Ax_i$ for $i = 1, \ldots,
4$.  But $\{ x_1, \ldots, x_4 \}$ is the set of extreme vectors of
$K_1$, so $\hat{A}$ belongs to the extreme ray of $\pi(K_1,K)$
generated by $A$. Hence $A$ and $\hat{A}$ are positive multiples
of each other. But then the digraph $({\cal E},{\cal
P}(\hat{A},K))$ is not given by Figure $1$ (and also $\hat{A}$ is
not $K$-primitive). So we arrive at a contradiction.
\hfill$\blacksquare$

Let $K_1, K_2$ be linearly isomorphic proper cones.  If $D$ is a
digraph that can be realized as $({\cal E}, {\cal P}(A_1,K_1))$
for some $K_1$-nonnegative matrix $A_1$, then clearly (up to graph
isomorpism) $D$ can also be realized as $({\cal E}, {\cal
P}(A_2,K_2))$ for some $K_2$-nonnegative matrix $A_2$. On the
other hand, if $K_1, K_2$ are assumed to be combinatorially
equivalent only, then the same cannot be said.

\begin{remark} \rm
 Let $K_1, K_2$ be combinatorially equivalent proper cones. Then:
\begin{enumerate}
\item[(i)] If $G$ is a digraph such that $G = ({\cal E}(K_1),{\cal
P}(A_1,K_1))$ for some $K_1$-primitive matrix $A_1$, then there
need not exist a $K_2$-primitive matrix $A_2$ such that $({\cal
E}(K_2),{\cal P}(A_2,K_2))$ is isomorphic with $G$. \item[(ii)]
The values of $\gamma(K_1), \gamma(K_2)$ need not be the same.
\end{enumerate}
\end{remark}

Since any two $3$-dimensional polyhedral cones with the same
number of extreme rays are combinatorially equivalent, the
preceding remark follows from Theorem \ref{coro6} and Theorem
\ref{theorem6.4}.

\setcounter{equation}{0}
\section{The higher-dimensional case}

In this section we are going to establish the following main
result of our paper.

\begin{theorem}\label{maintheorem}  For any pair of positive integers $m,n, 3 \le n \le m$,
$\fn{max}\{\gamma(K): K\in {\cal P}(m,n)\}$ equals $(n-1)(m-1)+1$
when $m$ is even or $m$ and $n$ are both odd, and is at least
$(n-1)(m-1)$ and at most $(n-1)(m-1)+1$ when $m$ is odd and $n$ is
even.
\end{theorem}

First, we explain how the proof starts.
 By Corollary \ref{coro4} the inequality $\gamma(K) \le (n-1)(m-1)+1$ holds for any $K\in {\cal P}(m,n)$.
 By the same corollary, when the inequality becomes equality, there exists a
 $K$-primitive matrix $A$ such that the digraph $({\cal E},{\cal
 P}(A,K))$ is given by Figure $1$.  Then a straightforward
 computation shows that $A^{(n-1)(m-1)-1}x_1$ is a positive linear combination of
 the $n-1$ extreme vectors $x_{m-n+1}, x_{m-n+2}, \ldots, x_{m-1}$, whereas $A^{(n-1)(m-1)}x_1$ is a positive linear
 combination of the $n-1$ extreme vectors $x_{m-n+2},x_{m-n+3}, \ldots,
 x_{m-1},x_m$.  On the other hand, by Lemma \ref{lemma4.2}(ii) we also have $\gamma(A,x_1) = \gamma(A) \le
 (n-1)(m-1)+1$.  Note that the latter inequality implies that $A^{(n-1)(m-1)+1}x_1 \in
 \fn{int}K$.
 So, to establish Theorem \ref{maintheorem}, it suffices to
 show that when $m$ is even or when $m$ and $n$ are both odd (respectively, when $m$ is odd and
 $n$ is even), it
 is possible to construct a polyhedral cone $K\in {\cal P}(m,n)$ and a $K$-primitive matrix $A$ such that the digraph $({\cal E},{\cal
 P}(A,K))$ is given by Figure $1$ and $x_{m-n+2}+x_{m-n+3}+ \cdots +
 x_m\in \partial K$ (respectively, $x_{m-n+1}+x_{m-n+2}+ \cdots + x_{m-1}\in \partial K$).

 To begin with, we show that for every pair $m,n, 3\le n \le m$, the
 digraph given by Figure $1$ can always be realized as $({\cal
 E},{\cal P}(A,K))$ for some $K\in {\cal P}(m,n)$ and a $K$-primitive
 matrix $A$.

\begin{lemma}\label{theorem7.1}
For every pair of positive integers $m, n, 3 \le n \le m$, there
is a polyhedral cone $K\in {\cal P}(m,n)$ for which there exists a
$K$-primitive matrix $A$ such that the digraph $({\cal E},{\cal
P}(A,K))$ is given by Figure $1$.
\end{lemma}

{\it Proof}. First, we treat the case when $m, n$ are both odd.
Write $m = 2k+1$ and $n = 2p+1$.  Then $1\le p \le k$.  Choose any
$c\in (0,c_m)$, where $c_m \in (0,1)$ is the positive real number
defined in Lemma \ref{lemma6.2}, and let $h(t)$ denote the
polynomial $t^m-ct-(1-c)$.  By part (iii) of the said lemma $h(t)$
has a complex root of the form $r_1e^{i\theta_1}$, where
$\theta_1\in (\frac{2\pi}{m},\frac{2\pi}{m-1})$ and $0 < r_1 < 1$.
In addition, $h(t)$ has no real roots other than $1$.  So $h(t)$
has $k$ pairs of non-real conjugate complex roots, say, $r_je^{\pm
i\theta_j}, j = 1, \ldots, k$ (where $r_1e^{i\theta_1}$ is the one
already mentioned). Let $K$ be the polyhedral cone in ${\Bbb R}^n$
given by $K = \fn{pos}\{ x_1, \ldots, x_m \}$, where for $j = 1,
\ldots, m$,
$$
x_j = \left[ \begin{array}{c}r_1^{j-1}\cos
(j-1)\theta_1\\r_1^{j-1}\sin (j-1)\theta_1\\ \vdots \\
r_p^{j-1}\cos (j-1)\theta_p\\r_p^{j-1}\sin (j-1)\theta_p\\1
\end{array} \right].$$
\noindent It is clear that $K$ is a pointed cone.  A sufficient
condition for $K$ to be a full cone is that the $n \times n$
matrix whose $j$th column is $x_j$, for $j = 1, \ldots, n$, is
nonsingular. Upon pre-multiplying the latter matrix by the
 $n \times n$ matrix
$$\left[
\begin{array}{rr}1&i\\1&-i
\end{array}\right] \oplus \cdots \oplus \left[
\begin{array}{rr}1&i\\1&-i
\end{array}\right] \oplus [1],$$
\noindent we obtain the Vandermonde matrix
$$\left[
\begin{array}{cccc}1&z_1&\cdots&z_1^{n-1}\\1&\bar{z_1}&\cdots&\bar{z_1}^{n-1}
\\\vdots&\vdots&\vdots&\vdots\\1&z_p&\cdots&z_p^{n-1}\\1&\bar{z_p}&\cdots&\bar{z_p}^{n-1}\\1&1&\cdots&1
\end{array}\right], $$
where $z_j = r_je^{i\theta_j}$ for $j = 1, \ldots, p$, which is
nonsingular, as the roots of the polynomial $h(t)$ are simple (see
Lemma \ref{lemma6.2}).  So $K$ is a full cone.

Now take $A$ to be the $n \times n$ matrix
$$r_1\left[ \begin{array}{rr}\cos\theta_1&-\sin\theta_1\\
\sin\theta_1&\cos\theta_1 \end{array}\right] \oplus \cdots \oplus
r_p\left[ \begin{array}{rr}\cos\theta_p&-\sin\theta_p\\
\sin\theta_p&\cos\theta_p \end{array}\right] \oplus [1].$$
\noindent As can be readily checked, $Ax_j = x_{j+1}$ for $j = 1,
\ldots, m-1$.  Also, we have
$$
   Ax_m = \left[
\begin{array}{c}r_1^m\cos m\theta_1\\r_1^m\sin
m\theta_1\\ \vdots \\ r_p^m\cos m\theta_p\\r_p^m\sin m\theta_p\\1
\end{array} \right].$$
\noindent  But the assumption that $r_je^{i\theta_j}$ are roots of
$h(t)$ for $j = 1, \ldots, p$ implies that relations (\ref{okay1})
and (\ref{okay2}), which appear in the proof of Lemma
\ref{lemma6.2}, are satisfied for the pair $(r,\theta) =
(r_j,\theta_j)$ for every such $j$.  So we have $Ax_m =
(1-c)x_1+cx_2$.  Hence $A$ is $K$-nonnegative.  It remains to show
that $x_1, \ldots, x_m$ are precisely the pairwise distinct
extreme vectors of $K$ (the polyhedral cone generated by them) and
the face $\Phi(x_1+x_2)$ contains (up to multiples) only $x_1,x_2$
as its extreme vectors. Once this is done, it will follow that the
digraph $({\cal E},{\cal P})$ is given by Figure $1$.

For $j = 1, \ldots, m$, denote by $u_j$ the subvector of $x_j$
formed by its $1$st, $2$nd and last components. Since $\theta_1
\in (\frac{2\pi}{m},\frac{2\pi}{m-1})$ and $r_1 = r_{\theta_1}$,
by Theorem \ref{coro6}(iii) $u_1, \ldots, u_m$ are precisely
 the pairwise distinct extreme
vectors of the polyhedral cone $\fn{pos}\{u_1, \ldots, u_m \}$. So
it is clear that each $x_j$ cannot be written as a nonnegative
linear combination of the remaining $x_l$$^\prime$s or, in other
words, $x_1, \ldots, x_m$ are precisely the extreme vectors of
$K$.  Also, the proof of Theorem \ref{coro6}(iii) guarantees that
$u_1, u_2$ are neighborly extreme vectors of the $3$-dimensional
polyhedral cone $\fn{pos} \{ u_1, \ldots, u_m \}$, which means
that there is no representation of $u_1+u_2$ as a positive linear
combination of $u_1, \ldots, u_m$, in which at least one of the
vectors $u_3, \ldots, u_m$ is involved.  As a consequence, there
is also no representation of $x_1+x_2$ as a positive linear
combination of $x_1, \ldots, x_m$, in which at least one of the
vectors $x_3, \ldots, x_m$ is involved.  Hence, the face of $K$
generated by $x_1+x_2$ is $2$-dimensional, as desired.

Now we consider the problem of constructing the desired pair
$(K,A)$ for even $m$.   Choose any $c\in (0,1)$.   By Lemma
\ref{lemma6.2} the polynomial $h(t) = t^m-ct-(1-c)$ has a complex
root $r_1e^{i\theta_1}$ with $\theta_1\in
(\frac{2\pi}{m},\frac{2\pi}{m-1})$ and $0 < r_1 < 1$. Furthermore,
$h(t)$ has two distinct real roots, namely, $1$ and, say $a$. We
write $m$ as $2k+2$ and let the non-real complex roots of $h(t)$
be $r_je^{\pm i\theta_j} \mbox{ for } j = 1, \ldots, k$ (where
$r_1e^{i\theta_1}$ is the root just mentioned). Now write $n$ as
$2p+2$ or $2p+1$ (with $1 \le p \le k$), depending on whether $n$
is even or odd. Let $K$ be the polyhedral cone in ${\Bbb R}^n$
given by:
$$K = \fn{pos}\{ x_1, \ldots, x_m \},$$
where for $j = 1, \ldots, m,$
$$
x_j = \left[ \begin{array}{c}r_1^{j-1}\cos
(j-1)\theta_1\\r_1^{j-1}\sin (j-1)\theta_1\\ \vdots \\
r_p^{j-1}\cos (j-1)\theta_p\\r_p^{j-1}\sin (j-1)\theta_p\\
a^{j-1}\\1
\end{array} \right] \mbox{ or }
 \left[ \begin{array}{c}r_1^{j-1}\cos (j-1)\theta_1\\r_1^{j-1}\sin
(j-1)\theta_1\\ \vdots \\ r_p^{j-1}\cos
(j-1)\theta_p\\r_p^{j-1}\sin (j-1)\theta_p\\1
\end{array} \right],$$
depending on whether $n$ is even or odd.

Now take $A$ to be the $n\times n$ matrix
$$r_1\left[ \begin{array}{rr}\cos\theta_1&-\sin\theta_1\\
\sin\theta_1&\cos\theta_1 \end{array}\right] \oplus \cdots \oplus
r_p\left[ \begin{array}{rr}\cos\theta_p&-\sin\theta_p\\
\sin\theta_p&\cos\theta_p \end{array}\right] \oplus [a] \oplus
[1]$$ or
$$r_1\left[ \begin{array}{rr}\cos\theta_1&-\sin\theta_1\\
\sin\theta_1&\cos\theta_1 \end{array}\right] \oplus \cdots \oplus
r_p\left[ \begin{array}{rr}\cos\theta_p&-\sin\theta_p\\
\sin\theta_p&\cos\theta_p \end{array}\right] \oplus [1],$$ again
depending on whether $n$ is even or odd. Using the same argument
as before, we can show that $K$ is a proper cone, $x_1, \ldots,
x_m$ are its extreme vectors, $A$ is $K$-primitive and $({\cal
E},{\cal P})$ is given by Figure $1$.

 When $m$ is odd and
$n$ is even, take any $c\in (c_m,1)$.  According to Lemma
\ref{lemma6.2}, the polynomial $h(t) = t^m-ct-(1-c)$ has simple
roots, three of which are real, namely, $1$ and say $a_1,a_2$,
with $a_1 < a_2$, where $a_1,a_2\in (-1,0)$.  Let the remaining
$2k$ (where $m-3 = 2k$) non-real roots be $r_je^{\pm i\theta_j}, j
= 1, \ldots, k$. (Note that in this case we have $m \ge 5$.)  By
the same lemma, we may assume that $\theta_1\in
(\frac{2\pi}{m},\frac{2\pi}{m-1})$. Write $n$ as $2p+2$. Then
$1\le p \le k$. Let $K$ be the polyhedral cone in ${\Bbb R}^n$
given by:
$$K = \fn{pos}\{ x_1, \ldots, x_m \},$$
where for $j = 1, \ldots, m,$
$$
x_j = \left[ \begin{array}{c}r_1^{j-1}\cos
(j-1)\theta_1\\r_1^{j-1}\sin (j-1)\theta_1\\ \vdots \\
r_p^{j-1}\cos (j-1)\theta_p\\r_p^{j-1}\sin (j-1)\theta_p\\
a_1^{j-1}\\1
\end{array} \right].
$$
Now let $A$ be the matrix
$$r_1\left[ \begin{array}{rr}\cos\theta_1&-\sin\theta_1\\
\sin\theta_1&\cos\theta_1 \end{array}\right] \oplus \cdots \oplus
r_p\left[ \begin{array}{rr}\cos\theta_p&-\sin\theta_p\\
\sin\theta_p&\cos\theta_p \end{array}\right] \oplus [a_1] \oplus
[1].$$ \noindent The subsequent arguments are similar to those for
the previous cases. We omit the details.  \hfill$\blacksquare$

Our next step is to show the following crucial lemma.

\begin{lemma}\label{lemmacrucial}  Let $m,n$ be positive integers such that $3\le n
\le m$, where $n$ is odd or $n,m$ are both even.  Let $K_0$ be the
cone in ${\Bbb R}^n$ with extreme vectors $y_1, \ldots, y_m$ given
by {\rm:}
$$ y_j=\left[\begin{array}{c}
\cos (j-1)\frac{2\pi}{m} \\
\sin (j-1)\frac{2\pi}{m} \\
\cos (j-1)\frac{4\pi}{m} \\
\sin (j-1)\frac{4\pi}{m} \\
\vdots \\
\cos (j-1)\frac{(n-1)\pi}{m} \\
\sin (j-1)\frac{(n-1)\pi}{m} \\
1
\end{array}\right], 1\le j \le m,
$$
when $n$ is odd {\rm;} and by
$$
y_j=\left[\begin{array}{c}
\cos (j-1)\frac{2\pi}{m} \\
\sin (j-1)\frac{2\pi}{m} \\
\cos (j-1)\frac{4\pi}{m} \\
\sin (j-1)\frac{4\pi}{m} \\
\vdots \\
\cos (j-1)\frac{(n-2)\pi}{m} \\
\sin (j-1)\frac{(n-2)\pi}{m} \\
(-1)^{j-1}\\1
\end{array}\right], 1\le j \le m,
$$
when $n,m$ are both even. Then $\sum_{j=1}^{n-1}y_j$ lies in
$\partial K_0$ and generates a simplicial face of $K_0$.
\end{lemma}

Before we consider the proof of Lemma \ref{lemmacrucial}, we show
how it can be used to finish the proof of Theorem
\ref{maintheorem}.

We first treat the case when $m$ is even or $m,n$ are both odd.
Recall the proof of Lemma \ref{theorem7.1}. In constructing the
polyhedral cone $K$ and the $K$-primitive matrix $A$, certain real
numbers $r_j$'s and $\theta_j$'s are involved. For the argument to
work there, $r_je^{\pm i\theta_j}$ ($j = 1, \ldots, p$, where $p =
\frac{n-1}{2} \mbox{ or } \frac{n-2}{2}$, depending on whether $n$
is odd or even) can be any $p$ pairs of non-real conjugate complex
roots of the polynomial $h(t)$, but $\theta_1$ has to be chosen
from $(\frac{2\pi}{m},\frac{2\pi}{m-1})$. Now in the proof of
Theorem \ref{maintheorem} in order that the argument works we need
to be more careful about the choices of the $r_j$'s and
$\theta_j$'s. Since $h(t)$ tends to the polynomial $t^m-1$ as $c$
tends to $0$ (from the right), when $c$ is sufficiently close to
$0$, for each $j = 1, \ldots, \frac{m}{2}-1$, there is precisely
one root of $h(t)$ within an $\varepsilon$-neighborhood of
$e^{\pm\frac{2\pi j}{m}i}$, where $\varepsilon$ is a sufficiently
small positive number such that the $\varepsilon$-neighborhoods of
$e^{\pm\frac{2\pi j}{m}i}$ for $j = 1, \ldots, \frac{m}{2}$ are
pairwise mutually disjoint. We label that particular root as
$r_je^{i\theta_j}$. Then it is readily seen that the extreme
vectors $x_1, \ldots, x_m$ of $K$ can be made as close as we
please to, respectively, the extreme vectors $y_1, \ldots, y_m$ of
$K_0$.  By Lemma \ref{lemmacrucial} the vector $y_1+y_2+ \cdots +
y_{n-1}$ lies in $\partial K_0$ and generates a simplicial face of
$K_0$. The same is also true for the vector $y_{m-n+2}+y_{m-n+3}+
\cdots + y_m$ of $K_0$, as it is clear that there is an
automorphism of $K_0$ that takes $y_j$ to $y_{j+m-n+1}$ for $j =
1, \ldots, m$ (where for $r
> m, y_r$ is taken to be $y_s$ with $s \equiv r (\fn{mod} m), 1\le
s \le m$). Hence $\fn{span}\{ y_{m-n+2}, y_{m-n+3}, \ldots,
y_{m}\}$ is a supporting hypersubspace for $K_0$ and the remaining
extreme vectors $y_1, \ldots, y_{m-n+1}$ all lie in the same open
half-space determined by this hypersubspace. By continuity, it
follows that when $c$ is sufficiently close to $0$, $\fn{span}\{
x_{m-n+2}, x_{m-n+3}, \ldots, x_m \}$ is a supporting
hypersubspace for $K$ and the remaining extreme vectors of $K$ all
lie on the same open half-space determined by this hypersubspace.
The latter implies that $x_{m-n+2}+ x_{m-n+3}+ \cdots + x_m \in
\partial K$, which is what we want.

When $m$ is odd and $n$ is even, we still use the construction for
$K$ and $A$ as given in the proof of Lemma \ref{theorem7.1}.
However, instead of $K_0$ we make use of the polyhedral cone $K_1$
of ${\Bbb R}^n$ with extreme vectors $y_1, \ldots, y_{m-1}$ given
by:
$$ y_j=\left[\begin{array}{c}
\cos (j-1)\frac{2\pi}{m-1} \\
\sin (j-1)\frac{2\pi}{m-1} \\
\cos (j-1)\frac{4\pi}{m-1} \\
\sin (j-1)\frac{4\pi}{m-1} \\
\vdots \\
\cos (j-1)\frac{(n-2)\pi}{m-1} \\
\sin (j-1)\frac{(n-2)\pi}{m-1} \\
(-1)^{j-1} \\1
\end{array}\right], 1\le j \le m-1.
$$
In defining the $x_j$'s and the matrix $A$, we choose the
$\theta_j$'s ($j = 1, \ldots, \frac{n-2}{2}$) in such a way that
$r_je^{\pm i\theta_j}$ is the root of $h(t)$ closest to $e^{\pm
\frac{2\pi j}{m-1}i}$, noting that $r_je^{i\theta_j},
r_je^{-i\theta_j}, j = 1, \ldots, \frac{m-3}{2}, a_1, a_2$ and
$1$, which are the roots of the polynomial $h(t)$, approach to,
respectively, the roots $e^{\frac{2\pi j}{m-1}i}, e^{-\frac{2\pi
j}{m-1}i}, j = 1, \ldots, \frac{m-3}{2}, -1,0$ and $1$ of the
polynomial $t^m-t$, as $c$ approach to $1$ (from the left). Then
$x_j$ tends to $ y_j$ for $j = 1, \ldots, m-1$ and also $x_m$
tends to $y_1$ as $c$ tends to $1$.
 Since $m-1$ is even, by what we have just
done before (for the case when $m$ is even), the vector
$y_{m-n+1}+ y_{m-n+2}+ \ldots +y_{m-1}$ generates an
$(n-1)$-dimensional simplicial face of $K_1$.  So $\fn{span}\{
y_{m-n+1}, y_{m-n+2}, \ldots, y_{m-1} \}$ is a supporting
hypersubspace for $K_1$ and the remaining extreme vectors $y_1,
\ldots, y_{m-n}$ all lie in the same open half-space determined by
this hypersubspace. By continuity, when $c$ is sufficiently close
to $1$, $\fn{span}\{ x_{m-n+1}, x_{m-n+2}, \ldots, x_{m-1} \}$ is
a supporting hypersubspace for $K$ and the remaining extreme
vectors of $K$ all lie in the same open half-space determined by
this hypersubspace (noting that $x_m$ also lies in this open
half-space as $x_m,x_1$ both tend to $y_1$).  Hence, $x_{m-n+1}+
x_{m-n+2} + \ldots + x_{m-1}\in \partial K$, as desired.\\

Now we come to the proof of Lemma \ref{lemmacrucial}.  The proof
is fairly long and takes several steps.

The lemma clearly holds for the case $m = n$.  Hereafter we assume
that $m > n$.

We want to show that $\sum_{j=1}^{n-1}y_j$ lies in $\partial K_0$
and generates a simplicial face.
 For the purpose, it suffices to find a vector $v$
such that $\langle v,y_j \rangle$ equals $0$ for $j=1,\ldots,n-1$
and is positive for $j = n, \ldots, m$.  Note that the vectors
$y_1,y_2,\ldots,y_{n-1}$ are linearly independent, because by the
argument given in the proof of Lemma \ref{theorem7.1} one can show
that the $n\times n$ matrix
\[
 P := \left[\begin{array}{cccc}
y_1 & y_2 & \cdots & y_n
\end{array}\right]
\]
is nonsingular.  So the desired vector $v$ is, up to a positive
scalar multiple, uniquely determined.
  Let $\tilde{v} = (C_{1n}, \ldots,
C_{nn})^T$, where $C_{ln}$ denotes the $(l,n)$-cofactor of $P$. By
elementary properties of the determinant function, for each $p$,
we have $\langle \tilde{v},y_p \rangle =  \det Q_p$, where $Q_p$
denotes the $n\times n$ matrix obtained from $P$ by replacing its
$n$th column by $y_p$.  Note that $\langle v,y_p \rangle = 0$ for
$p = 1, \ldots, n-1$. It remains to show that for $p = n,n+1,
\ldots, m, \det Q_p$ are all nonzero and have the same sign
--- if $\det Q_p$'s are all positive, take $v = \tilde{v}$; if $\det Q_p$'s are all
negative, take $v = -\tilde{v}$.

We find it more convenient to work in the complex domain.  Denote
$e^{2\pi i/m}$ by $\omega$. Pre-multiplying the matrix $Q_p$ by
\[
\underbrace{\left[\begin{array}{rr} 1 & i \\ 1 & -i
\end{array}\right]\oplus \left[\begin{array}{rr} 1 & i \\ 1 & -i
\end{array}\right]\oplus \cdots \oplus \left[\begin{array}{rr} 1 &
i \\ 1 & -i \end{array}\right]}_{\frac{n-2}{2} \mbox{ times}}
\oplus I_2
\] when $n$ is even, or by
\[
\underbrace{\left[\begin{array}{rr} 1 & i \\ 1 & -i
\end{array}\right]\oplus \left[\begin{array}{rr} 1 & i \\ 1 & -i
\end{array}\right]\oplus \cdots \oplus \left[\begin{array}{rr} 1 &
i \\ 1 & -i \end{array}\right]}_{\frac{n-1}{2} \mbox{ times}}
\oplus (1)
\] when $n$ is odd, we obtain the matrix
\[
\left[\begin{array}{cccccc}
1 & \omega & \omega^2 & \cdots & \omega^{n-2} & \omega^{p-1} \\
1 & \bar{\omega} & \bar{\omega}^2 & \cdots & \bar{\omega}^{n-2} & \bar{\omega}^{p-1} \\
1 & \omega^2 & \omega^4 & \cdots & \omega^{2(n-2)} & \omega^{2(p-1)} \\
1 & \bar{\omega}^2 & \bar{\omega}^4 & \cdots & \bar{\omega}^{2(n-2)} & \bar{\omega}^{2(p-1)} \\
\vdots & \vdots & \vdots & & \vdots & \vdots \\
1 & \omega^\frac{n-2}{2} & \omega^{n-2} & \cdots & \omega^\frac{(n-2)^2}{2} & \omega^{\frac{(n-2)(p-1)}{2}} \\
1 & \bar{\omega}^\frac{n-2}{2} & \bar{\omega}^{n-2} & \cdots & \bar{\omega}^\frac{(n-2)^2}{2} & \bar{\omega}^{\frac{(n-2)(p-1)}{2}} \\
1 & -1 & 1 & \cdots & (-1)^{n-2} & (-1)^{p-1} \\
1 & 1 & 1 & \cdots & 1 & 1
\end{array}\right]
\]
when $m,n$ are both even, or the matrix
\[
\left[\begin{array}{cccccc}
1 & \omega & \omega^2 & \cdots & \omega^{n-2} & \omega^{p-1} \\
1 & \bar{\omega} & \bar{\omega}^2 & \cdots & \bar{\omega}^{n-2} & \bar{\omega}^{p-1} \\
1 & \omega^2 & \omega^4 & \cdots & \omega^{2(n-2)} & \omega^{2(p-1)} \\
1 & \bar{\omega}^2 & \bar{\omega}^4 & \cdots & \bar{\omega}^{2(n-2)} & \bar{\omega}^{2(p-1)} \\
\vdots & \vdots & \vdots & & \vdots & \vdots \\
1 & \omega^\frac{n-1}{2} & \omega^{n-1} & \cdots & \omega^{\frac{(n-1)(n-2)}{2}} & \omega^{\frac{(n-1)(p-1)}{2}} \\
1 & \bar{\omega}^\frac{n-1}{2} & \bar{\omega}^{n-1} & \cdots & \bar{\omega}^{\frac{(n-1)(n-2)}{2}} & \bar{\omega}^{\frac{(n-1)(p-1)}{2}} \\
1 & 1 & 1 & \cdots & 1 & 1
\end{array}\right]
\]
when $n$ is odd.  Denote the determinant of the matrix by $e_p$.

Note that for each $p = n,n+1, \ldots, m, e_p$ is equal to $\det
Q_p$ times a nonzero constant, which depends on $n$ but not on
$p$.  More specifically, the said nonzero constant is $\pm
2^{\frac{n-2}{2}} \mbox{ if } n \equiv 2$ (mod 4), $\pm
2^{\frac{n-2}{2}}i \mbox{ if } n \equiv 0$ (mod 4), $\pm
2^{\frac{n-1}{2}} \mbox{ if } n \equiv 1$ (mod 4), and $\pm
2^{\frac{n-1}{2}}i \mbox{ if } n \equiv 3$ (mod 4).  We want to
show all $e_p$s are nonzero real numbers with the same sign if
$n\equiv 2$ (mod 4) or $n\equiv 1$ (mod 4), and if $n\equiv 0$
(mod 4) or $n\equiv 3$ (mod 4) then all $ie_p$s are also nonzero
real numbers with the same sign.

We now consider a kind of generalized Vandermonde determinant on
the indeterminates $t_1, \ldots, t_n$.  Let $n\ge 3$. For every
integer $p \ge n-1$, let $f_p(t_1, \ldots, t_n)$ denote the
polynomial function
\[ \extrarowheight=3pt
f_p(t_1, \ldots, t_n)=\det\left[\begin{array}{cccccc}
1 & t_1 & t^2_1 & \cdots & t^{n-2}_1 & t^p_1 \\
1 & t_2 & t^2_2 & \cdots & t^{n-2}_2 & t^p_2 \\
1 & t_3 & t^2_3 & \cdots & t^{n-2}_3 & t^p_3 \\
\vdots & \vdots & \vdots & & \vdots & \vdots \\
1 & t_n & t^2_n & \cdots & t^{n-2}_n & t^p_n
\end{array}\right].
\]

\noindent By the $k$th {\it complete {\rm(}homogeneous{\rm)}
symmetric polynomial} in $t_1, \ldots,t_n$, denoted by $h_k(t_1,
t_2, \ldots, t_n)$, we mean the polynomial which is the sum of all
monomials in $t_1, \ldots, t_n$ of degree $k$, where each monomial
appears exactly once.  For instance,
$$h_2(t_1,t_2,t_3) = t_1t_2+t_1t_3+t_2t_3+t_1^2+t_2^2+t_3^2.$$
By definition, $h_0(t_1, \ldots, t_n) \equiv 1$.  It is convenient
to define $h_r(t_1, \ldots, t_n)$ to be $0$ for $r < 0$.

\begin{claim}\rm  For every integer $p \ge n-1$, we have
\begin{equation}
f_p(t_1,t_2,\ldots,t_n)=h_{p+1-n}(t_1,t_2,\ldots,t_n)\prod_{1\le
i< j\le n} (t_j-t_i). \label{eqn3}
\end{equation}
\end{claim}

\noindent The above claim can be deduced from a formula that
expresses a Schur function $s_{\lambda}(t_1, \ldots, t_n)$ as a
determinant involving complete symmetric polynomials $h_r(t_1,
\ldots, t_n)$, known as the {\it Jacobi-Trudi determinant} (see
Macdonald \cite [p.25] {Mac} or Sagan \cite[p.154-159]{Sag}).
Recall that if $\lambda = (\lambda_1, \lambda_2, \ldots)$ is a
{\it partition} of length $\le n$ (i.e. a finite or infinite
sequence of nonnegative integers in non-increasing order such that
the number of nonzero $\lambda_i$ is at most $n$) then the
quotient
$$ \frac{\fn{det}[t_i^{j-1+\lambda_{n-j+1}}]_{1\le i,j \le
n}}{\prod_{1\le i < j \le n}(t_j-t_i)}$$ is a symmetric polynomial
in $t_1, \ldots, t_n$.  It is called the {\it Schur function
associated with} $\lambda$ and is denoted by $s_{\lambda}(t_1,
\ldots, t_n)$.  For $p \ge n-1$, the polynomial $f_p(t_1,
\ldots,t_n)$ introduced above satisfies
$$ \frac{f_p(t_1, \ldots, t_n)}{\prod_{1\le i < j \le n}(t_j-t_i)} =
s_{\lambda}(t_1, \ldots, t_n),$$ where $\lambda$ is the partition
$(p-(n-1),0, \ldots)$.  According to the said determinantal
formula, for any partition $\lambda$ of length $\le n$, we have
$$ s_{\lambda}(t_1, \ldots, t_n) =
\fn{det}[h_{\lambda_i-i+j}(t_1, \ldots, t_n)]_{1\le i,j \le n}.$$
A little calculation shows that, when $\lambda = (p-(n-1),0,
\ldots), [h_{\lambda_i-i+j}]_{1\le i,j \le n}$ is an upper
triangular matrix whose $(1,1)$-entry is $h_{p-(n-1)}$ and all of
whose other diagonal entries are equal to $1$.  It follows that we
have $s_{\lambda}(t_1, \ldots, t_n) = h_{p+1-n}(t_1, \ldots,
t_n)$, as desired. \hfill$\blacksquare$\\

For any set $S$ of $n$ complex numbers, we denote by $h_j(S)$ the
value of $h_j(t_1, \ldots, t_n)$ evaluated at an $n$-tuple formed
by all of the members of $S$, taken in any order.  Similarly, we
use $\sigma_j(S)$ to denote the {\it $j$th elementary symmetric
function} on $S$.

\begin{claim}\rm Let $S$ be a nonempty proper subset of ${\Bbb Z}_m$, the set of
all $m$th roots of unity.  Denote by $S^c$ the complement of $S$
in ${\Bbb Z}_m$.  Then
$$h_j(S) = (-1)^j\sigma_j(S^c)$$
for $j = 1, \ldots, m-n$, where $n$ is the cardinality of $S$.
\end{claim}

{\it Proof}.  For any set $T$ of $n$ complex numbers, as is known
(see, for instance, \cite {Mac}), the generating function for the
elementary symmetric functions on $T$ is given by:
$$E(t;T) = \sum_{r=0}^n \sigma_r(T)t^r = \prod_{x\in T}(1+xt),$$
(where $\sigma_0(T)$ is taken to be $1$).  Also, the generating
function for the complete symmetric polynomials on $T$ is given
by:
$$H(t;T) = \sum_{r=0}^\infty h_r(T)t^r = \prod_{x\in
T}\frac{1}{1-xt}.$$

\noindent In view of the relation $\prod_{x\in {\Bbb Z}_m}(1-xt) =
1-t^m$, for the given set $S$, we have \[
\sum_{r=0}^{m-n}(-1)^r\sigma_r(S^c)t^r = E(-t,S^c) = \prod_{x\in
S^c}(1-xt) = \frac{1-t^m}{\prod_{x\in S}(1-xt)} =
(1-t^m)\Sigma_{r=0}^\infty h_r(S)t^r. \]
 \noindent By comparing
the coefficients, our claim follows. \hfill$\blacksquare$\\

In view of Claim 1 and the discussion preceding the claim, it
remains to show the following:

\begin{claim}\rm
For $r = 1, \ldots, m-n$,
 $h_r(S) > 0$, where $S$ is the subset of the set of $m$th roots of unity given
 by:
$$S =
\left\{\begin{array}{l}\{\omega,\bar{\omega},\ldots,
\omega^{\frac{n-2}{2}},\bar{\omega}^{\frac{n-2}{2}},-1,1\}
\mbox{ when }m,n\mbox{ are both even }\\
\{\omega,\bar{\omega}, \ldots,
\omega^{\frac{n-1}{2}},\bar{\omega}^{\frac{n-1}{2}},1\} \mbox{
when } n \mbox{ is odd. }\end{array}\right . $$
\end{claim}

\noindent Now, by Claim 2 we have
\[
\prod_{x\in S^c}(t-x) =
t^{m-n}+\sum_{j=1}^{m-n}(-1)^j\sigma_j(S^c)t^{m-n-j} =
t^{m-n}+\sum_{j=1}^{m-n}h_j(S)t^{m-n-j}.
\]
To establish Claim 3, it suffices to show that the coefficients of
the polynomial $\prod_{x\in S^c}(t-x)$ are all positive. We
complete our argument by applying the following interesting
nontrivial result due to Barnard et al.$\,$\cite [Theorem
1]{B--D--P--W}. (Or see \cite [Theorem 2.4.5] {B--E} or \cite
[Theorem 4]{E--G} ).

\begin{thm} Let $p(t)$ be a polynomial of degree $n$, $p(0) = 1$, with
nonnegative coefficients and zeros $a_1, \ldots, a_n$.  For $\tau
\ge 0$ write
$$p_{\tau}(t) = \prod_{1\le j \le n\atop |\fn{arg}(a_j)| >
\tau}(1-t/a_j),$$ where $\fn{arg}(z)$ is defined so that
$\fn{arg}(z)\in [-\pi,\pi)$.  If $p_{\tau}(t) \ne p(t)$, then the
coefficients of $p_{\tau}(t)$ are all positive.
\end{thm}

Let us consider the case when $m, n$ are both even first. In this
case we have $S^c =
\{\omega^{\frac{n}{2}},\bar{\omega}^{\frac{n}{2}}, \ldots,
\omega^{\frac{m}{2}-1}, \bar{\omega}^{\frac{m}{2}-1} \}$.Take
$p(t)$ to be the following polynomial, which has nonnegative
coefficients and constant term $1$:
$$ \prod_{1\le j\le m-1 \atop j\ne \frac{m}{2}}(t-\omega^j) = \frac{t^m-1}{(t-1)(t+1)} = t^{m-2}+t^{m-4}+ \cdots +t^2+1.$$

\noindent Choose any positive number $\tau$ from
$(\frac{n-2}{m}\pi,\frac{n}{m}\pi)$.  Then
$$p_{\tau}(t) = \prod_{j=\frac{n}{2}}^{\frac{m}{2}-1}
(1-\frac{t}{\omega^j})(1-\frac{t}{\bar{\omega}^j}) =
\prod_{j=\frac{n}{2}}^{\frac{m}{2}-1}(t-\omega^j)(t-\bar{\omega}^j).$$
By Theorem D, $p_{\tau}(t)$ is a polynomial with positive
coefficients.   So $\prod_{x\in S^c}(t-x)$ equals $p_{\tau}(t)$
and is a polynomial with positive coefficients.

When $m$ is even and $n$ is odd, we take $p(t)$ to be the same
polynomial, but choose $\tau$ from
$(\frac{n-1}{m}\pi,\frac{n+1}{m}\pi)$.  A little calculation shows
that in this case $\prod_{x\in S^c}(t-x)$ is equal to
$(t+1)p_{\tau}(t)$ and so it also has positive coefficients.

When $m, n$ are both odd, we take $p(t)$ to be the polynomial
$t^{m-1}+t^{m-2}+ \cdots + 1$ and apply a similar argument.

The proof for Theorem \ref{maintheorem} is complete.\hfill$\blacksquare$\\

{\bf Conjecture 7.4}. For any positive integers $m,n, 3 \le n \le
m$,
$$\fn{max}\{ \gamma(K): K\in {\cal P}(m,n)\} = (n-1)(m-1),$$ when $m$
is odd and $n$ is even.\\

By Theorem \ref{theorem2}(I) the Conjecture is confirmed for the
minimal cone case.

\setcounter{equation}{0}
\section{Uniqueness of exp-maximal cones and their exp-maximal primitive matrices}

Given positive integers $m,n$ with $3\le n \le m$, up to linear
isomorphism, how many exp-maximal cones are there in ${\cal
P}(m,n)$ ? For a given exp-maximal cone $K$ in ${\cal P}(m,n)$, up
to cone-equivalence modulo positive scalar multiplication, how
many exp-maximal $K$-primitive matrices are there ?  In this
section we are going to address these questions for the cases $m =
n, m = n+1$ and $n = 3$.

Since there is, up to linear isomorphism, only one simplicial cone
of a given dimension, we need not treat the problem of identifying
exp-maximal cones in ${\cal P}(m,n)$ for the case $m = n$. The
problem of identifying exp-maximal minimal cones has already been
carried out in Section $5$.  According to Theorem \ref{theorem2},
a cone $K\in {\cal P}(n+1,n)$ is exp-maximal if and only if $K$ is
an indecomposable minimal cone with a balanced relation for its
extreme vectors or $n$ is even and $K$ is the direct sum of a ray
and an $(n-1)$-dimensional indecomposable minimal cone with a
balanced relation for its extreme vectors. But, by Theorem
\ref{theorem2.5} for every positive integer $n \ge 3$, there is,
up to linear isomorphism, only one $n$-dimensional indecomposable
minimal cone with a balanced relation for its extreme vectors, so
$n$-dimensional exp-maximal minimal cones are known completely: up
to linear isomorphism, there are one such cone when $n$ is odd and
two such cones when $n$ is even.  Some work on identifying
exp-maximal $3$-dimensional cones has also been done in Section
$6$.  By Theorem \ref{coro6}(iii), for every positive integer $m
\ge 3$, there are uncountably infinitely many exp-maximal cones
$K_{\theta}$ in ${\cal P}(m,3)$, one for each
$\theta\in(\frac{2\pi}{m},\frac{2\pi}{m-1})$.  Later in this
section, we will show that when $m \ge 4$, up to linear
isomorphism, the $K_{\theta}$'s are precisely all the exp-maximal
cones in ${\cal P}(m,3)$ and, moreover, when $m \ge 6$, for
different values of $\theta$, the corresponding $K_{\theta}$'s are
not linearly isomorphic.  (At present, we {\it do not} know
whether the latter assertion can be extended to cover the case $m
= 5$.)

Once the exp-maximal cones in ${\cal P}(m,n)$ have been
identified, the next natural problem to consider is to identify,
for each typical exp-maximal cone $K$, all the exp-maximal
$K$-primitive matrices.  So we will also identify, up to
cone-equivalence and scalar multiples, the exp-maximal primitive
matrices for exp-maximal minimal cones, for exp-maximal
$3$-dimensional cones, as well as for simplicial cones.  As a
matter of fact, in identifying exp-maximal minimal cones in
Section $5$, we have already provided (implicitly in Lemma
\ref{lemma 5.1}), up to cone-equivalence and scalar multiples, all
the exp-maximal primitive matrices for minimal cones.  Also, as we
will show, for $m \ge 4$, if $K\in {\cal P}(m,3)$ is an
exp-maximal polyhedral cone and $A$ is an exp-maximal
$K$-primitive matrix then, up to scalar multiples, $A$ is
cone-equivalent to $A_{\theta}$ for some $\theta\in
(\frac{2\pi}{m},\frac{2\pi}{m-1})$, where $A_{\theta}$ is the
matrix introduced in Theorem \ref{coro6}, and moreover, provided
that $m \ge 6$, $A_{\theta}$ is (up to scalar multiples) the only
exp-maximal $K_{\theta}$-primitive matrix. However, there are
exp-maximal primitive matrices for the cone ${\Bbb R}^3_+$ that
are (up to cone-equivalence and scalar multiples) not of the form
$A_{\theta}$, where $\theta\in (\frac{2\pi}{3},\pi)$.

 We deal with
the minimal cone case first. We begin with a result, which is true
not only for the minimal cone case.

\begin{lemma}\label{lemma8.1}  Let $K\in {\cal P}(m,n)$ be indecomposable.  If $A, \tilde{A}$ are different $K$-nonnegative
matrices such that the digraphs $({\cal E},{\cal P}(A,K))$ and
$({\cal E},{\cal P}(\tilde{A},K))$ are given both by Figure $1$ or
both by
 Figure $2$, then $\tilde{A}$ and $A$ are not cone-equivalent.
\end{lemma}

{\it Proof}.  Let $A$ and $\tilde{A}$ be cone-equivalent
$K$-nonnegative matrices such that the digraphs $({\cal E},{\cal
P}(A,K))$ and $({\cal E},{\cal P}(\tilde{A},K))$ are given both by
Figure $1$ or both by Figure $2$.  Then there exists an
automorphism $P$ of $K$ such that $P\tilde{A} = AP$.  We contend
that for $j = 1, \ldots, m, P$ maps $x_j$ to a positive multiple
of itself. Once this is proved, it will follow that $P\in
\Phi(I)$. But $K$ is indecomposable, so by \cite [Theorem
3.3]{L--S} $\Phi(I)$ is an extreme ray of $\pi(K)$; hence $P$ is a
positive multiple of $I$ and we have $\tilde{A} = A$, as desired.

To prove our contention, we first deal with the case when the
digraphs $({\cal E},{\cal P}(A,K))$ and $({\cal E},{\cal
P}(\tilde{A},K))$ are both given by Figure $1$.  As $P$ is an
automorphism of $K$, $P$ permutes the extreme rays of $K$ among
themselves.  According to the proof of Lemma \ref{lemma4.2}(ii),
the maximum value of $\gamma(A,x)$, for $x = x_1, \ldots, x_m$, is
attained at $x_1$ only.  When $A$ is replaced by $\tilde{A}$, the
same can be said.  Since $A$ and $\tilde{A}$ are cone-equivalent,
by Fact \ref{fact2.5}(v) $P$ must map $x_1$ to a positive multiple
of itself. Making use of the relation $P\tilde{A} = AP$ and the
fact that $Ax_i$ is a positive multiple of $x_{i+1}$ for $i = 1,
\ldots, m-1$ and proceeding inductively, we readily show that $P$
maps each $x_i$ to a positive multiple of itself, which is our
contention.

If $({\cal E},{\cal P}(A,K))$ and $({\cal E},{\cal
P}(\tilde{A},K))$ are both given by Figure $2$, we readily show
that $Px_2$ must be a positive multiple of itself and then we can
proceed in a similar manner.  \hfill$\blacksquare$

Next, a remark on the automorphisms of an indecomposable minimal
cone is in order.

Let $K$ be an indecomposable minimal cone in ${\Bbb R}^n$ with
extreme vectors $x_1, \ldots, x_{n+1}$ that satisfy the relation
$$ x_1+ \cdots + x_p = x_{p+1}+ \cdots + x_{n+1}.$$
Let $\sigma$ be a permutation on $\{1, \ldots, n+1\}$ that maps
$\{ 1, \ldots, p\}$ and $\{ p+1, \ldots, n\}$ each onto itself, or
interchanges the first set with the second set (in which case $n$
is odd and $p = \frac{n+1}{2}$). Then $\sigma$ determines a
(unique) automorphism $P_{\sigma}$ of $K$ with spectral radius $1$
($x_1+ \cdots + x_p$ being the corresponding eigenvector) which is
given by: $P_{\sigma}x_j = x_{\sigma(j)}$ for $j = 1, \ldots,
n+1$. Conversely, every automorphism of $K$ whose spectral radius
is $1$ arises in this way.

Now we consider the exp-maximal primitive matrices for
indecomposable exp-maximal minimal cones first.  We give two
results, one for the odd dimensional case and the other for the
even dimensional case.  The relations {\rm(\ref{rel2})},
{\rm(\ref{eq2})}, {\rm(\ref{rel3})}, {\rm(\ref{eq3})},
{\rm(\ref{rel5})}, {\rm(\ref{eq5})} that are mentioned in these
results have already appeared in Lemma \ref{lemma 5.1}.

\begin{theorem}\label{theorem8.2} Let $K$ be an $n$-dimensional indecomposable exp-maximal minimal
cone, where $n$ is odd.  Suppose that the extreme vectors $x_1,
\ldots, x_{n+1}$ of $K$ satisfy relation {\rm(\ref{rel2})} with
{\rm(}with $m = n+1${\rm)}.  For every $\alpha
> 0$, let $A_{\alpha}$ be the exp-maximal $K$-primitive matrix
given by {\rm(\ref{eq2})} {\rm(}but with $A$ replaced by
$A_{\alpha}${\rm)}.  Then\,{\rm:}
\begin{enumerate}

\item[{\rm(i)}]  $\Phi(A_{\alpha})$ is a $2$-dimensional face,
independent of the choice of the positive scalar $\alpha${\rm;}
its relative interior consists of positive multiples of matrices
of the form $A_{\tilde{\alpha}}$.
 \item[{\rm(ii)}] Every
exp-maximal $K$-primitive matrix is cone-equivalent to a positive
multiple of some $A_{\alpha}$ and thus is a positive multiple of a
matrix of the form $P_{\sigma}^{-1}A_{\alpha}P_{\sigma}$, where
$P_{\sigma}$ is the automorphism of $K$ given by $P_{\sigma}x_j =
x_{\sigma(j)}$ for $j = 1, \ldots, n+1$, $\sigma$ being a
permutation on $\{1, \ldots, n+1\}$ that maps $\{1,3, \ldots,
n-2,n\}$ onto itself or onto $\{2,4, \ldots, n-1,n+1\}$.
 \item[{\rm(iii)}] For distinct positive
scalars $\alpha_1,\alpha_2$, $A_{\alpha_1}$ and $A_{\alpha_2}$ or
their positive multiples are not cone-equivalent.
\end{enumerate}
\end{theorem}

{\it Proof} (i)  It is clear that $\Phi(A_{\alpha}) =
\Phi(A_{\tilde{\alpha}})$ for any $\alpha, \tilde{\alpha} > 0$ as
$\Phi(A_{\alpha}x_j) = \Phi(A_{\tilde{\alpha}}x_j)$ for each $j$
and $K$ is polyhedral.  We are going to show that
$\Phi(A_{\alpha})$ is equal to the $2$-dimensional face generated
by the extreme matrices $B, C$ determined respectively by:
$$Bx_i = x_{i+1} \mbox{ for } i = 1, \ldots, m,
\mbox{ where }x_{m+1} \mbox{ is taken to be } x_1,$$ and
$$Cx_1 = x_2 = Cx_m,  Cx_i = 0 \mbox{ for } i = 2,  \ldots,
m-1.$$ It is readily checked that $B$ and $C$ each preserve
relation (\ref{rel2}); so they are well-defined and
$K$-nonnegative.  Since $A_{\alpha} = B+\alpha C$, we have $B,
C\in \Phi(A_{\alpha})$ and hence $\fn{pos}\{ B,C\} \subseteq
\Phi(A_{\alpha})$.  To complete the argument, we contend that
every matrix in $\fn{ri}\Phi(A_{\alpha})$ is a positive multiple
of some $A_{\tilde{\alpha}}$ (and hence belongs to $\fn{pos}\{ B,C
\}$). Once this is proved, the desired reverse inclusion follows
as $\Phi(A_{\alpha}) = \fn{cl}[\fn{ri} \Phi(A_{\alpha})]$.

Consider any $K$-nonnegative matrix $\tilde{A}$ that satisfies
$\Phi(\tilde{A}) = \Phi(A_{\alpha})$. Since $\Phi(\tilde{A}x_m) =
\Phi(A_{\alpha}x_m)$ and $\Phi(A_{\alpha}x_m)$ is the
$2$-dimensional face of $K$ generated by $x_1,x_2$, after
normalizing $\tilde{A}$, we may assume that $\tilde{A}x_m =
x_1+\tilde{\alpha}x_2$ for some $\tilde{\alpha} > 0$.  Similarly,
we may assume that $\tilde{A}x_i = \tilde{\alpha}_{i+1}x_{i+1}$
for $i = 1, \ldots, m-1$. Substituting the values of the
$\tilde{A}x_i$'s into the relation obtained from (\ref{rel2}) by
applying $\tilde{A}$, we obtain
$$\tilde{\alpha}_2x_2+\tilde{\alpha}_4x_4+ \cdots +
\tilde{\alpha}_{m-2}x_{m-2}+\tilde{\alpha}_mx_m =
\tilde{\alpha}_3x_3+ \tilde{\alpha}_5x_5+ \cdots +
\tilde{\alpha}_{m-1}x_{m-1}+x_1+\tilde{\alpha}x_2.$$ Since
(\ref{rel2}) is, up to multiples, the only relation for the
extreme vectors of $K$, we have
$$\tilde{\alpha}_i = 1 = \tilde{\alpha}_2-\tilde{\alpha}
\mbox{ for } i = 3, \ldots, m.$$  Hence $\tilde{A}$ is given by
$$\tilde{A}x_1 = (1+\tilde{\alpha})x_2, \tilde{A}x_i = x_{i+1} \mbox{ for } i
= 2, \ldots, m-1, \mbox{ and } \tilde{A}x_m =
x_1+\tilde{\alpha}x_2,$$ for some $\tilde{\alpha} > 0$.  This
proves that, after normalization, $\tilde{A}$ equals some
$A_{\tilde{\alpha}}$, which is our contention.

(ii)  Let $A$ be an exp-maximal $K$-primitive matrix.  In view of
Theorem \ref{theorem2}(II)(ii), the digraph $({\cal E},{\cal
P}(A,K))$ is given by Figure $1$, except that $x_1, \ldots, x_m$
are to be replaced by $x_{\sigma(1)}, \ldots, x_{\sigma(m)}$
respectively, where $\sigma$ is some permutation on $\{1, \ldots,
m\}$. By Lemma \ref{lemma 5.1}(i) we can find positive scalars
$\alpha$ and $\lambda_j, j = 1, \ldots, m$, such that, after
normalizing $A$, the relation on $\fn{Ext}K$ and the matrix $A$
are given respectively by the relations obtained from (\ref{rel2})
and (\ref{eq2}) by replacing each $x_j$ by
$\lambda_jx_{\sigma(j)}$. Since the relation on $x_1, \ldots, x_m$
is, up to multiples, unique, it follows that $\lambda_1, \ldots,
\lambda_m$ are the same and moreover $\sigma$ maps the set $\{
1,3, \ldots, n\}$ onto itself or onto $\{ 2,4, \ldots, n+1\}$. It
is readily checked that we have $A =
P_{\sigma}^{-1}A_{\alpha}P_{\sigma}$.   So, after normalization,
$A$ is equivalent to some $A_{\alpha}$.

(iii) As can be readily seen, if $\alpha_1,\alpha_2$ are distinct
positive scalars, then $A_{\alpha_1}$ and $A_{\alpha_2}$ are
linearly independent. So by Lemma \ref{lemma8.1} a positive
multiple of $A_{\alpha_1}$ cannot be cone-equivalent to a positive
multiple of $A_{\alpha_2}$. \hfill$\blacksquare$

\begin{theorem}\label{theorem8.3} Let $K$ be an $n$-dimensional indecomposable exp-maximal minimal
cone, where $n$ is even.  Suppose that the extreme vectors $x_1,
\ldots, x_{n+1}$ of $K$ satisfy relation {\rm(\ref{rel3})}
{\rm(}with $m = n+1${\rm)}.  For every $\alpha
> 0$, let $A_{\alpha}$ be the exp-maximal $K$-primitive matrix
given by {\rm(\ref{eq3})} {\rm(}but with $A$ replaced by
$A_{\alpha}${\rm)}. For every $\alpha, \beta > 0$, let
$A_{\alpha,\beta}$ be the $K$-nonnegative matrix defined by {\rm:}
\begin{eqnarray*}
   A_{\alpha,\beta}x_1 &=& \beta x_2, \\
   A_{\alpha,\beta}x_3 &=& (1+\alpha)x_1+(1+\beta)x_2,  \\
    A_{\alpha,\beta}x_n &=& x_3+\alpha x_1, \\
   A_{\alpha,\beta}x_i &=& \left\{ \begin{array}{cl}x_{i+3}&\mbox{ when }i
   \mbox{ is even}, i \ne n\\x_{i-1}&\mbox{ when }i \mbox{ is odd}, i \ne
   1,3.\end{array}\right.
\end{eqnarray*}
 Then {\rm:}
\begin{enumerate}
\item[{\rm(i)}]  $\Phi(A_{\alpha})$ is a $2$-dimensional face,
independent of the choice of the positive scalar $\alpha${\rm;}
its relative interior consists of positive multiples of matrices
of the form $A_{\tilde{\alpha}}$.
 \item[{\rm(ii)}] $\Phi(A_{\alpha,\beta})$ is
a $3$-dimensional simplicial face, independent of the choice of
the positive scalars $\alpha, \beta${\rm;} its relative interior
consists of positive multiples of matrices of the form
$A_{\tilde{\alpha},\tilde{\beta}}$. \item[{\rm(iii)}] Every
exp-maximal $K$-primitive matrix is cone-equivalent to a positive
multiple of some $A_{\alpha}$ or some $A_{\alpha,\beta}$ and thus
is a positive multiple of a matrix of the form
$P_{\sigma}^{-1}A_{\alpha}P_{\sigma}$ or
$P_{\sigma}^{-1}A_{\alpha,\beta}P_{\sigma}$, where $P_{\sigma}$ is
the automorphism of $K$ given by $P_{\sigma}x_j = x_{\sigma(j)}$
for $j = 1, \ldots, n+1$, $\sigma$ being a permutation on $\{1,
\ldots, n+1\}$ that maps the set $\{ 1,2,4, \ldots, n-2,n\}$ onto
itself.
 \item[{\rm(iv)}] The $A_{\alpha}$'s, $A_{\alpha,\beta}$'s or their positive multiples are pairwise not
cone-equivalent.
\end{enumerate}
\end{theorem}

{\it Proof} (i) Use the same argument as that for Theorem
\ref{theorem8.2}(i).

 (ii) Use the same kind of argument as that for Theorem
\ref{theorem8.2}(i); in this case we can show that
$\Phi(A_{\alpha,\beta})$ is the $3$-dimensional simplicial face
generated by the extreme matrices $B,C,D$ given respectively by:
$$Bx_1 = 0, Bx_3 = x_1+x_2, Bx_n = x_3, \mbox{ and }$$
$$
Bx_i = \left\{\begin{array}{cl}x_{i+3}&\mbox{ when }i \mbox{ is
even}, i \ne n\\x_{i-1}&\mbox{ when }i \mbox{ is odd}, i \ne
1,3\end{array}\right.
$$
$$Cx_3 = x_1 = Cx_n, Cx_i = 0 \mbox{ for } i \ne 3,n$$
and
$$Dx_1 = x_2 = Dx_3, Dx_i = 0 \mbox{ for } i \ne 1,3.$$

(iii)  We first show that $A_{\alpha,\beta}$ is exp-maximal
$K$-primitive.  Let $\tilde{K}$ denote the $n$-dimensional
indecomposable minimal cone with extreme vectors $x_1, \ldots,
x_{n+1}$ that satisfy relation (\ref{rel5}) (with $m = n+1$), and
let $\tilde{A}$ be the $\tilde{K}$-nonnegative matrix defined by
(\ref{eq5}) (but with $A$ replaced by $\tilde{A}$).  By Lemma
\ref{lemma?}(ii) and Theorem \ref{theorem2}(I), $\tilde{A}$ is
exp-maximal $\tilde{K}$-primitive.
 Let $\pi$ be the permutation on $\{1,2, \ldots, n+1\}$
given by $\pi(1) = 3, \pi(2) = 1$, and
$$\pi(j) = \left\{\begin{array}{cl}j-1&\mbox{ when } j \mbox{ is odd}, j\ne
1\\j+1&\mbox{ when } j\mbox{ is even}, j \ne
2.\end{array}\right.$$ Let $P$ be the matrix from $\fn{span}K$ to
$\fn{span}\tilde{K}$ given by $Px_j = x_{\pi(j)}$.  Then, as can
be readily checked, $P$ is an isomorphism which maps $K$ onto
$\tilde{K}$ and, moreover, we have $A_{\alpha,\beta} =
P_{\pi}^{-1}AP_{\pi}$.  So $A_{\alpha,\beta}$ is cone-equivalent
to $\tilde{A}$ and hence is exp-maximal $K$-primitive.

Let $A$ be an exp-maximal $K$-primitive matrix.  In view of
Theorem \ref{theorem2}(III)(ii), the digraph $({\cal E},{\cal
P}(A,K))$ is given by Figure $1$ or Figure $2$ (with $m = n+1$),
except that $x_1, \ldots, x_{n+1}$ are to be replaced by
$x_{\sigma(1)}, \ldots, x_{\sigma(n+1)}$ respectively, where
$\sigma$ is some permutation on $\{1, \ldots, n+1\}$. In the
former case, following the argument given in the proof for Theorem
\ref{theorem8.2}(ii), we can show that $A$ is a positive multiple
of a matrix of the form $P_{\sigma}^{-1}A_{\alpha}P_{\sigma}$
where $P_{\sigma}$ is the automorphism of $K$ given by
$P_{\sigma}x_j = x_{\sigma(j)}$ for $j = 1, \ldots, n+1$, $\sigma$
being a permutation on $\{1, \ldots, n+1\}$ that maps the set $\{
1,2,4, \ldots, n-2,n\}$ onto itself, noting that $\sigma$ cannot
interchange the sets $\{1,2,4, \ldots, m-3,m-1\}$ and $\{3,5,
\ldots, m-2,m\}$ as their cardinality differ by $1$.  So, in this
case, $A$ is cone-equivalent to a positive multiple of some
$A_{\alpha}$.

When the digraph $({\cal E},{\cal P}(A,K))$ is given by Figure $2$
(but with vertices relabelled as indicated above), by Lemma
\ref{lemma 5.1}(ii) we can find positive scalars $\alpha, \beta$
and $\lambda_j, j = 1, \ldots, m$, such that, after normalizing
$A$, the relation on $\fn{Ext}K$ and the matrix $A$ are given by
the relations obtained from (\ref{rel5}) and (\ref{eq5})
respectively by replacing each $x_j$ by $\lambda_jx_{\sigma(j)}$.
But the relation on $\fn{Ext}K$, which is also given by
(\ref{rel3}), is unique (up to multiples), it follows that all the
$\lambda_j$'s are the same.  So the relation on $\fn{Ext}K$ and
the matrix $A$ are given respectively by (\ref{rel5}) and
(\ref{eq5}), but with each $x_j$ replaced by $x_{\sigma(j)}$.  But
then $A$ is cone-equivalent to the matrix $\tilde{A}$, which was
introduced at the beginning of the proof for this part, and hence
is also cone-equivalent to $A_{\alpha,\beta}$, as desired.

(iv)  As done in the proof for Theorem \ref{theorem8.3}(iii),
 if $\alpha_1,\alpha_2$ are distinct
positive scalars, then a positive multiple of $A_{\alpha_1}$ and a
positive multiple of $A_{\alpha_2}$ are linearly independent, and
so by Lemma \ref{lemma8.1} they are not cone-equivalent.  For a
similar reason, a positive multiple of $A_{\alpha_1,\beta_1}$ is
also not cone-equivalent to a positive multiple of
$A_{\alpha_2,\beta_2}$, provided that $(\alpha_1,\beta_1)\ne
(\alpha_2,\beta_2)$.  Moreover, a matrix of the form $A_{\alpha}$
and one of the form $A_{\beta,\gamma}$, or their positive
multiples, are also not cone-equivalent, because $({\cal E},{\cal
P}(A_{\alpha}))$ is given by Figure $1$ whereas $({\cal E},{\cal
P}(A_{\beta,\gamma}))$ is given by Figure $2$, and the two
digraphs are not isomorphic. \hfill$\blacksquare$

Now we consider the exp-maximal primitive matrices for a
decomposable exp-maximal minimal cone.  In this case, Lemma
\ref{lemma8.1} no longer applies.  What we have is the following:

\begin{lemma}\label{lemma8.4} Let $K\in {\cal P}(n+1,n)$ be an exp-maximal decomposable minimal cone with extreme vectors $x_1, \ldots,
x_{n+1}$ {\rm(}where $n$ is even{\rm)}. Suppose that $K =
\fn{pos}\{ x_2 \}\oplus \fn{pos}\{
 x_1,x_3, x_4, \ldots, x_{n+1}\}$, where $x_1,x_3, x_4, \ldots,
 x_{n+1}$
 satisfy the relation given by {\rm(\ref{rel1})} {\rm(}with $m = n+1${\rm)}.  Let $A$ and
$\tilde{A}$ be the $K$-nonnegative matrices defined respectively
 by {\rm(\ref{eq1})} and by the relation obtained from {\rm(\ref{eq1})} by replacing $A, \alpha, \beta$
 by $\tilde{A}, \tilde{\alpha}, \tilde{\beta}$ respectively.  Then for any $\omega > 0$,
$\tilde{A}$ and $\omega A$ are cone-equivalent if and only if
$\omega = 1$ and $\alpha\beta = \tilde{\alpha}\tilde{\beta}$.
\end{lemma}

{\it Proof}.  ``Only if" part:  First, note that the given
assumptions guarantee that the digraphs $({\cal E},{\cal P}(A,K))$
and $({\cal E},{\cal P}(\tilde{A},K))$ are both given by Figure
$2$ (see Lemma \ref{lemma?}(iii)).

Suppose that $\tilde{A}$ and $\omega A$ are cone-equivalent.  Let
$P$ be an automorphism of $K$ such that $P\tilde{A} = \omega AP$.
By the argument given in the proof of Lemma \ref{lemma8.1}, we can
show that $P$ takes $x_2$ to a positive multiple of itself. Then
from the relation $P\tilde{A}x_2 = \omega APx_2$ we infer that $P$
also maps $x_3$ to a positive multiple of itself. Proceeding
inductively, we can show that $P$ maps each $x_j$ to a positive
multiple of itself.  Say, we have $Px_i = \lambda_ix_i$ for $i =
1, \ldots, m$. Substituting the values of the $Px_i$'s into the
relation obtained from (\ref{rel1}) by applying $P$ and using the
fact that, up to multiples, (\ref{rel1}) is the only relation for
the extreme vectors of $K$, we conclude that all the
$\lambda_j$'s, for $j = 1, \ldots, m, j\ne 2$, are equal.  Denote
their common value by $\lambda$.  Then $P$ is given by $Px_i =
\lambda x_i$ for $i = 1, \ldots, m, i\ne 2$, and $Px_2 = \mu x_2$,
where $\mu$ denotes $\lambda_2$.  Now by the given assumptions on
$A$ and $\tilde{A}$ we have
$$P\tilde{A}x_1 = P(\tilde{\alpha}x_2+x_3) = \tilde{\alpha}\mu
x_2+\lambda x_3 \mbox{ and } \omega APx_1 = \omega A(\lambda x_1)
= \omega \lambda(\alpha x_2+x_3).$$  But $P\tilde{A}x_1 = \omega
APx_1$, so we obtain $\omega = 1$ and $\tilde{\alpha}/\alpha =
\lambda/\mu$. Then the relation $P\tilde{A} = \omega AP$ reduces
to $P\tilde{A} = AP$. Similarly, from the relation $P\tilde{A}x_2
= APx_2$ we obtain $\beta/\tilde{\beta} = \lambda/\mu$. Hence we
have $\alpha \beta = \tilde{\alpha}\tilde{\beta}$.

Conversely, suppose that $\alpha \beta =
\tilde{\alpha}\tilde{\beta}$.  Choose positive scalars $\lambda,
\mu$ such that $\lambda/\mu = \tilde{\alpha}/\alpha (=
\beta/\tilde{\beta})$. Let $P$ be the automorphism of $K$
determined by $Px_2= \mu x_2$ and $Px_i = \lambda x_i$ for all $i
\ne 2$.  Then, as can be readily checked, $P\tilde{A}x_i = APx_i$
for every $i$.  Hence we have $P\tilde{A} = AP$, i.e., $A$ and
$\tilde{A}$ are cone-equivalent. \hfill$\blacksquare$

In view of Lemma \ref{lemma8.4} and using the kind of argument as
given in the proofs for Theorem \ref{theorem8.2} and Theorem
\ref{theorem8.3}, we can establish the following, whose proof we
omit:

\begin{theorem}\label{theorem8.5} Let $K\in {\cal P}(n+1,n)$ be an exp-maximal decomposable minimal cone with extreme vectors $x_1, \ldots,
x_{n+1}$ {\rm(}where $n$ is even{\rm)}. Suppose that $K =
\fn{pos}\{ x_2 \}\oplus \fn{pos}\{
 x_1,x_3, x_4, \ldots, x_{n+1}\}$, where $x_1,x_3, x_4, \ldots,
 x_{n+1}$
 satisfy the relation given by {\rm(\ref{rel1})} {\rm(}with $m = n+1${\rm)}.  For every $\alpha,\beta > 0$, let
$A_{\alpha,\beta}$ be the exp-maximal $K$-primitive matrix defined
by {\rm(\ref{eq1})} {\rm(}but with $A$ replaced by
$A_{\alpha,\beta}${\rm)}, and assume that $x_1, \ldots, x_m$
satisfy relation {\rm(\ref{rel1})}. Then {\rm:}
\begin{enumerate}
\item[{\rm(i)}] $\Phi(A_{\alpha,\beta})$ is a $3$-dimensional
simplicial face, independent of the choice of the positive scalars
$\alpha,\beta${\rm;} its relative interior consists of positive
multiples of matrices of the form
$A_{\tilde{\alpha},\tilde{\beta}}$. \item[{\rm(ii)}] Every
exp-maximal $K$-primitive matrix is cone-equivalent to a positive
multiple of some $A_{1,\beta}$. \item[{\rm(iii)}] For distinct
positive scalars $\beta_1, \beta_2$, the matrices $A_{1,\beta_1},
A_{1,\beta_2}$, or their positive multiples, are pairwise not
cone-equivalent.
\end{enumerate}
\end{theorem}

In view of the preceding theorems, we can conclude that for every
exp-maximal minimal cone $K$, indecomposable or not, there are
uncountably infinitely many exp-maximal $K$-primitive matrices
which are pairwise linearly independent and non-cone-equivalent.

Next, we consider the $3$-dimensional cone case.  We will need the
following known result (\cite [Lemma 5.3]{G--K--T}):

\begin{thm} Let $\{ x_1, \ldots, x_n \}$ be a basis for ${\Bbb R}^n$.  Let
$x_0 = \sum_{i=1}^n\alpha_ix_i$ where each $\alpha_i$ is different
from $0$.  Let $A$ and $B$ are $n\times n$ real matrices. Suppose
that $A$ is nonsingular and also that $Bx_j$ is a multiple of
$Ax_j$ for $j = 0, 1, \ldots, n$.  Then $B$ is a multiple of $A$.
\end{thm}

\begin{theorem}\label{theorem8.6}  Let $m \ge 4$ be a positive integer.  For each
$\theta\in (\frac{2\pi}{m},\frac{2\pi}{m-1})$, let $K_\theta$ and
$A_\theta$ be respectively the exp-maximal cone and exp-maximal
$K_\theta$-primitive matrix as defined in Theorem
\ref{coro6}{\rm(iii)}.
\begin{enumerate}
  \item[{\rm(i)}] If
$K\in {\cal P}(m,3)$ is an exp-maximal polyhedral cone and $A$ is
an exp-maximal $K$-primitive matrix with $\rho(A) = 1$, then there
exists a unique $\theta\in (\frac{2\pi}{m},\frac{2\pi}{m-1})$ such
that $A$ is cone-equivalent to $A_\theta$. \item[{\rm(ii)}] When
$m \ge 6, A_\theta$ is, up to positive scalar multiples, the only
exp-maximal $K_\theta$-primitive matrix.
\end{enumerate}
\end{theorem}

{\it Proof}. (i) Since $A$ is exp-maximal, by Theorem
\ref{theorem1}(i), after relabelling its vertices, we may assume
that the digraph $({\cal E}, {\cal P}(A,K))$ is given by Figure
$1$. Let $v\in \fn{int}K^*$ be a Perron vector of $A^T$ and denote
by $C$ the complete cross-section $\{ x\in K: \langle x,v \rangle
= 1\}$ of $K$, which is a polygon with $m$ extreme points.  By the
proof of Lemma \ref{lemma4.2}(i) and in view of Lemma
\ref{lemma6.1}, we may assume that $x_1, \ldots, x_m$ are
precisely all the extreme points of $C$ and $x_i,x_{i+1}$ are
neighborly extreme points for $i = 1, \ldots, m$ (where $x_{m+1}$
is taken to be $x_1$); also, $Ax_j = x_{j+1}$ for $j = 1, \ldots,
m-1$ and $Ax_m = (1-c)x_1+cx_2$ for some $c\in (0,1)$, and
moreover $t^m-ct-(1-c)$ is an annihilating polynomial for $A$. Let
$u$ be the Perron vector of $A$ that belongs to $C$ and denote by
$\hat{C}$ the polygon $C-u$.  Note that $\fn{span}\hat{C}$ equals
$(\fn{span}\{v\})^\bot$ and so it is invariant under $A$ (as $v$
is an eigenvector of $A^T$).  Let $\lambda_1, \lambda_2$ denote
the eigenvalues of the restriction of $A$ to
$(\fn{span}\{v\})^\bot$. It is clear that the eigenvalues of $A$
are $1$ (the Perron root) and $\lambda_1,\lambda_2$.  By the
Perron-Frobenius theory, $|\lambda_j| < 1$ for $j = 1,2$.  We
contend that $\lambda_1, \lambda_2$ form a conjugate pair of
non-real complex numbers.

For $j = 1, \ldots, m$, denote by $y_j$ the point $x_j-u$.
Clearly, $y_1, \ldots, y_m$ are all the extreme points of
$\hat{C}$ and $y_i,y_{i+1}$ are neighborly extreme points for $i =
1, \ldots, m$ (where $y_{m+1}$ is taken to be $y_1$); and $0\in
\fn{ri}\hat{C}$ as $u\in \fn{ri}C$.
 Since $Au = u$,
the action of $A$ on $C$ induces a corresponding action on
$\hat{C}$: we have $Ay_j = y_{j+1} \mbox{ for }j= 1, \ldots, m-1,
\mbox{ and } Ay_m = (1-c)y_1+cy_2$. Take note of the following
consequences of the fact that the digraph $({\cal E}, {\cal
P}(A,K))$ is given by Figure $1$: $A$ maps no extreme points of
$C$ into $\fn{ri}C$ and maps the relative interior of precisely
one $2$-dimensional face of $C$ into $\fn{ri}C$. The preceding
assertion is still true if $C$ is replaced by $\hat{C}$.

Since $t^m-c t-(1-c)$ is an annihilating polynomial for $A$, we
have $\lambda_j^m = c\lambda_j+(1-c)\cdot 1$. Suppose that
$\lambda_1, \lambda_2$ are real.  Then clearly we have $\lambda_1,
\lambda_2\in (-1,0)$.  Let $w \in (\fn{span}\{ v \})^\bot$ be an
eigenvector of $A$ corresponding to $\lambda_1$ and suppose that
$\fn{span} \{ w \}$ meets the relative boundary of $\hat{C}$ at
the points $z_1,z_2$. Then $z_1, z_2$ are scalar multiples of $w$
but with opposite signs; say, we have $z_1 = a_1w$ and $z_2 =
a_2w$ with $|a_2| \ge |a_1|$.  A little calculation shows that
$Az_1 = \alpha z_2$ for some $\alpha\in (0,1)$; as $0\in
\fn{ri}\hat{C}$, this implies that $Az_1\in \fn{ri}\hat{C}$. If
$z_1$ is an extreme point of $\hat{C}$, we already obtain a
contradiction, as $A$ sends no extreme point of $\hat{C}$ to
$\fn{ri}\hat{C}$.  So suppose that $z_1$ lies in the relative
interior of a side of the polygon $\hat{C}$.  Now the point
$Az_2$, which is a positive multiple of $z_1$, either lies in
$\fn{int}\hat{C}$ or is equal to $z_1$. Suppose $Az_2\in
\fn{int}\hat{C}$. Since $A$ sends no extreme point of $\hat{C}$ to
$\fn{int}\hat{C}$, it follows that $z_2$ is not an extreme point
and so it lies in the relative interior of a side of the polygon
$\hat{C}$. But then $A$ maps the relative interior of two
different sides of $\hat{C}$ into its relative interior, which is
a contradiction. So we must have $Az_2 = z_1$. If $z_2$ is an
extreme point, then necessarily $z_2 = y_m$ and the side of
$\hat{C}$ that contains $z_1$ is the line segment
$\overline{y_1y_2}$. On the other hand, since $Ay_1 = y_2$ and
$Ay_2 = y_3$, we obtain $Az_1\in \fn{ri}\overline{y_2y_3}$, which
contradicts the fact that $Az_1\in \fn{ri}\hat{C}$.  So $z_2$ must
lie in the relative interior of a side of $\hat{C}$.  Then
necessarily the side of $\hat{C}$ that contains $z_2$ is
$\overline{y_{m-2}y_{m-1}}$, whereas the side that contains $z_1$
is $\overline{y_{m-1}y_m}$. Since $z_2\in
\fn{ri}\overline{y_{m-2}y_{m-1}}$, we have $z_1 = Az_2\in
\fn{ri}\overline{y_{m-1}y_m}$ and hence $Az_1\in
\fn{ri}\fn{conv}\{ y_m,y_1,y_2 \}$.  Note that $Az_1$ also belongs
to $\fn{ri}\fn{conv} \{ y_{m-2},y_{m-1},y_m \}$ as it lies in
$\fn{ri}\overline{z_1z_2}$ and $z_1\in
\fn{ri}\overline{y_{m-1}y_m}, z_2\in
\fn{ri}\overline{y_{m-2}y_{m-1}}$.
 But $\fn{ri}\fn{conv}\{ y_m,y_1,y_2 \}\cap
\fn{ri}\fn{conv} \{ y_{m-2},y_{m-1},y_m \} = \emptyset$ as $m \ge
4$, so we arrive at a contradiction.

In the above we have shown that the eigenvalues
$\lambda_1,\lambda_2$ of $A$ form a conjugate pair of non-real
complex numbers, say, they are $\pm re^{i\theta}$.  Then clearly
$r < 1$. We can now choose a basis $\{ u_1, u_2 \}$ for
$(\fn{span}\{ v \})^\bot$ with $u_1 = y_1$ such that
\begin{eqnarray*}
  Au_1 &=& r\cos \theta\ u_1+ r\sin \theta\ u_2, \\
  Au_2 &=& -r\sin \theta\ u_1+r\cos \theta\ u_2.
\end{eqnarray*}
For $j = 1, \ldots, m$, let $y_j = \alpha_ju_1+\beta_ju_2$.  Then
we have
$${\alpha_{j+1} \choose \beta_{j+1}} =
r\left (\begin{array}{rr}\cos \theta & -\sin \theta \\ \sin
\theta& \cos \theta \end{array} \right ){\alpha_j \choose
\beta_j}$$ for $j = 1, \ldots, m-1$, as $Ay_j = y_{j+1}$ for every
such $j$.  Since $y_1, \ldots, y_m$ form the consecutive vertices
of a polygon in $(\fn{span}\{ v \})^\bot$, the points
$(\alpha_j,\beta_j)^T, j = 1, \ldots, m$, also form the
consecutive vertices of a polygon in ${\Bbb R}^2$.  Now it should
be clear that we have $(m-1)\theta < 2\pi$ and $2\pi < m\theta$,
which implies that $\theta\in (\frac{2\pi}{m},\frac{2\pi}{m-1})$.

As $(\alpha_1,\beta_1)^T = (1,0)^T$, a little calculation gives
$Ay_m = r^m \cos m\theta u_1+ r^m\sin m\theta u_2$.  Since $u_1,
u_2$ are linearly independent, the relation $Ay_m = (1-c)y_1+cy_2$
leads to relations (\ref{okay1}) and (\ref{okay2}) that appear in
the proof of Lemma \ref{lemma6.2}.  Eliminating $c$ from these two
relations, we obtain
$$ \frac{\sin (m-1)\theta}{\sin \theta}r^m - \frac{\sin
m\theta}{\sin \theta}r^{m-1}+1 = 0;$$ hence $r$ equals
$r_{\theta}$, the unique positive real root of the polynomial
$g_\theta(t)$.  It is clear that $\theta$ is unique, because we
have $c = r_{\theta}^{m-1}\frac{\sin m\theta}{\sin \theta}$ (see
Corollary \ref{coro6.3}(ii)).  Now let $P$ be the $3\times 3$
matrix given by: $Pu_j = e_j$, where $u_3 = u$, the Perron vector
of $A$, and $e_j$ is the $j$th standard unit vector of ${\Bbb
R}^3$.  It is readily checked that $P$ is a nonsingular matrix
that maps $K$ onto $K_\theta$. Moreover, we have $PA =
A_{\theta}P$.  So the cones $K$ and $K_\theta$ are linearly
isomorphic, and the cone-preserving maps $A$ and $A_\theta$ are
cone-equivalent.

(ii) Let $B$ be an exp-maximal $K_\theta$-primitive matrix.  Then
$\gamma(B) = 2m+1$ and by Theorem \ref{theorem1}(i) the digraph
$({\cal E},{\cal P}(B,K_\theta))$ is, apart from the labelling of
its vertices, given by Figure $1$.  For simplicity, we denote
$x_j(\theta)$ by $x_j$ for $j = 1, \ldots, m$.  Also, we adopt the
convention that for any integer $j\notin \{1, \ldots, m\}, x_j$ is
taken to be $x_k$, where $k$ is the unique integer that satisfies
$1\le k \le m, k \equiv j (\fn{mod} m)$.  According to Lemma
\ref{lemma6.1}, adjacent vertices of the digraph $({\cal E},{\cal
P}(B,K_\theta))$ correspond to neighboring extreme rays of
$K_\theta$.  Using an argument similar to the one given in the
proof of Theorem \ref{theorem6.4}, we can show that there exists
$p, 1\le p \le m$ such that one of the following holds:

(I) $Bx_j$ is a positive multiple of $x_{j+1}$ for $j = p, p+1,
\ldots, p+m-2$ and $Bx_{p+m-1}$ is a positive linear combination
of $x_{p+m}$ and $x_{p+m+1}$; or

(II) $Bx_j$ is a positive multiple of $x_{j-1}$ for $j = p, p-1,
\ldots, p-m+2$ and $Bx_{p-m+1}$ is a positive linear combination
of $x_{p-m}$ and $x_{p-m-1}$.

We consider the case when (I) holds first.   Then $Bx_j$ is a
positive multiple of $A_\theta x_j$ for all $j = 1, \ldots, m$
except for $j\equiv 1,p+m-1 (\fn{mod}m)$. Since there are $m-2$
such $x_j$'s and $m-2 \ge 4$ (as $m \ge 6$), by Theorem E, $B$
must be a positive multiple of $A_\theta$, which is what we want.
(In fact, then necessarily $p = 1$.)

When (II) holds, we find that $Bx_j$ is a positive multiple of
$A_\theta^{-1} x_j$ for all $j = 1, \ldots, m$, except for $j
\equiv 1, p-m+1 (\fn{mod}m)$.  By Theorem E again we conclude that
$B$ is a positive multiple of $A_\theta^{-1}$, which is
impossible, as $A_\theta^{-1}$ is not $K_\theta$-nonnegative.
\hfill$\blacksquare$\\

\begin{coro}\label{coro8.9} Let $m\ge 6$ be a positive integer, and let $K_\theta\in {\cal P}(m,3)$
be the exp-maximal cone as defined before. Then{\rm:}
\begin{enumerate}
\item[{\rm(i)}] The automorphism group of $K_\theta$ consists of
scalar matrices only. \item[{\rm(ii)}] For any $\theta_1,
\theta_2\in (\frac{2\pi}{m},\frac{2\pi}{m-1}), \theta_1\ne
\theta_2$, the cones $K_{\theta_1}, K_{\theta_2}$ are not linearly
isomorphic.
\end{enumerate}
\end{coro}

{\it Proof}.  (i) We first establish the following:

{\bf Assertion}. Every automorphism of $K_{\theta}$ that commutes
with $A_{\theta}$ is a scalar matrix.

{\it Proof}.  Let $P$ be an automorphism of $K_{\theta}$ that
commutes with $A$.  For simplicity, denote $x_j(\theta)$ by $x_j$
for $j = 1, \ldots, m$.  Since $P$ is an automorphism of
$K_{\theta}, P$ permutes the extreme rays of $K_{\theta}$ among
themselves.  So $Px_1$ is a positive multiple of $x_p$ for some
$p\in \langle m \rangle$.
 In view of the relation $A_{\theta}Px_1 =
PA_{\theta}x_1$, we see that $Px_2$ is a positive multiple of
$A_{\theta}x_p$.  But $A_{\theta}x_p$ is a positive multiple of
$x_{p+1}$ if $p < m$ and is a positive linear combination of $x_1$
and $x_2$ if $p = m$, so we must have $1 \le p < m$ and $Px_2$ is
a positive multiple of $x_{p+1}$.  Suppose that $2 \le p$.  By
considering the relation $A_{\theta}Px = PA_{\theta}x$ for $x =
x_2, \ldots, m-p+1$, and in this order, and proceeding
inductively, we find that $Px_j = x_{p+j-1}$ for $j = 1, \ldots,
m+1-p$, and in particular we have $Px_{m-p+1} = x_m$.  Then from
the relation $A_{\theta}Px_{m-p+1} = PA_{\theta}x_{m-p+1}$ we
infer that $Px_{m+2-p}$ is a positive linear combination of $x_1$
and $x_2$, which is a contradiction.  Thus, we have $p = 1$ and by
the same kind of argument we can then show that $Px_j$ is a
positive multiple of $x_j$ for $j = 1, \ldots, m$.  Hence $P\in
\Phi(I)$.  Then we can invoke \cite[Theorem 3.3]{L--S} to conclude
that $\Phi(I)$ is an extreme ray of $\Phi(I)$, and so $P$ is a
scalar matrix. \\

Now back to the proof of (i).  Let $P$ be an automorphism of
$K_{\theta}$. Since $A_\theta$ is an exp-maximal
$K_\theta$-primitive matrix, so is $P^{-1}A_{\theta}P$.  By
Theorem \ref{theorem8.6}(ii) we have $P^{-1}A_{\theta}P = \alpha
A_{\theta}$ for some $\alpha > 0$.  As $P^{-1}A_{\theta}P$ and
$A_{\theta}$ are similar and $A_{\theta}$ is nonsingular,
necessarily $\alpha = 1$.  So we have $A_{\theta}P = PA_{\theta}$,
and by the above Assertion it follows that $P$ is a scalar matrix.

(ii) Let $\theta_1, \theta_2\in (\frac{2\pi}{m},\frac{2\pi}{m-1})$
be such that the cones $K_{\theta_1}, K_{\theta_2}$ are linearly
isomorphic, say, $P$ is a linear isomorphism that maps
$K_{\theta_2}$ onto $K_{\theta_1}$.  Since $A_{\theta_1}$ is
exp-maximal $K_{\theta_1}$-primitive, clearly
$P^{-1}A_{\theta_1}P$, which is cone-equivalent to $A_{\theta_1}$,
is exp-maximal $K_{\theta_2}$-primitive.  In view of Theorem
\ref{theorem8.6}(ii) and the fact $\rho(A_{\theta_j}) = 1$ for $j
= 1,2$, we have $P^{-1}A_{\theta_1}P = A_{\theta_2}$.  Now for $j
= 1,2$, the eigenvalues of $A_{\theta_j}$ are $1$ and
$r_{\theta_j}e^{\pm i\theta_j}$.  So we must have $\theta_1 =
\theta_2$. \hfill$\blacksquare$\\

Note that part (i), and hence also part (ii), of Theorem
\ref{theorem8.6} is not true for $m = 3$. This is because every
$A_{\theta}$ has a pair of conjugate non-real complex eigenvalues,
whereas an exp-maximal ${\Bbb R}^3_+$-primitive matrix need not
have non-real eigenvalues (see Remark \ref{remark8.8} below).

As can be readily checked, parts (i) and (ii) of Corollary
\ref{coro8.9} are also both invalid when $m = 3 \mbox{ or }4$.

We {\it do not} know whether Theorem \ref{theorem8.6}(ii) or
Corollary \ref{coro8.9}(i),(ii) can be extended to cover the case
$m = 5$. However, we are going to show that the two problems are
equivalent.

Upon close examination, one finds that the Assertion given in the
proof of Corollary \ref{coro8.9}, in fact, holds also for every
positive integer $m \ge 4$. Furthermore, if Theorem
\ref{theorem8.6}(ii) is true for $m = 5$, then the arguments for
part (i) and (ii) of the corollary still work for the case $m =
5$.

Now suppose that Corollary \ref{coro8.9} is extendable to the case
$m = 5$, but Theorem \ref{theorem8.6}(ii) is not.  Then for some
$\theta\in (\frac{2\pi}{5},2\frac{\pi}{4})$ there exists an
exp-maximal $K_{\theta}$-primitive matrix $A$, different from
$A_{\theta}$, such that $\rho(A) = 1$.  By Theorem
\ref{theorem8.6}(i), there exists $\phi\in
(\frac{2\pi}{5},\frac{2\pi}{4})$ such that $K_{\theta}$ is
cone-equivalent to $K_{\phi}$ and $A$ is cone-equivalent to
$K_{\phi}$.  So there exists an isomorphism $P$ such that
$PK_{\phi} = K_{\theta}$ and $A_{\phi} = P^{-1}AP$.  In view of
Corollary \ref{coro8.9}(ii), necessarily $\phi = \theta$.  Hence,
$P$ is an automorphism of $K_{\theta}$ and by Corollary
\ref{coro8.9}(i), $P$ is a scalar matrix.  Therefore, $A =
A_{\phi} = A_{\theta}$, which is a contradiction.

As expected, we have the following result, which describes all the
exp-maximal $K$-primitive matrices for $K$ in ${\cal P}(n,n), n
\ge 3$:

\begin{remark}\label{theorem8.7}\rm
Let $K\in {\cal P}(n,n), n \ge 3$ and let $A$ be a $K$-primitive
matrix with $\rho(A) = 1$.  Then $A$ is exp-maximal $K$-primitive
if and only if there exists $c\in (0,1)$ such that $A$ is
cone-equivalent to $C_h{\rm(}\in \pi({\Bbb R}^n_+){\rm)}$, the
companion matrix of the polynomial $h(t) = t^n-ct-(1-c)$, i.e.,
$$C_h =
\left[\begin{array}{ccccc}
0 & & & & 1-c \\
1 & 0 & \abb{\Large \textbf{0}}{10,1} & & c \\
& 1 & \ddots & & 0 \\
& \abb{\Large \textbf{0}}{-5,-12} & \ddots & \ddots & \vdots \\
& & & 1 & 0
\end{array}\right]_{n\times n}
.$$
\end{remark}

{\it Proof}.  The ``if" part is obvious as $C_h$ is exp-maximal
${\Bbb R}^n_+$-primitive.  To show the ``only if" part, we may
assume that $K = {\Bbb R}^n_+$.  Then there exists a permutation
matrix $P$ such that $P^{-1}AP = B$, where $B$ is a nonnegative
matrix with the same zero-nonzero pattern as the companion matrix
$C_h$. It is not difficult to find a diagonal matrix $D$ with
positive diagonal entries such that $D^{-1}BD = C_h$ for some
$c\in (0,1)$.  But $PD$ is an automorphism of ${\Bbb R}^n_+$, so
it follows that $A$ is cone-equivalent to $C_h$.
\hfill$\blacksquare$\\

For completeness, let us mention the following result, which is
not difficult to prove:

\begin{remark}\label{remark8.8}\rm  Let $K\in {\cal P}(3,3)$ and let $A$ be an
a $K$-primitive matrix with $\rho(A) = 1$.  Then $A$ is
exp-maximal $K$-primitive if and only if $A$ is cone-equivalent to
one of the following:
\begin{enumerate}
\item[(i)] $A_{\theta}\in \pi(K_{\theta})$, where $\theta\in
(\frac{2\pi}{3},\pi)$, and $K_{\theta}\in {\cal P}(3,3)$ and
$A_{\theta}\in \pi(K_{\theta})$ are defined in the same way as
before; \item[(ii)] $\fn{diag}(\alpha_1,\alpha_2,1)\in
\pi(\tilde{K})$, where for some $c\in (\frac{3}{4},1), \alpha_1,
\alpha_2$ are the (distinct) real roots, other than $1$, of the
polynomial $t^3-ct-(1-c)$ and $\tilde{K}$ is the polyhedral cone
in ${\Bbb R}^3$ generated by the extreme vectors $x_1 = (1,1,1)^T,
x_2 = (\alpha_1,\alpha_2,1)^T$ and $x_3 =
(\alpha_1^2,\alpha_2^2,1)^T$; \item[(iii)]
$\left[\begin{array}{rrr}-\frac{1}{2}&1&0\\0&-\frac{1}{2}&0\\0&0&1\end{array}\right]\in
\pi(K_0)$, where $K_0 = \fn{pos}\{ x_1,Ax_1,A^2x_1\}$ with $x_1 =
(1,1,1)^T.$
\end{enumerate}
\end{remark}

\setcounter{equation}{0}
\section{Remarks and open questions}

A positive integer $\kappa$ is called the {\it critical exponent}
of a normed space $E$ (or of the norm on $E$) if the inequalities
$\|A^k\| = \|A\| = 1$ imply that $\rho(A) = 1$, and if $\kappa$ is
the smallest number with the indicated property.  It is known that
not every norm in a finite-dimensional space has a critical
exponent.  An example of one such norm can be found in \cite
{B--L}.  Borrowing the latter example, we are going to show that
there exists a proper cone which does not have finite exponent.

\begin{example}\rm
 Let $\|\cdot\|$ denote the norm of ${\Bbb R}^2$ whose unit closed ball
 is defined by the inequalities:
 $$ 3\xi_1-2 \le \xi_2 \le \xi_1^3, \mbox{ if } -2\le \xi_1 \le
 -1,$$
$$3\xi_1-2 \le \xi_2 \le 3\xi_1+2, \mbox{ if } |\xi_1| \le 1,$$
and
$$ \xi_1^3 \le \xi_2 \le 3\xi_1+2, \mbox{ if } 1 \le \xi_1 \le
2.$$  (See Figure 6.)  Let $K$ be the proper cone in $\IR^3$ given
by: $K=\{\alpha {x \choose 1}: \alpha\ge 0$ and $\|x\| \le 1\}$.
For every positive integer k, let $B_k$ denote the $2\times 2$
diagonal matrix $\fn{diag}(2^{-1/k},2^{-3/k})$.
 As shown in \cite [p.67]{B--L}, $B_k$ has the property that $\|B_k\|=\|B^k_k\|=1$ but
$\|B^{k+1}_k\| <1$. Let $A_k =B_k \oplus (1)$. Then it is easy to
see that $A_k$ is $K$-primitive and $\gamma(A_k)=k+1$. Since $k$
can be arbitrarily large, this shows that for this $K$ we have
$\gamma(K)=\infty$. It is also of interest to note that the
$K$-primitive matrices $A_k$ obtained in this example are, in
fact, all extreme matrices of the cone $\pi(K)$. The point is, $K$
is an indecomposable cone and each of the $A_k$$^\prime$s maps
infinitely many extreme rays of $K$ onto extreme rays.
\end{example}

\[
{\linethickness{0.7pt} \qbezier(-10,0)(10,11)(12,56)
\qbezier(-10,0)(-30,-11)(-32,-56) } \put(-65,0){\vector(1,0){125}}
\put(-10,-63){\vector(0,1){125}} \put(12,56){\line(-1,-2){37}}
\put(-32,-56){\line(1,2){37}} \put(65,-4){$\xi_1$}
\put(-13,68){$\xi_2$} \put(15,17){$\xi_2=\xi_1^3$}
\put(-25,-85){\rm Figure 6.}
\]

Let $E_n$ denote the set of values attained by the exponents of
primitive matrices of order $n$. Dulmage and Mendelsohn
\cite{D--M} have found intervals in the set $\{1,2, \ldots,
(n-1)^2+1 \}$ containing no integer which is the exponent of a
primitive matrix of order $n$.  These intervals have been called
{\it gaps} in $E_n$.  The problem of determining $E_n$ or the gaps
is an intricate problem, but it has been completely resolved.
(See, for instance, \cite{B--R}.)

For a given polyhedral cone (or a proper cone) $K$, we can
consider a similar problem --- to determine the set of values
attained by the exponents of $K$-primitive matrices.  We expect
that for every polyhedral cone $K$ of dimension greater than $2$
there are gaps in the set of values attained by the exponents of
$K$-primitive matrices (but at present we {\it do not} have a
proof for this claim). As an illustration, consider an
indecomposable minimal cone $K\in {\cal P}(m,n)$ with a balanced
relation for its extreme vectors, where $n$ is an odd integer $\ge
5$.  For a $K$-primitive matrix $A$, if the digraph $({\cal
E},{\cal P})$ is given by Figure $1$ or Figure $2$ then
$\gamma(A)$ equals $n^2-n+1$ or $n^2-n$ (see Theorem
\ref{theorem2}(II)) and Lemma \ref{lemma?}(ii)).  On the other
hand, if the digraph is not given by Figure $1$ or Figure $2$,
then by Lemma \ref{lemma3} the length of the shortest circuit in
$({\cal E},{\cal P})$ is at most $n-1(= m-2)$ and by Remark
\ref{remark extra} it follows that $\gamma(A) \le (n-1)^2+2$. So
in this case any integer lying in the closed interval
$[n^2-2n+4,n^2-n-1]$ cannot be attained as the exponent of some
$K$-primitive matrix.\\

Perhaps, a less difficult problem is the following:\\

{\bf Question 9.1.} Let $m \ge 4$ be a positive integer. Determine
the set of integers that can be attained as the exponent of a
$K$-primitive matrix for some $n$-dimensional polyhedral cone $K$
with $m$ extreme rays, where $3\le n\le m$.\\

Let $K\in {\cal P}(m,n)$ be an exp-maximal non-simplicial
polyhedral cone and let $A$ be an exp-maximal $K$-primitive
matrix.  According to Theorem \ref{maintheorem}, $\gamma(A)$
equals $(n-1)(m-1)+1$ or $(n-1)(m-1)$. In view of Theorem
\ref{theorem1}(i), in either case we have $n = m_A$.  If
$\gamma(A) = (n-1)(m-1)+1$, then by Theorem \ref{theorem1}(i)
again the digraph $({\cal E},{\cal P})$ is given by Figure $1$. If
$\gamma(A) = (n-1)(m-1)$, then by Theorem \ref{theorem1}(ii)
either $({\cal E},{\cal P})$ is given by Figure $1$ or Figure $2$,
or $m_A = 3$.  The last possibility cannot happen, because then we
have $\gamma(A) = 2(m-1)$ and $n = m_A = 3$, in contradiction with
Theorem \ref{coro6}(i).  If $({\cal E},{\cal P})$ is given by
Figure $1$, then by Lemma \ref{lemma4.2}(iii) $K$ is
indecomposable.  If $({\cal E},{\cal P})$ is given by Figure $2$,
then again by Lemma \ref{lemma4.2}(iii), $K$ is either
indecomposable or is an even-dimensional minimal cone which is the
direct sum of a ray and an indecomposable minimal cone with a
balanced relation for its extreme vectors.  Conversely, if $K$ is
an even-dimensional minimal cone with the said property, then by
Theorem \ref{theorem2}(III)(i) $K$ is exp-maximal.  We have thus
characterized all decomposable non-simplicial exp-maximal
polyhedral cones.\\

{\bf Question 9.2.} Identify indecomposable exp-maximal cones in
${\cal P}(m,n)$ for $m > n$, $m \ne n+1$ and $n \ne 3$.\\

In this work we are able to identify the exp-maximal cones and the
corresponding exp-maximal primitive matrices only for the extreme
cases $m = n, m = n+1$ and $n = 3$.  To deal with the other cases,
one may consider the following question first:\\

{\bf Question 9.3.}  Given positive integers $m,n$ with $3\le n
\le m$, characterize the $n\times n$ real matrices $A$ with the
property that there exists $K\in {\cal P}(m,n)$ such that $A$ is
$K$-nonnegative and $({\cal E},{\cal P}(A,K))$ is given by Figure
$1$.\\

In below we provide the answers to the preceding question for two
special cases, namely, $n = 3, m \ge 4$, and $m = n$. For
convenience, we normalize the matrices under consideration.

\begin{remark}\label{remark9.1}\rm  Let $A$ be a $3\times 3$ real matrix with $\rho(A) =
1$, and let $m\ge 4$ be a given positive integer.  Then there
exists $K\in {\cal P}(m,3)$ such that $A$ is $K$-nonnegative and
$({\cal E}, {\cal P}(A,K))$ is given by Figure $1$ (and hence $A$
is exp-maximal $K$-primitive) if and only if the eigenvalues of
$A$ are $1, r_{\theta}e^{\pm i\theta}$, where $\theta\in
(\frac{2\pi}{m},\frac{2\pi}{m-1})$, $r_\theta$ is the unique
positive real root of the polynomial $g_{\theta}(t)$ given in
Lemma \ref{lemma6.2}.
\end{remark}

{\it Proof}.  ``Only if" part:  Since $({\cal E},{\cal P}(A,K))$
is isomorphic to Figure $1$, by Lemma \ref{lemma6.1}(ii) and
Theorem \ref{coro6}(i), $A$ is exp-maximal $K$-primitive.  By
Theorem \ref{theorem8.6}(i) $A$ is cone-equivalent to $A_{\theta}$
for some $\theta\in (\frac{2\pi}{m},\frac{2\pi}{m-1})$.  Hence,
the eigenvalues of $A$ are the same as those for $A_{\theta}$,
i.e., $1, r_{\theta}e^{\pm i\theta}$.

``If" part:  Since $A$ and $A_{\theta}$ are similar, both being
similar to the diagonal matrix
$\fn{diag}(1,r_{\theta}e^{i\theta},r_{\theta}e^{-i\theta})$, there
exists a $3\times 3$ real nonsingular matrix $P$ such that
$P^{-1}AP = A_{\theta}$.  Take $K = PK_{\theta}$.  As can be
readily checked, $A$ is $K$-nonnegative and $A$ is cone-equivalent
to $A_{\theta}$.  Hence $({\cal E}, {\cal P}(A,K))$ is isomorphic
to Figure $1$. \hfill$\blacksquare$

\begin{remark}\label{remark9.2}\rm  Let $A$ be an $n\times n$ real matrix with
$\rho(A) = 1$.  Then there exists $K\in {\cal P}(n,n)$ such that
$A$ is $K$-nonnegative and the digraph $({\cal E},{\cal P}(A,K))$
is given by Figure $1$ (and hence $A$ is exp-maximal
$K$-primitive) if and only if $A$ is non-derogatory and the
characteristic polynomial of $A$ is of the form $t^n-ct-(1-c)$,
where $c\in (0,1)$.
\end{remark}

{\it Proof}. ``Only if" part:  Suppose that there exists $K\in
{\cal P}(n,n)$ such that $A$ is $K$-nonnegative and the digraph
$({\cal E},{\cal P}(A,K))$ is given by Figure $1$.  By Lemma
\ref{lemma4.2} $A$ is non-derogatory and has an annihilating
polynomial of the form of the form $t^n-ct-(1-c)$, where $c\in
(0,1)$. Since $A$ is $n\times n$, clearly $t^n-ct-(1-c)$ is also
the characteristic polynomial of $A$.

 ``If" part:  In this case, $A$ is similar to $C_h$, the companion matrix of
$h(t)$.  So there exists an $n\times n$ real matrix $P$ such that
$P^{-1}AP = C_h$.  Let $K = P{\Bbb R}^n_+$.  Then $A$ is
$K$-nonnegative and is cone-equivalent to $C_h$.  By Theorem
\ref{theorem8.7} $A$ is exp-maximal $K$-primitive.  But the cone
$K$ is simplicial, so the digraph $({\cal E}, {\cal P}(A,K))$ is
isomorphic to Figure $1$. \hfill$\blacksquare$

We would like to add that in the ``if" part of Remark
\ref{remark9.2}, except for the case $c = c_n, n$ being odd, we
may omit the assumption that $A$ be non-derogatory. The point is,
then by Lemma \ref{lemma6.2}(i) $A$ has simple eigenvalues and
hence is, necessarily, non-derogatory.  As an illustration for the
exceptional case, consider the matrix $A =
\fn{diag}(-\frac{1}{2},-\frac{1}{2},1)$.  Its characteristic
polynomial is $t^3-\frac{3}{4}t-\frac{1}{4}$, which is of the said
form.  (Recall that $c_3 = \frac{3}{4}$.) Since $A$ is not
non-derogatory, by Lemma \ref{lemma4.2}(i), there cannot exist
$K\in {\cal P}(n,n)$ such that $A$ is $K$-nonnegative and the
digraph $({\cal E},{\cal P}(A,K))$ is given by Figure $1$.

Note that in Question 9.3 we do not pose the same question for
Figure $2$, because the two questions are equivalent; that is, the
existence of a pair $(K,A)$ with $({\cal E},{\cal P}(A,K))$ given
by Figure $1$ guarantees the existence of a pair $(K,A)$ with
$({\cal E},{\cal P}(A,K))$ given by Figure $2$, and vice versa.

To see this, suppose that $A$ is $K$-nonnegative and $({\cal E},
{\cal P}(A,K))$ is given by Figure $2$.  Then $\Phi(x_1+x_2)$ is a
$2$-dimensional face of $K$ and $Ax_m$ lies in its relative
interior. Let $\hat{K}$ denote the polyhedral cone generated by
$Ax_m, x_2, x_3, \ldots, x_m$. It is readily shown that $Ax_m,
x_2, \ldots, x_m$ are precisely (up to multiples) all the extreme
vectors of $\hat{K}$.  Furthermore, $A$ is $\hat{K}$-nonnegative
and $({\cal E},{\cal P}(A,\hat{K}))$ is isomorphic to Figure $1$
(under the isomorphism given by: $\Phi(Ax_m)\mapsto \Phi(x_m),
\Phi(x_j)\mapsto \Phi(x_{j-1})$ for $j = 2, \ldots, m$).
Conversely, suppose that $A$ is $K$-nonnegative and $({\cal E},
{\cal P}(A,K))$ is given by Figure $1$.  Let $\tilde{K} =
\fn{pos}\{ (1-\alpha)x_1+\alpha x_m, x_1, x_2, \ldots, x_{m-1}\}$.
 It is not difficult to show that for $\alpha > 1$,
sufficiently close to $1$, $(1-\alpha)x_1+\alpha x_m, x_1, x_2,
\ldots, x_{m-1}$ are precisely all the extreme vectors of
$\tilde{K}$. Furthermore, $A$ is $\tilde{K}$-nonnegative and
$({\cal E},{\cal P}(A,\tilde{K}))$ is isomorphic to Figure $2$
(under the isomorphism given by: $\Phi((1-\alpha)x_1+\alpha
x_m)\mapsto \Phi(x_1), \Phi(x_j)\mapsto \Phi(x_{j+1})$ for $j = 1.
\ldots, m-1$).

Now it should be clear that Lemma \ref{theorem7.1} is still true
if Figure $1$ is replaced by Figure $2$.

Below are two other questions on this topic that one may
explore:\\

{\bf Question 9.4.}  If $K$ is an $n$-dimensional minimal cone
such that the relation for its extreme vectors has $p$ vectors on
one side and $q$ vectors on the other side, where $p,q \ge 2, p+q
\le
n+1$, what is $\gamma(K)$ ?\\

{\bf Question 9.5.}  Given positive integers $m, n$ with $m \ge
n$, determine
$$\fn{min}\{ \gamma(K): K\in {\cal P}(m,n) \}.$$

We suspect that for the polyhedral cone $K_0\in {\cal P}(m,3)$
with extreme vectors $x_1, \ldots, x_m$ given by $x_j = (\cos
\frac{2j\pi}{m}, \sin \frac{2j\pi}{m},1)^T$ for $j = 1, \ldots,
m$, we have $$\gamma(K_0) = \fn{min}\{ \gamma(K): K\in {\cal
P}(m,3) \}.$$\\

{\bf Acknowledgement} Part of this work was carried out in the
period from 2004, November to 2005, March when the second author
visited the Department of Mathematics at the Technion as a Lady
Davis visiting professor.  He would like to take this opportunity
to thank the Technion for its hospitality.\\

\newpage

\end{document}